\documentclass[a4paper,twoside]{article}
\usepackage{amssymb}
\usepackage{amsbsy}
\usepackage{amsmath}
\usepackage{theorem}
\usepackage{epsfig}
\usepackage{graphicx}

\numberwithin{equation}{section}
\newtheorem{theorem}{Theorem}[section]
\newtheorem{corollary}[theorem]{Corollary}
\newtheorem{conjecture}[theorem]{Conjecture}
\newtheorem{lemma}[theorem]{Lemma}
\newtheorem{proposition}[theorem]{Proposition}

{\theorembodyfont{\rmfamily}
\newtheorem{definition}[theorem]{Definition}
}
{\theorembodyfont{\rmfamily}
\newtheorem{remark}[theorem]{Remark}
}
{\theorembodyfont{\rmfamily}

}
{\theorembodyfont{\rmfamily}
\newtheorem{example}[theorem]{Example}
}
{\theorembodyfont{\rmfamily}
\newtheorem{examples}[theorem]{Examples}
}
{\theorembodyfont{\rmfamily}
\newtheorem{final remarks and questions}[theorem]{Final remarks and questions}
}
\begin{document}

\title{\textbf{On the String-Theoretic Euler Number }\\
\textbf{of 3-dimensional A-D-E Singularities}}
\author{{\thanks{\noindent{}{\scriptsize {GSRT-fellow supported by the E.U. and the
Greek Ministry of Research and Technology.}\newline
{2000 MS-Classification. Primary 14Q15, 32S35, 32S45; Secondary 14B05,
14E15, 32S05, 32S25.}}}\ }\textbf{Dimitrios I. Dais}\ \ and\ \ $^{\diamond}$%
\textbf{Marko Roczen\medskip} \\
\noindent$^{\ast}${\scriptsize {Mathematics Department, Section of Algebra
and Geometry,}}\\
{\scriptsize {University of Ioannina, GR-45110 Ioannina, Greece}}\\
{\scriptsize {e-mail: ddais@cc.uoi.gr}}\\
$^{\diamond}$ {\scriptsize {Institut f\"{u}r Mathematik, Humboldt
Universit\"{a}t zu Berlin,}}\\
{\scriptsize {Rudower Chaussee 25, D-10099 Berlin, Germany}}\\
{\scriptsize {e-mail: roczen@mathematik.hu-berlin.de}}}
\date{}
\maketitle

\begin{abstract}
\noindent The string-theoretic E-functions $E_{\text{str}}\left(
X;u,v\right) $ of normal complex varieties $X$ having at most log-terminal
singularities are defined by means of snc-resolutions. We give a direct
computation of them in the case in which $X$ is the underlying space of the
three-dimensional $\mathbf{A}$-$\mathbf{D}$-$\mathbf{E}$ singularities by
making use of a canonical resolution process. Moreover, we compute the
string-theoretic Euler number for several compact complex threefolds with
prescribed $\mathbf{A}$-$\mathbf{D}$-$\mathbf{E}$ singularities.
\end{abstract}

%
%
%

%
%
%

\section{{\ }Introduction\label{INTRO}}

\noindent The \textit{string-theoretic} (or \textit{stringy}) \textit{Hodge
numbers} $h_{\text{str}}^{p,q}\left( X\right) $ of normal, projective
complex varieties $X$ with at most Gorenstein quotient or toroidal
singularities were introduced in \cite{BD} in an attempt to determine a
suitable mathematical formulation (and generalization) for the numbers which
are encoded into the Poincar\'{e} polynomial of the chiral and antichiral
rings of the physical ``integer charge orbifold theory'', due to the
LG/CY-correspondence of Vafa, Witten, Zaslow and others. (See \cite{VAFA},
\cite[\S3-5]{WITTEN}, \cite[\S4]{ZASLOW}). These numbers are generated by
the so-called $E_{\text{str}}$-\textit{polynomials} and, as it was shown in
\cite{BD} and \cite{BB}, they are the right quantities to establish several
\textit{mirror-symmetry identities }for Calabi-Yau varieties. In fact, as
long as a stratification (separating singularity types) for such an $X$ is
available, the key-point is how one defines the $E_{\text{str}}$-polynomial
\textit{locally} at these special Gorenstein singular points (by
``measuring'', in a sense, how far they are from admitting of crepant
resolutions).\medskip

\noindent Recently Batyrev \cite{BATYREV1} generalized this definition and
made it work also for the case in which one allows $X$ to have at most
log-terminal singularities. In this general framework, ones has to introduce
appropriate $E_{\text{str}}$-\textit{functions} $E_{\text{str}}\left(
X;u,v\right) $ instead which may be not even rational. The treatment of
varieties $X$ with $e_{\text{str}}\left( X\right) =$ lim$_{u,v\rightarrow
1}\,E_{\text{str}}\left( X;u,v\right) \notin\mathbb{Z}$ is therefore
unavoidable. Nevertheless, as it turned out, this new language is a very
important tool as it unifies the considerations of certain invariants
associated to a wide palette of ``MMP-singularities'' and leads to the use
of more flexible manipulations, as for example in the study of the behaviour
of \textit{log-flips}, and in the proof of cohomological \textit{Mckay
correspondence - }both on the level of counting dimensions and on the level
of determining the motivic Gorenstein volume. (See \cite[1.6, 4.11 and 8.4]
{BATYREV2} and \cite[Thm. 5.1]{DENEF-LOESER}).\medskip

\noindent In the present paper we deal with the evaluation of the $E_{\text{%
str}}$-functions and string-theoretic Euler numbers for the
three-dimensional $\mathbf{A}$-$\mathbf{D}$-$\mathbf{E}$ singularities, and
emphasize some distinctive features of the computational methodology.\bigskip

\noindent\textsf{(a) }\textbf{Log-terminal singularities. }Let $X$ be a
normal complex variety, i.e., a normal, integral, separated scheme of finite
type over $\mathbb{C}$. Suppose that $X$ is $\mathbb{Q}$-Gorenstein, i.e.,
that a positive integer multiple of its canonical Weil divisor $K_{X}$ is a
Cartier divisor. $X$ is said to have at most \textit{log-terminal}
(respectively,\textit{\ canonical} / \textit{terminal}) singularities if
there exists an \textit{snc}-desingularization $\varphi:\widetilde{X}%
\longrightarrow X$, i.e., a desingularization of $X$ whose exceptional locus
$\frak{Ex}\left( \varphi\right) =\cup_{i=1}^{r}D_{i}$ consists of \textbf{s}%
mooth prime divisors $D_{1},D_{2},\ldots,D_{r}$ with only \textbf{n}ormal
\textbf{c}rossings, such that the ``discrepancy'' w.r.t. $\varphi$, which is
the difference between the canonical divisor of $\widetilde{X}$ and the
pull-back of the canonical divisor of $X$, is of the form
\begin{equation*}
K_{\widetilde{X}}-\varphi^{\ast}\left( K_{X}\right) =\sum_{i=1}^{r}\ a_{i}\
D_{i}
\end{equation*}
with all the $a_{i}$'s $>-1$ ( $\geq0$ / $>0$).

\begin{examples}
(i) The quotients $\mathbb{C}^{2}/G$, $\ $for $G$ a linearly acting finite
subgroup of GL$\left( 2,\mathbb{C}\right) $ (resp. of SL$\left( 2,\mathbb{C}%
\right) )$, have at most log-terminal (resp. canonical) isolated
singularities.\smallskip \newline
(ii) All $\mathbb{Q}$-Gorenstein toric varieties have at most log-terminal
(but not necessarily isolated) singularities.\smallskip
\end{examples}

\noindent\textsf{(b)} $E$\textbf{-polynomials. }As it was shown by Deligne
in \cite[\S8]{DELIGNE}, the cohomology groups $H^{i}\left( X,\mathbb{Q}%
\right) $ of any complex variety $X$ are equipped with a functorial \textit{%
mixed Hodge structure} (MHS).\textit{\ }The same remains true if one works
with cohomologies $H_{c}^{i}\left( X,\mathbb{Q}\right) $ with compact
supports. There exist namely an increasing \textit{weight-filtration}
\begin{equation*}
\mathcal{W}_{\bullet}:\ \ 0=W_{-1}\subset W_{0}\subset W_{1}\subset
\cdots\subset W_{2i-1}\subset W_{2i}=H_{c}^{i}\left( X,\mathbb{Q}\right)
\end{equation*}
and a decreasing \textit{Hodge-filtration }
\begin{equation*}
\mathcal{F}^{\,\bullet}:\ \ H_{c}^{i}\left( X,\mathbb{C}\right)
=F^{0}\supset F^{1}\supset\cdots\supset F^{i}\supset F^{i+1}=0,
\end{equation*}
such that $\mathcal{F}^{\,\bullet}$ induces a natural filtration
\begin{equation*}
\begin{array}{l}
F^{p}\left( Gr_{k}^{\mathcal{W}_{\bullet}}(H_{c}^{i}\left( X,\mathbb{C}%
\right) )\right) = \\
(W_{k}\left( H_{c}^{i}\left( X,\mathbb{C}\right) \right) \cap F^{p}\left(
H_{c}^{i}\left( X,\mathbb{C}\right) \right) +W_{k-1}\left( H_{c}^{i}\left( X,%
\mathbb{C}\right) \right) )\ /\ W_{k-1}\left( H_{c}^{i}\left( X,\mathbb{C}%
\right) \right)
\end{array}
\end{equation*}
(denoted again by $\mathcal{F}^{\,\bullet}$) on the complexification of the
graded pieces
\begin{equation*}
Gr_{k}^{\mathcal{W}_{\bullet}}(H_{c}^{i}\left( X,\mathbb{Q}\right)
)=W_{k}/W_{k-1}.
\end{equation*}
Let now
\begin{equation*}
h^{p,q}\left( H_{c}^{i}\left( X,\mathbb{C}\right) \right) :=\text{ dim}_{%
\mathbb{C}}Gr_{\mathcal{F}^{\bullet}}^{p}Gr_{p+q}^{\mathcal{W}_{\bullet
}}\left( H_{c}^{i}\left( X,\mathbb{C}\right) \right)
\end{equation*}
denote the corresponding \textit{Hodge numbers} by means of which one
defines the so-called $E$\textit{-polynomial} of $X$:
\begin{equation*}
E\left( X;u,v\right) :=\sum_{p,q}\ e^{p,q}\left( X\right) \ u^{p}v^{q}\in%
\mathbb{Z}\left[ u,v\right] ,
\end{equation*}
where
\begin{equation*}
e^{p,q}\left( X\right) :=\sum_{i\geq0}\ \left( -1\right) ^{i}\ h^{p,q}\left(
H_{c}^{i}\left( X,\mathbb{C}\right) \right) \ .
\end{equation*}
The $E$-polynomials are to be viewed as ``generating functions'' encoding
our invariants. For instance, the topological Euler characteristic $e\left(
X\right) $ is $E\left( X;1,1\right) $. In fact, the $E$-polynomial behaves
similarly; e.g, for locally closed subvarieties $Y,Y_{1},Y_{2}$ of $X$,
\begin{equation}
E\left( X\mathbb{r}Y;u,v\right) =E\left( X;u,v\right) -E\left( Y;u,v\right) ,
\label{E-DIFF}
\end{equation}
\begin{equation}
E\left( Y_{1}\cup Y_{2};u,v\right) =E\left( Y_{1};u,v\right) +E\left(
Y_{2};u,v\right) -E\left( Y_{1}\cap Y_{2};u,v\right)  \label{INC-EX}
\end{equation}
and
\begin{equation}
E\left( X;u,v\right) =E\left( F;u,v\right) \cdot E\left( Z;u,v\right)
\label{E-FIB}
\end{equation}
whenever $F$ denotes the fiber of a Zariski locally trivial fibration $%
X\longrightarrow Z$.

\begin{example}
If $Y\longrightarrow X$ is the blow-up of a $d$-dimensional complex manifold
$X$ at a point $x\in X$ and $D\cong \mathbb{P}_{\mathbb{C}}^{d-1}$ the
exceptional divisor, then $E\left( Y;u,v\right) $ equals
\begin{equation}
E\left( X\mathbb{r}\left\{ x\right\} ;u,v\right) +E\left( D;u,v\right)
=E\left( X;u,v\right) +uv+\left( uv\right) ^{2}+\cdots +\left( uv\right)
^{d-1}\ \smallskip   \label{BLOWUP}
\end{equation}
\end{example}

\noindent\textsf{(c)} $E_{\text{str}}$\textbf{-functions.} Allowing the
existence of log-terminal singularities in order to pass to stringy
invariants, one takes essentialy into account the ``discrepancy
coefficients''.

\begin{definition}
\label{EST-DEF}Let $X$ be a normal complex variety with at most log-terminal
singularities, $\varphi :\widetilde{X}\longrightarrow X$ \ an
snc-desingularization of $X$ as in \textsf{(a)}, $D_{1},D_{2},\ldots ,D_{r}$
the prime divisors of the exceptional locus, and $I:=\left\{ 1,2,\ldots
,r\right\} $. For any subset $J\subseteq I$ define
\begin{equation*}
D_{J}:=\left\{
\begin{array}{ll}
\widetilde{X}, & \text{if \ }J=\varnothing  \\
\, & \, \\
\bigcap_{j\in J}\ D_{j}, & \text{if \ }J\neq \varnothing
\end{array}
\ \ \ \ \ \ \ \text{and }\right. \ \ \ \ \ \ \ D_{J}^{\circ }:=D_{J}\,%
\mathbb{r}\bigcup_{j\in I\mathbb{r}J}\ D_{j}\ .
\end{equation*}
The algebraic function
\begin{equation}
E_{\text{str}}\left( X;u,v\right) :=\sum_{J\subseteq I}\ E\left(
D_{J}^{\circ };u,v\right) \ \prod_{j\in J}\ \frac{uv-1}{\left( uv\right)
^{a_{j}+1}-1}  \label{E-STR}
\end{equation}
(under the convention for $\prod_{j\in J}$ to be $1$, if $J=\varnothing ,$
and $E\left( \varnothing ;u,v\right) :=0$) is called the \textit{%
string-theoretic} $E$\textit{-function }of $X$.
\end{definition}

\noindent The main result of \cite{BATYREV1} says that:

\begin{theorem}
\label{INDEP}The string-theoretic $E$-function $E_{\text{\emph{str}}}\left(
X;u,v\right) $ is independent of the choice of the snc-desingularization $%
\varphi :\widetilde{X}\longrightarrow X.$
\end{theorem}

\begin{remark}
\label{REM1}(i) The proof of \ \ref{INDEP} relies on ideas of Kontsevich
\cite{KONTSEVICH}, Denef and Loeser by making use of the interpretation of
the defining formula (\ref{E-STR}) as some kind of ``motivic non-Archimedean
integral'' over the space of arcs of $\widetilde{X}$. (For an introduction
to \textit{motivic integration }and\textit{\ measures}, we refer to Craw
\cite{CRAW} and Looijenga \cite{LOOIJENGA}).\smallskip \newline
(ii) To define (\ref{E-STR}) it is sufficient for $\varphi :\widetilde{X}%
\longrightarrow X$ to fulfil the snc-condition only for those $D_{i}$'s for
which $a_{i}\neq 0.$ \smallskip \newline
(iii) If \ $X$ admits a \textit{crepant} desingularization $\pi :\widehat{X}%
\longrightarrow X$, i.e., $K_{\widehat{X}}=\pi ^{\ast }K_{X}$ with $\widehat{%
X}$ smooth, then $E_{\text{str}}\left( X;u,v\right) =E(\widehat{X}%
;u,v).\smallskip $\newline
\noindent (iv) In general $E_{\text{str}}\left( X;u,v\right) $ may be not a
rational function in the two variables $u,v$. Nevertheless, if $X$ has at
most Gorenstein singularities, then the discrepancy coefficients $%
a_{1},\ldots ,a_{r}$ are non-negative integers and
\begin{equation*}
E_{\text{str}}\left( X;u,v\right) \in \mathbb{Z}[\![u,v]\!]\cap \mathbb{Q(}%
u,v).
\end{equation*}
(Of course, for $X$ projective, stringy Hodge numbers $h_{\text{str}%
}^{p,q}\left( X\right) $ can be defined only if $E_{\text{str}}\left(
X;u,v\right) \in \mathbb{Z}\left[ u,v\right] $).\smallskip \newline
(v) The \textit{existence} of snc-desingularizations of any $X$ is
guaranteed by Hironaka's main theorems \cite{HIRONAKA}. But since definition
\ref{EST-DEF} is intrinsic in its nature, it is practically fairly difficult
to compute $E_{\text{str}}\left( X;u,v\right) $ precisely without having
\textit{at least one} snc-desingularization of $X$ at hand, accompanied
firstly with the intersection graph of $D_{1},\ldots ,D_{r}$ and secondly
with the knowledge of their analytic structure.
\end{remark}

\begin{definition}
One defines the rational number
\begin{equation}
e_{\text{str}}\left( X\right) :=\text{ }\underset{u,v\rightarrow 1}{\text{lim%
}}E_{\text{str}}\left( X;u,v\right) =\sum_{J\subseteq I}\ e\left(
D_{J}^{\circ }\right) \ \prod_{j\in J}\ \frac{1}{a_{j}+1}\   \label{e-STR}
\end{equation}
as the \textit{string-theoretic Euler number} of $X$. Moreover, the \textit{%
string-theoretic index }ind$_{\text{str}}\left( X\right) $ of $X$ is defined
to be the positive integer
\begin{equation*}
\text{ind}_{\text{str}}\left( X\right) :=\text{min}\left\{ \ l\in \mathbb{Z}%
_{\geq 1}\mathbb{\ }\left| \ \ e_{\text{str}}\left( X\right) \in \frac{1}{l}%
\,\mathbb{Z\ }\right. \right\} \ .
\end{equation*}
\end{definition}

\begin{examples}
(i) For $\mathbb{Q}$-Gorenstein toric varieties $X$, ind$_{\text{str}}\left(
X\right) =1$, and $e_{\text{str}}\left( X\right) $ is equal to the
normalized volume of the defining fan. Moreover, for Gorenstein toric
varieties $X$, $E_{\text{str}}\left( X;u,v\right) $ is a
polynomial.\smallskip \newline
(ii) Normal algebraic surfaces $X$ with at most log-terminal singularities
have ind$_{\text{str}}\left( X\right) =1$. There exist, however, normal
complex varieties $X$ of dimension $d\geq 3$ with at most Gorenstein
canonical singularities having ind$_{\text{str}}\left( X\right) >1$.
\end{examples}

\noindent\noindent\noindent Batyrev formulated in \cite[5.9]{BATYREV1} the
following conjecture:$\allowbreak$

\begin{conjecture}[On the range of the string-theoretic index]
\label{BAT-CONJ}Let $X$ be a $d$-dimensional normal complex variety having
at most Gorenstein canonical singularities. Then\smallskip\ \emph{ind}$_{%
\text{\emph{str}}}\left( X\right) $ is bounded by a constant $C\left(
d\right) $ depending only on $d$.
\end{conjecture}

\begin{remark}
As it will be clear by Theorem \ref{MAIN}, Conjecture \ref{BAT-CONJ} is
\textit{not} true in general. Nevertheless, there exist several classes of
examples of such $X$'s with string-theoretic index bounded by a constant
which depends exclusively on the dimension. (See e.g. \cite[5.1, 5.10]
{BATYREV1} for the case in which $X$ is the cone over a $(d-1)$-dimensional
smooth projective Fano variety being equipped with a projective embedding
defined by a suitable very ample line bundle)\smallskip . The problem of
chacterizing those $X$'s having bounded ind$_{\text{str}}\left( X\right) $
is still open.
\end{remark}

\noindent \noindent \noindent \textsf{(d) }\textbf{The }$\mathbf{A}$\textbf{-%
}$\mathbf{D}$\textbf{-}$\mathbf{E}$\textbf{'s. }The $d$-dimensional
analogues of the classical hypersurface $\mathbf{A}$\textbf{-}$\mathbf{D}$%
\textbf{-}$\mathbf{E}$ singularities \cite{DUVAL} have underlying spaces of
the form
\begin{equation}
\fbox{$
\begin{array}{ccc}
& \  &  \\
& X_{f}:=X_{f}^{\left( d\right) }:=\text{Spec}\left( \mathbb{C}\left[
x_{1},\ldots ,x_{d+1}\right] \ /\ \text{\textbf{(}}f\text{\textbf{)}}\right)
,\ d\geq 2, &  \\
& \  &  \\
& \text{ with\ \ \ \ }f\left( x_{1},\ldots ,x_{d+1}\right) :=g\left(
x_{1},x_{2}\right) +g^{\prime }(x_{3},\ldots ,x_{d+1}) &  \\
& \  &
\end{array}
$}\medskip  \label{HYP}
\end{equation}
where $g\left( x_{1},x_{2}\right) $ is the defining polynomial of a simple
curve singularity
\begin{equation*}
X_{g}:=\,\text{Spec}\left( \mathbb{C}\left[ x_{1},x_{2}\right] \ /\ \text{%
\textbf{(}}g\text{\textbf{)}}\right)
\end{equation*}
in the affine plane with

\begin{equation*}
\begin{tabular}{|c|c|}
\hline
\textbf{Types} & $g\left( x_{1},x_{2}\right) $ \\ \hline\hline
$\mathbf{A}_{n}$ & $
\begin{array}{c}
\  \\
x_{1}^{n+1}+x_{2}^{2},\ n\geq1 \\
\
\end{array}
$ \\ \hline
$\mathbf{D}_{n}$ & $
\begin{array}{c}
\  \\
x_{1}^{n-1}+x_{1}x_{2}^{2},\ n\geq4 \\
\
\end{array}
$ \\ \hline
$\mathbf{E}_{6}$ & $
\begin{array}{c}
\  \\
x_{1}^{3}+x_{2}^{4} \\
\
\end{array}
$ \\ \hline
$\mathbf{E}_{7}$ & $
\begin{array}{c}
\  \\
x_{1}^{3}+x_{1}x_{2}^{3} \\
\
\end{array}
$ \\ \hline
$\mathbf{E}_{8}$ & $
\begin{array}{c}
\  \\
x_{1}^{3}+x_{2}^{5} \\
\
\end{array}
$ \\ \hline
\end{tabular}
\medskip
\end{equation*}
and $g^{\prime}(x_{3},\ldots,x_{d+1}):=\sum_{j=3}^{d+1}x_{j}^{2}$ is nothing
but the defining quadratic polynomial of the affine $\left( d-2\right) $%
-dimensional quadric
\begin{equation*}
X_{g^{\prime}}:=X_{g^{\prime}}^{(d-2)}:=\text{Spec}\left( \mathbb{C}\left[
x_{3},\ldots,x_{d+1}\right] \ /\ \text{\textbf{(}}g^{\prime}\text{\textbf{)}}%
\right) .
\end{equation*}

\newpage

\begin{remark}
\label{REM2}The $d$-dimensional $\mathbf{A}$\textbf{-}$\mathbf{D}$\textbf{-}$%
\mathbf{E}$ singularities have lots of interesting properties:\smallskip
\newline
\noindent (i) Herszberg \cite{HERSZBERG} and Treger \cite[Thm. 1]{TREGER}
proved that they are \textit{absolutely isolated}, i.e., that they can be
resolved by blowing up successively a finite number of closed points; in
fact, up to analytic isomorphism, they are the only absolutely isolated
singularities of multiplicity $2$.\smallskip \newline
\noindent (ii) Generalizing the classical result of Artin \cite{ARTIN},
Burns \cite[3.3-3.4]{BURNS} showed that they are \textit{rational}, i.e.,
that for any desingularization $\pi :Y\longrightarrow X_{f}^{\left( d\right)
}$ in dimension $d\geq 2$, we have $\left( R^{i}\pi _{\ast }\mathcal{O}%
_{Y}\right) _{\mathbf{0}}=0$ for all $i\geq 1$. In particular, this means
that they have to be canonical (resp. terminal) of index $1$ for $d\geq 2$
(resp. for $d\geq 3$); cf. Reid \cite{REID}.\smallskip \newline
\noindent (iii) Finally, \ Arnold's results \cite{ARNOLD} (see also
\cite[8.26-8.27]{DIMCA1}) imply that they are the only simple (i.e., ``$0$%
-modular'') hypersurface singularities.
\end{remark}

\noindent {}These properties lead us to the conclusion that $X_{f}^{\left(
d\right) }$'s might belong to the class of the best possible candidates for
performing concrete computations for the string-theoretic invariants. On the
other hand, we should stress that none of the above general techniques
mentioned in \ref{REM2} (i)-(ii) are ``constructive'' enough in the sense of
\ref{REM1} (v). That's why we restrict ourselves in this paper to the
three-dimensional case, and based on a canonical snc-resolution being
constructed by Giblin \cite{GIBLIN} and independently by the second-named
author in \cite{ROCZEN1}, \cite{ROCZEN2}, we work out the needed details to
prove the following:

\begin{theorem}
\label{MAIN}The rational, string-theoretic $E$-functions of the underlying
spaces $X=X_{f}^{\left( 3\right) }$ of the $3$-dimensional $\mathbf{A}$-$%
\mathbf{D}$-$\mathbf{E}$-singularities are functions in $w=uv$ given by the
following formulae\emph{:\medskip }{\small \newline
\emph{(i) }$\mathbf{Type}$ $\mathbf{A}_{n}$, $n$ $\mathbf{even.}$
\begin{equation*}
\fbox{$
\begin{array}{l}
E_{\text{\emph{str}}}\left( X;u,v\right) =w^{3}+w-1+\sum\limits_{i=2}^{\frac{%
n}{2}}\frac{\left( w-1\right) \left( w^{2}-1\right) }{w^{i+1}-1}+\frac{%
\left( w-1\right) w^{2}}{w^{n+3}-1} \\
\  \\
+\left( w-1\right) \left( w^{2}-1\right) \left[ \sum\limits_{i=1}^{\frac{n}{2%
}-1}\frac{1}{(w^{i+1}-1)(w^{i+2}-1)}+\frac{1}{(w^{\frac{n}{2}%
+1}-1)(w^{n+3}-1)}\right]
\end{array}
$}
\end{equation*}
\medskip \emph{(ii) }$\mathbf{Type}$ $\mathbf{A}_{n}$, $n$ $\mathbf{odd.}$
\begin{equation*}
\fbox{$
\begin{array}{l}
E_{\text{\emph{str}}}\left( X;u,v\right) =\left( w-1\right) \,(\,w+1)^{2}+w
 + \lfloor \frac{1}{n} \rfloor \\
\  \\
+\left( w^{2}-1\right) \left[ \sum\limits_{i=2}^{\frac{n-1}{2}}\frac{\left(
w-1\right) }{w^{i+1}-1}+\frac{w}{w^{\frac{n+3}{2}}-1}+\sum\limits_{i=1}^{%
\frac{n-1}{2}}\frac{\left( w-1\right) }{(w^{i+1}-1)(w^{i+2}-1)}\right]
\cdot \lceil \frac{n-1}{n} \rceil
\end{array}
$}
\end{equation*}

\emph{(iii) }$\mathbf{Type}$ $\mathbf{D}_{n}$, $n$ $\mathbf{even.}$%
\begin{equation*}
\fbox{$
\begin{array}{l}
E_{\text{\emph{str}}}\left( X;u,v\right) =\left( w-1\right) \,(w^{2}+3w+1)
\\
\  \\
+\left( w-1\right) \,(\,w+1)^{2}\left[ \frac{2}{w^{n}-1}+\sum\limits_{i=3}^{%
\frac{n}{2}+1}\frac{1}{w^{2(n+4-2i)}-1}\right]  \\
\, \\
+2\left( w-1\right) \,(\,1+4w+w^{2})\left[ \sum\limits_{i=1}^{\frac{n}{2}-1}%
\frac{1}{w^{\left( \frac{n}{2}-i+1\right) }-1}\right]  \\
\  \\
+\left( 1+w\right) \left[ 4\left( \frac{w-w^{n}}{w^{n}-1}\right) \left(
\frac{w-w^{\frac{n}{2}}}{w^{\frac{n}{2}}-1}\right) +\sum\limits_{i=1}^{\frac{%
n}{2}-1}\left( \frac{w-w^{\left( \frac{n}{2}-i+1\right) }}{w^{\left( \frac{n%
}{2}-i+1\right) }-1}\right) ^{2}\right]  \\
\, \\
+\left( 1+w\right) \left[ 2\sum\limits_{(\kappa ,\lambda )}\left( \frac{%
w-w^{\kappa +1}}{w^{\kappa +1}-1}\right) \left( \frac{w-w^{\lambda +1}}{%
w^{\lambda +1}-1}\right) -7\left( \tfrac{n}{2}-1\right) \right]  \\
\  \\
+\sum\limits_{(\kappa ,\lambda ,\mu )}\left( \frac{w-w^{\kappa +1}}{%
w^{\kappa +1}-1}\right) \left( \frac{w-w^{\lambda +1}}{w^{\lambda +1}-1}%
\right) \left( \frac{w-w^{\mu +1}}{w^{\mu +1}-1}\right)  \\
\  \\
+2\sum\limits_{(\kappa ^{\prime },\lambda ^{\prime },\mu ^{\prime })}\left(
\frac{w-w^{\kappa ^{\prime }+1}}{w^{\kappa ^{\prime }+1}-1}\right) \left(
\frac{w-w^{\lambda ^{\prime }+1}}{w^{\lambda ^{\prime }+1}-1}\right) \left(
\frac{w-w^{\mu ^{\prime }+1}}{w^{\mu ^{\prime }+1}-1}\right) +2n-5
\end{array}
$}
\end{equation*}
where the pairs $(\kappa ,\lambda )$ of the fourth sum are taken from the
set
\begin{equation*}
\begin{array}{l}
\{\left. (\tfrac{n}{2}-i,\tfrac{n}{2}-\left( i+1\right) )\ \right| \ 1\leq
i\leq \tfrac{n}{2}-2\} \\
\  \\
\cup \{\left. (\tfrac{n}{2}-i,2\left( n-2i\right) -1)\ \right| \ 1\leq i\leq
\tfrac{n}{2}-1\} \\
\  \\
\cup \{\left. (\tfrac{n}{2}-\left( i+1\right) ,2\left( n-2i\right) -1)\
\right| \ 1\leq i\leq \tfrac{n}{2}-2\},
\end{array}
\end{equation*}
the triples $(\kappa ,\lambda ,\mu )$ of the fifth sum from the set
\begin{equation*}
\begin{array}{l}
\{\left. (\tfrac{n}{2}-i,\tfrac{n}{2}-i,2\left( n-2i\right) -1)\ \right| \
1\leq i\leq \tfrac{n}{2}-1\} \\
\  \\
\cup \{\left. (\tfrac{n}{2}-\left( i+1\right) ,\tfrac{n}{2}-\left(
i+1\right) ,2\left( n-2i\right) -1)\ \right| \ 1\leq i\leq \tfrac{n}{2}-2\},
\end{array}
\end{equation*}
and the triples $(\kappa ^{\prime },\lambda ^{\prime },\mu ^{\prime })$ of
the sixth sum from the set
\begin{equation*}
\begin{array}{l}
\{\left. (\tfrac{n}{2}-i,\tfrac{n}{2}-\left( i+1\right) ,2\left( n-2i\right)
-1)\ \right| \ 1\leq i\leq \tfrac{n}{2}-2\} \\
\  \\
\cup \{(n-1,\tfrac{n}{2}-1,\tfrac{n}{2}-1)\}\ .
\end{array}
\end{equation*}

\emph{(iv) }$\mathbf{Type}$ $\mathbf{D}_{n}$, $n$ $\mathbf{odd.}$%
\begin{equation*}
\fbox{$
\begin{array}{l}
E_{\text{\emph{str}}}\left( X;u,v\right) =\left( w-1\right) \,(\,w+1)^{2} \\
\  \\
+\left( w-1\right) \,(\,w+1)^{2}\left[ \frac{1}{w^{n-1}-1}+\frac{1}{w^{n}-1}%
+\sum\limits_{i=3}^{\frac{n+1}{2}}\frac{1}{w^{2(n+3-2i)}-1}\right]  \\
\, \\
+2\left( w-1\right) \,(\,1+4w+w^{2})\left[ \sum\limits_{i=1}^{\frac{n-3}{2}}%
\frac{1}{w^{\left( \frac{n-1}{2}-i+1\right) }-1}\right]  \\
\  \\
+2\left( 1+w\right) \left( \frac{w-w^{\frac{n-1}{2}}}{w^{\frac{n-1}{2}}-1}%
\right) \left[ \frac{w-w^{n}}{w^{n}-1}+\frac{w-w^{n-1}}{w^{n-1}-1}\right]
\\
\  \\
+\left( 1+w\right) \left[ \left( \frac{w-w^{n}}{w^{n}-1}\right) \left( \frac{%
w-w^{n-1}}{w^{n-1}-1}\right) +\sum\limits_{i=1}^{\frac{n-3}{2}}\left( \frac{%
w-w^{\left( \frac{n-1}{2}-i+1\right) }}{w^{\left( \frac{n-1}{2}-i+1\right)
}-1}\right) ^{2}\right]  \\
\  \\
+\left( 1+w\right) \left[ 2\sum\limits_{(\kappa ,\lambda )}\left( \frac{%
w-w^{\kappa +1}}{w^{\kappa +1}-1}\right) \left( \frac{w-w^{\lambda +1}}{%
w^{\lambda +1}-1}\right) -\frac{7}{2}\left( n-1\right) +6\right]  \\
\  \\
+\sum\limits_{(\kappa ,\lambda ,\mu )}\left( \frac{w-w^{\kappa +1}}{%
w^{\kappa +1}-1}\right) \left( \frac{w-w^{\lambda +1}}{w^{\lambda +1}-1}%
\right) \left( \frac{w-w^{\mu +1}}{w^{\mu +1}-1}\right)  \\
\  \\
+2\sum\limits_{(\kappa ^{\prime },\lambda ^{\prime },\mu ^{\prime })}\left(
\frac{w-w^{\kappa ^{\prime }+1}}{w^{\kappa ^{\prime }+1}-1}\right) \left(
\frac{w-w^{\lambda ^{\prime }+1}}{w^{\lambda ^{\prime }+1}-1}\right) \left(
\frac{w-w^{\mu ^{\prime }+1}}{w^{\mu ^{\prime }+1}-1}\right) +2\left(
n-1\right) -4\,
\end{array}
$}
\end{equation*}
where the pairs $(\kappa ,\lambda )$ are taken from the set
\begin{equation*}
\begin{array}{l}
\{\left. (\tfrac{n-1}{2}-i,\tfrac{n-1}{2}-\left( i+1\right) )\ \right| \
1\leq i\leq \tfrac{n-5}{2}\} \\
\  \\
\cup \{\left. (\tfrac{n-1}{2}-i,2\left( n-2i\right) -3)\ \right| \ 1\leq
i\leq \tfrac{n-3}{2}\} \\
\  \\
\cup \{\left. (\tfrac{n-1}{2}-\left( i+1\right) ,2\left( n-2i\right) -3)\
\right| \ 1\leq i\leq \tfrac{n-5}{2}\},
\end{array}
\end{equation*}
the triples $(\kappa ,\lambda ,\mu )$ from the set
\begin{equation*}
\begin{array}{l}
\{(n-1,\tfrac{n-3}{2},\tfrac{n-3}{2})\}\cup \{\left. (\tfrac{n-1}{2}-i,%
\tfrac{n-1}{2}-i,2\left( n-2i\right) -3)\ \right| \ 1\leq i\leq \tfrac{n-3}{2%
}\} \\
\  \\
\cup \{\left. (\tfrac{n-1}{2}-\left( i+1\right) ,\tfrac{n-1}{2}-\left(
i+1\right) ,2\left( n-2i\right) -3)\ \right| \ 1\leq i\leq \tfrac{n-5}{2}\},
\end{array}
\end{equation*}
and the triples $(\kappa ^{\prime },\lambda ^{\prime },\mu ^{\prime })$ from
the set
\begin{equation*}
\begin{array}{l}
\{\left. (\tfrac{n-1}{2}-i,\tfrac{n-1}{2}-\left( i+1\right) ,2\left(
n-2i\right) -3)\ \right| \ 1\leq i\leq \tfrac{n-5}{2}\} \\
\  \\
\cup \{(n-1,n-2,\tfrac{n-3}{2})\}\ .
\end{array}
\end{equation*}

\emph{(v) }$\mathbf{Type}$ $\mathbf{E}_{6}.$%
\begin{equation*}
\fbox{$
\begin{array}{l}
E_{\text{\emph{str}}}\left( X;u,v\right) =w^{3}-1+\tfrac{w+1}{w^{2}+1}+%
\tfrac{\left( w+1\right) ^{2}\left( w-1\right) }{w^{7}-1}+\tfrac{\left(
w+1\right) ^{2}\left( w-1\right) }{w^{10}-1}+\tfrac{2(1+4w+w^{2})}{w+1} \\
\  \\
+\left( 1+w\right) \left[ \sum\limits_{(\kappa ,\lambda )}\left( \frac{%
w-w^{\kappa +1}}{w^{\kappa +1}-1}\right) \left( \frac{w-w^{\lambda +1}}{%
w^{\lambda +1}-1}\right) -9\right]  \\
\  \\
+\sum\limits_{(\kappa ,\lambda ,\mu )}\left( \frac{w-w^{\kappa +1}}{%
w^{\kappa +1}-1}\right) \left( \frac{w-w^{\lambda +1}}{w^{\lambda +1}-1}%
\right) \left( \frac{w-w^{\mu +1}}{w^{\mu +1}-1}\right) +5
\end{array}
$}
\end{equation*}
where the pairs $(\kappa ,\lambda )$ of the first sum are taken from the set
\begin{equation*}
\{(1,1),(1,3),(3,1),(1,6),(6,1),(1,9),(9,1),(3,6),(6,9)\}
\end{equation*}
and the triples $(\kappa ,\lambda ,\mu )$ of the second sum from the set
\begin{equation*}
\{(1,1,9),(1,6,9),(1,9,6),(1,3,6),(1,6,3)\}\ .
\end{equation*}

\emph{(vi) }$\mathbf{Type}$ $\mathbf{E}_{7}.$%
\begin{equation*}
\fbox{$
\begin{array}{l}
E_{\text{\emph{str}}}\left( X;u,v\right) =\left( w-1\right) \,(\,w+1)^{2}%
\left[ 1+\frac{1}{w^{6}-1}+\frac{1}{w^{10}-1}+\frac{1}{w^{12}-1}+\frac{1}{%
w^{14}-1}\right]  \\
\  \\
\allowbreak +2\left( w-1\right) \,(\,1+4w+w^{2})\left[ \frac{1}{w^{2}-1}+%
\frac{1}{w^{3}-1}+\frac{1}{w^{5}-1}\right]  \\
\  \\
+\left( 1+w\right) \left[ \sum\limits_{(\kappa ,\lambda )}\left( \frac{%
w-w^{\kappa +1}}{w^{\kappa +1}-1}\right) \left( \frac{w-w^{\lambda +1}}{%
w^{\lambda +1}-1}\right) -21\right]  \\
\  \\
+\sum\limits_{(\kappa ,\lambda ,\mu )}\left( \frac{w-w^{\kappa +1}}{%
w^{\kappa +1}-1}\right) \left( \frac{w-w^{\lambda +1}}{w^{\lambda +1}-1}%
\right) \left( \frac{w-w^{\mu +1}}{w^{\mu +1}-1}\right) +12
\end{array}
$}
\end{equation*}
where the pairs $(\kappa ,\lambda )$ are taken from the set
\begin{equation*}
\begin{array}{l}
\{(4,9),(9,4),(4,11),(11,4),(1,11),(11,1),(4,4) \\
(1,4),(4,1),(4,13),(13,4),(2,13),(13,2),(2,2) \\
(2,5),(5,2),(1,2),(2,1),(4,2),(2,4),(1,1)\}\,.
\end{array}
\end{equation*}
and the triples $(\kappa ,\lambda ,\mu )$ from the set
\begin{equation*}
\begin{array}{l}
\{(1,1,11),(1,2,4),(1,4,2),(1,4,11),(1,11,4),(2,2,5), \\
(2,2,13),(2,4,13),(2,13,4),(4,4,9),(4,4,11),(4,4,13)\}\ .
\end{array}
\end{equation*}

\emph{(vii) }$\mathbf{Type}$ $\mathbf{E}_{8}.$%
\begin{equation*}
\fbox{$
\begin{array}{l}
E_{\text{\emph{str}}}\left( X;u,v\right) =w^{3}-1+\allowbreak \left(
w-1\right) \,(\,w+1)^{2}\left[ \frac{1}{w^{12}-1}+\frac{1}{w^{16}-1}+\frac{1%
}{w^{20}-1}+\frac{1}{w^{24}-1}\right]  \\
\  \\
+2\left( w-1\right) \,(\,1+4w+w^{2})\left[ \frac{1}{w^{2}-1}+\frac{1}{w^{3}-1%
}+\frac{1}{w^{5}-1}+\frac{1}{w^{8}-1}\right]  \\
\, \\
+\left( 1+w\right) \left[ \sum\limits_{(\kappa ,\lambda )}\left( \frac{%
w-w^{\kappa +1}}{w^{\kappa +1}-1}\right) \left( \frac{w-w^{\lambda +1}}{%
w^{\lambda +1}-1}\right) -28\right]  \\
\  \\
+\sum\limits_{(\kappa ,\lambda ,\mu )}\left( \frac{w-w^{\kappa +1}}{%
w^{\kappa +1}-1}\right) \left( \frac{w-w^{\lambda +1}}{w^{\lambda +1}-1}%
\right) \left( \frac{w-w^{\mu +1}}{w^{\mu +1}-1}\right) +17
\end{array}
$}
\end{equation*}
where the pairs $(\kappa ,\lambda )$ are taken from the set
\begin{equation*}
\begin{array}{l}
\{(1,1),(1,2),(2,1),(1,4),(4,1),(1,11),(11,1), \\
(2,2),(2,4),(4,2),(2,7),(7,2),(2,19),(19,2), \\
(4,4),(4,7),(7,4),(4,11),(11,4),(4,23),(23,4), \\
(7,7),(7,15),(15,7),(7,19),(19,7),(7,23),(23,7)\}
\end{array}
\end{equation*}
and the triples $(\kappa ,\lambda ,\mu )$ from the set
\begin{equation*}
\begin{array}{l}
\{(1,1,11),(1,2,4),(1,4,2),(1,4,11),(1,11,4),(2,2,19), \\
(2,4,7),(2,7,4),(2,7,19),(2,19,7),(4,4,11),(4,4,23), \\
(4,7,23),(4,23,7),(7,7,15),(7,7,19),(7,7,23)\}\ .\medskip
\end{array}
\end{equation*}
}
\noindent In particular, the values of the corresponding
string-theoretic Euler numbers \emph{(\ref{e-STR})} are equal to
{\small
\begin{equation*}
\begin{tabular}{|c|c|}
\hline
$\mathbf{Types}$ & $e_{\text{\emph{str}}}\left( X\right) $ \\ \hline\hline
$\mathbf{A}_{n},$ $n$ \emph{even} & $
\begin{array}{c}
\  \\
2-\tfrac{3}{n+3} \\
\
\end{array}
$ \\ \hline
$\mathbf{A}_{n},$ $n$ \emph{odd} & $
\begin{array}{c}
\, \\
2 \\
\,
\end{array}
$ \\ \hline
$\mathbf{D}_{n},$ $n$ \emph{even} & $
\begin{array}{c}
\  \\
\begin{array}{c}
-\tfrac{80n^{4}-381n^{3}+96n^{2}-128}{16n^{3}} \\
\  \\
+\allowbreak {\sum\limits_{i=1}^{\frac{n}{2}-2}}\frac{2\left(
372-492n^{2}-32i-184n^{3}+20n^{4}+688in-160in^{3}+304in^{2}+208n+5n^{5}-50in^{4}\right)
}{\left( n-2i\right) ^{3}\left( n-2i+2\right) ^{2}}
\end{array}
\!\!\! \\
\
\end{array}
$ \\ \hline
$\mathbf{D}_{n},$ $n$ \emph{odd} & $
\begin{array}{c}
\  \\
\tfrac{-96n^{3}+765n^{2}-1562n+1085}{16\left( n-1\right) ^{2}} \\
\  \\
+{\sum\limits_{i=1}^{\frac{n-5}{2}}}\tfrac{2\left(
585n-129+130n^{2}-306i-214n^{3}-5n^{4}-200in+40in^{3}+484in^{2}+5n^{5}-50in^{4}\right)
}{\left( n+1-2i\right) ^{2}\left( n-1-2i\right) ^{3}} \\
\
\end{array}
$ \\ \hline
$\mathbf{E}_{6}$ & $
\begin{array}{c}
\  \\
\frac{67}{40}=\allowbreak 1.\,\allowbreak 675 \\
\
\end{array}
$ \\ \hline
$\mathbf{E}_{7}$ & $
\begin{array}{c}
\, \\
\frac{609\,851}{189\,000}\approx 3.\,\allowbreak 226\,7 \\
\,
\end{array}
$ \\ \hline
$\mathbf{E}_{8}$ & $
\begin{array}{c}
\, \\
\allowbreak \frac{315\,467}{230\,400}\approx \allowbreak 1.\,\allowbreak
369\,2 \\
\,
\end{array}
$ \\ \hline
\end{tabular}
\end{equation*}}

\noindent and the string-theoretic indices take the following values\textit{:}%

{\footnotesize
\begin{equation*}
\begin{tabular}{|c|c|}
\hline
$\mathbf{Types}$ & \emph{ind}$_{\text{\emph{str}}}\left( X\right) $ \\
\hline\hline
$\mathbf{A}_{n}$ & $
\begin{array}{c}
\  \\
\left\{
\begin{array}{ll}
1, & \text{\emph{if}\textit{\ \ }}n\equiv 1\,(\text{\emph{mod }}2) \\
n+3, & \text{\emph{if}\textit{\ \ }}n\equiv 2\ \text{\emph{or} }4\,(\text{%
\emph{mod }}6) \\
\frac{n}{3}+1, & \text{\emph{if}\textit{\ \ }}n\equiv 0\,(\text{\emph{mod }}%
6)
\end{array}
\right.  \\
\
\end{array}
$ \\ \hline
$\mathbf{D}_{n}$ & $
\begin{array}{c}
\begin{array}{c}
\  \\
\text{\emph{It belongs to the intervall}} \\
\  \\
\left\{
\begin{array}{ll}
(n,\ n^{3}\,\prod\limits_{i=1}^{\frac{n}{2}-2}\left( n-2i\right) ^{3}\left(
n-2i+2\right) ^{2}]\cap \mathbb{Z}, & \text{\emph{if}\textit{\ \ }}n\ \text{%
\emph{even}} \\
\, & \, \\
(n,16\left( n-1\right) ^{2}\prod\limits_{i=1}^{\frac{n-5}{2}}\left(
n+1-2i\right) ^{2}\left( n-1-2i\right) ^{3}]\cap \mathbb{Z}, & \text{\emph{%
if\ }\textit{\ }}n\ \text{\emph{odd}}
\end{array}
\right.
\end{array}
\  \\
\
\end{array}
$ \\ \hline
$\mathbf{E}_{6}$ & $
\begin{array}{c}
\  \\
2^{3}5 \\
\
\end{array}
$ \\ \hline
$\mathbf{E}_{7}$ & $
\begin{array}{c}
\, \\
2^{3}3^{3}5^{3}7 \\
\,
\end{array}
$ \\ \hline
$\mathbf{E}_{8}$ & $
\begin{array}{c}
\, \\
\allowbreak \allowbreak 2^{10}3^{2}5^{2} \\
\,
\end{array}
$ \\ \hline
\end{tabular}
\end{equation*}
}
\end{theorem}

\medskip
%
%
%

\section{The canonical desingularization procedure\label{CANDES}}

\noindent Throughout this section we shall omit the superscript $d(=3)$, use
the notation (\ref{HYP}), and write the defining equation as:
\begin{equation*}
X_{f}=\left\{ \left( x_{1},x_{2},x_{3},x_{4}\right) \in\mathbb{C}^{4}\
\left| \ \right. f\left( x_{1},x_{2},x_{3},x_{4}\right) =g\left(
x_{1},x_{2}\right) +x_{3}^{2}+x_{4}^{2}=0\right\} .
\end{equation*}
$\allowbreak\allowbreak$Let $\pi:\mathbf{Bl}_{\mathbf{0}}(\mathbb{C}%
^{4})\longrightarrow\mathbb{C}^{4}$ be the blow up of $\mathbb{C}^{4}$ at
the origin, with {\small
\begin{equation*}
\mathbf{Bl}_{\mathbf{0}}(\mathbb{C}^{4})=\left\{ \left( \left(
x_{1},x_{2},x_{3},x_{4}\right) ,\left( t_{1}:t_{2}:t_{3}:t_{4}\right)
\right) \in\mathbb{C}^{4}\times\mathbb{P}_{\mathbb{C}}^{3}\ \left|
\begin{array}{l}
\ x_{i}\,t_{j}=x_{j}\,t_{i}, \\
\ \forall i,j,\ 1\leq i,j\leq4
\end{array}
\right. \right\} ,
\end{equation*}
} $\mathcal{E}=\pi^{-1}\left( \mathbf{0}\right) =\left\{ \mathbf{0}\right\}
\times\mathbb{P}_{\mathbb{C}}^{3}$, and let $U_{i}\subset \mathbf{Bl}_{%
\mathbf{0}}(\mathbb{C}^{4})$ denote the open set given by $\left(
t_{i}\neq0\right) $. In terms of analytic coordinates we may write for $%
i\in\left\{ 1,2,3,4\right\} ,$ {\small
\begin{equation*}
U_{i}=\left\{ (\left( x_{1},x_{2},x_{3},x_{4}\right) ,(\xi_{1},..,\widehat{%
\xi_{i}},..,\xi_{4}))\in\mathbb{C}^{4}\times\mathbb{C}^{3}\ \left|
\begin{array}{l}
\ x_{j}=x_{i}\,\xi_{j}, \\
\forall j,\ j\in\left\{ 1,2,3,4\right\} \mathbb{r}\left\{ i\right\}
\end{array}
\right. \right\} ,
\end{equation*}
} where $\xi_{j}=\frac{t_{j}}{t_{i}}$, and $\widehat{\xi_{i}}$ \ means that
we omit $\xi_{i}$. Moreover, we may identify $U_{i}$ with a $\mathbb{C}^{4}$
with respect to the coordinates $x_{i},\xi_{1},\ldots,\widehat{\xi_{i}}%
,\ldots ,\xi_{4}$. The restriction $\pi\left| _{U_{i}}\right. $ \ is
therefore given by mapping {\small
\begin{equation*}
\begin{array}{c}
\mathbb{C}^{4}\ni(x_{i},\xi_{1},..,\widehat{\xi_{i}},..,\xi_{4}) \\
\downarrow\ \cong \\
\begin{array}{c}
(\left(
x_{i}\,\xi_{1},..,x_{i}\,\xi_{i-1},x_{i},x_{i}\,\xi_{i+1},..,x_{i}\,\xi_{4}%
\right) ,(\xi_{1}:..:\underset{i\text{-th pos.}}{\underbrace{1}}%
:..:\xi_{4}))\in U_{i} \\
\downarrow\ \pi\left| _{U_{i}}\right. \\
\left(
x_{i}\,\xi_{1},\ldots,x_{i}\,\xi_{i-1},x_{i},x_{i}\,\xi_{i+1},\ldots,x_{i}\,%
\xi_{4}\right)
\end{array}
\end{array}
\end{equation*}
} Note that $\mathcal{E}_{i}:=\mathcal{E}\cap U_{i}$ is described as the
coordinate hyperplane $\left( x_{i}=0\right) $; i.e., the open cover $%
\left\{ U_{i}\right\} _{1\leq i\leq4}$ \ of $\mathbf{Bl}_{\mathbf{0}}(%
\mathbb{C}^{4})$ restricts to $\mathcal{E}$ to provide the standard open
cover of $\mathbb{P}_{\mathbb{C}}^{3}\ $by affine spaces $\mathbb{C}^{3}$,
with $\{\xi_{j}\}_{j\in\left\{ 1,2,3,4\right\} \mathbb{r}\left\{ i\right\} }$
being the analytic coordinates of $\mathcal{E}_{i}$.\smallskip

\noindent\textit{Notation. \ }To work with a more convenient notation we
define
\begin{equation*}
\mathbf{Bl}_{\mathbf{0}}(\mathbb{C}^{4})=\bigcup_{i=1}^{4}\ U_{i},\ \ \ \ \
U_{i}=\text{Spec}\left( \mathbb{C}\left[ y_{i,1},y_{i,2},y_{i,3},y_{i,4}%
\right] \right) ,
\end{equation*}
by setting as coordinates for $U_{i}$'s:
\begin{equation*}
y_{i,k}:=\left\{
\begin{array}{ll}
x_{k}, & \text{for\ }i=k \\
\xi_{k}, & \text{for\ }i\neq k
\end{array}
\right.
\end{equation*}
\medskip

\noindent \underline{$\bullet $ \textbf{Step 1: The first blow-up.}} \
Blowing up $X_{f}$ at the origin, we take the diagram\medskip\
\begin{equation*}
\begin{array}{ccccc}
\mathcal{E} & \subset & \mathbf{Bl}_{\mathbf{0}}(\mathbb{C}^{4}) & \overset{%
\pi }{\longrightarrow } & \mathbb{C}^{4} \\
\cup &  & \cup &  & \cup \\
\mathcal{E}\cap \mathbf{Bl}_{\mathbf{0}}(X_{f}) & \subset & \mathbf{Bl}_{%
\mathbf{0}}(X_{f}) & \overset{\pi \left| _{\text{restr.}}\right. }{%
\longrightarrow } & X_{f}
\end{array}
\end{equation*}
and consider the strict transform\medskip\
\begin{equation*}
\mathbf{Bl}_{\mathbf{0}}(X_{f})=\overline{\pi ^{-1}(X_{f}\cap (\mathbb{C}^{4}%
\mathbb{r}\left\{ \mathbf{0}\right\} ))}=\overline{\pi ^{-1}(X_{f})\cap (%
\mathbf{Bl}_{\mathbf{0}}(\mathbb{C}^{4})\mathbb{r}\mathcal{E}))}\medskip
\end{equation*}
of $X_{f}$ in $\mathbb{C}^{4}$ under $\pi $, and the corresponding
exceptional (not necessarily prime) divisor $\mathcal{E}_{f}:=\mathcal{E}%
\cap \mathbf{Bl}_{\mathbf{0}}(X_{f})$ with respect to $\pi \left| _{\text{%
restr.}}\right. $\smallskip\ \bigskip

\noindent $\blacktriangleright $ \textbf{Local description of }$\mathbf{Bl}_{%
\mathbf{0}}(X_{f})$ \textbf{and} $\mathcal{E}_{f}.$ \ After pulling back $f$
by $\pi $ and restricting ourselves onto $U_{i}$, we get\medskip
\begin{equation*}
\pi ^{\ast }(f)\left| _{U_{i}}\right. =x_{i}^{2}\ \widetilde{f_{i}}%
=y_{i,i}^{2}\ \widetilde{f_{i}},
\end{equation*}
\ with $\ \widetilde{f_{i}}\in \mathbb{C}\left[
y_{i,1},y_{i,2},y_{i,3},y_{i,4}\right] $. More precisely, we obtain{\small %
\bigskip \medskip\
\begin{equation*}
\begin{tabular}{|c|c|c|}
\hline
\textbf{Types} & $\ \
\begin{array}{c}
\, \\
\ \ \widetilde{f_{1}} \\
\,
\end{array}
$ & $\ \
\begin{array}{c}
\, \\
\ \widetilde{f_{2}} \\
\,
\end{array}
$ \\ \hline\hline
$\mathbf{A}_{n}$ & $
\begin{array}{c}
\  \\
y_{1,1}^{n-1}+y_{1,2}^{2}+y_{1,3}^{2}+y_{1,4}^{2} \\
\
\end{array}
$ & $y_{2,1}^{n+1}y_{2,2}^{n-1}+1+y_{2,3}^{2}+y_{2,4}^{2}$ \\ \hline
$\mathbf{D}_{n}$ & $
\begin{array}{c}
\  \\
y_{1,1}^{n-3}+y_{1,1}y_{1,2}^{2}+y_{1,3}^{2}+y_{1,4}^{2} \\
\
\end{array}
$ & $y_{2,1}^{n-1}y_{2,2}^{n-3}+y_{2,1}y_{2,2}+y_{2,3}^{2}+y_{2,4}^{2}$ \\
\hline
$\mathbf{E}_{6}$ & $
\begin{array}{c}
\  \\
y_{1,1}+y_{1,1}^{2}y_{1,2}^{4}+y_{1,3}^{2}+y_{1,4}^{2} \\
\
\end{array}
$ & $y_{2,1}^{3}y_{2,2}+y_{2,2}^{2}+y_{2,3}^{2}+y_{2,4}^{2}$ \\ \hline
$\mathbf{E}_{7}$ & $
\begin{array}{c}
\  \\
y_{1,1}+y_{1,1}^{2}y_{1,2}^{3}+y_{1,3}^{2}+y_{1,4}^{2} \\
\
\end{array}
$ & $y_{2,1}^{3}y_{2,2}+y_{2,1}y_{2,2}^{2}+y_{2,3}^{2}+y_{2,4}^{2}$ \\ \hline
$\mathbf{E}_{8}$ & $
\begin{array}{c}
\  \\
y_{1,1}+y_{1,1}^{3}y_{1,2}^{5}+y_{1,3}^{2}+y_{1,4}^{2} \\
\
\end{array}
$ & $y_{2,1}^{3}y_{2,2}+y_{2,2}^{3}+y_{2,3}^{2}+y_{2,4}^{2}$ \\ \hline
\end{tabular}
\end{equation*}
} and

{\small
\begin{equation*}
\begin{tabular}{|c|c|c|}
\hline
\textbf{Types} & $\ \
\begin{array}{c}
\, \\
\ \ \widetilde{f_{3}} \\
\,
\end{array}
$ & $
\begin{array}{c}
\, \\
\ \widetilde{f_{4}} \\
\,
\end{array}
$ \\ \hline\hline
$\mathbf{A}_{n}$ & $
\begin{array}{c}
\  \\
y_{3,1}^{n+1}y_{3,3}^{n-1}+y_{3,2}^{2}+1+y_{3,4}^{2} \\
\
\end{array}
$ & $
\begin{array}{c}
\  \\
y_{4,1}^{n+1}y_{4,4}^{n-1}+y_{4,2}^{2}+y_{4,3}^{2}+1 \\
\
\end{array}
$ \\ \hline
$\mathbf{D}_{n}$ & $
\begin{array}{c}
\  \\
\!\!y_{3,1}^{n-1}y_{3,3}^{n-3}+y_{3,1}y_{3,2}^{2}y_{3,3}+1+y_{3,4}^{2}\!\!
\\
\
\end{array}
$ & $
\begin{array}{c}
\  \\
\!\!y_{4,1}^{n-1}y_{4,4}^{n-3}+y_{4,1}y_{4,2}^{2}y_{4,4}+y_{4,3}^{2}+1\!\!\!
\\
\
\end{array}
$ \\ \hline
$\mathbf{E}_{6}$ & $
\begin{array}{c}
\  \\
y_{3,1}^{3}y_{3,3}+y_{3,2}^{4}y_{3,3}^{2}+1+y_{3,4}^{2} \\
\
\end{array}
$ & $
\begin{array}{c}
\  \\
y_{4,1}^{3}y_{4,4}+y_{4,4}^{4}y_{4,4}^{2}+y_{4,3}^{2}+1 \\
\
\end{array}
$ \\ \hline
$\mathbf{E}_{7}$ & $
\begin{array}{c}
\  \\
y_{3,1}^{3}y_{3,3}+y_{3,1}y_{3,2}^{2}y_{3,3}^{2}+1+y_{3,4}^{2} \\
\
\end{array}
$ & $
\begin{array}{c}
\  \\
y_{4,1}^{3}y_{4,4}+y_{4,1}y_{4,2}^{2}y_{4,4}^{2}+y_{4,3}^{2}+1 \\
\
\end{array}
$ \\ \hline
$\mathbf{E}_{8}$ & $
\begin{array}{c}
\  \\
y_{3,1}^{3}y_{3,3}+y_{3,2}^{5}y_{3,3}^{3}+1+y_{3,4}^{2} \\
\
\end{array}
$ & $
\begin{array}{c}
\  \\
y_{4,1}^{3}y_{4,4}+y_{4,2}^{5}y_{4,4}^{3}+y_{4,3}^{2}+1 \\
\
\end{array}
$ \\ \hline
\end{tabular}
\medskip
\end{equation*}
} \noindent Locally, {\small
\begin{equation*}
\mathbf{Bl}_{\mathbf{0}}(X_{f})\left| _{U_{i}}\right. \overset{\cong }{%
\rightarrow}\left\{ \left( y_{i,1},y_{i,2},y_{i,3},y_{i,4}\right) \in\mathbb{%
C}^{4}\ \left| \ \ \widetilde{f_{i}}\left(
y_{i,1},y_{i,2},y_{i,3},y_{i,4}\right) =0\right. \right\} ,
\end{equation*}
} and using the restrictions of $\widetilde{f_{i}}$'s on the $\mathcal{E}%
_{i} $'s, $i=1,2,3,4$, we get the equations for $\mathcal{E}_{f}\left|
_{U_{i}}\right. $: {\small
\begin{equation*}
\mathbf{Bl}_{\mathbf{0}}(X_{f})\cap\mathcal{E}_{i}=\mathcal{E}_{f}\left|
_{U_{i}}\right. \!\overset{\cong}{\rightarrow}\left\{ \!\left(
y_{i,1},y_{i,2},y_{i,3},y_{i,4}\right) \!\in\mathbb{C}^{4}\left| \ y_{i,i}=%
\widetilde{f_{i}}\left( y_{i,1},y_{i,2},y_{i,3},y_{i,4}\right) =0\right.
\!\right\} .
\end{equation*}
}

\begin{lemma}[Local Reduction]
\label{LRL}The types of the singularities of $\mathbf{Bl}_{\mathbf{0}}(X_{f})
$ are given by the following table\emph{:\medskip } {\small
\begin{equation*}
\begin{tabular}{|c|c|c|}
\hline
$
\begin{array}{c}
\mathbf{Initial}\,\ \mathbf{types}\text{ } \\
\mathbf{of}\text{ }\mathbf{singularities}\text{ } \\
\mathbf{of\ }\text{ }X_{f}
\end{array}
$ & $
\begin{array}{c}
\mathbf{New\ \,singularities\ } \\
\mathbf{(and\ their\ types)} \\
\mathbf{on\ \ \ Bl}_{\mathbf{0}}(X_{f})
\end{array}
$ & $
\begin{array}{c}
\mathbf{located\ \ in\ the\ } \\
\mathbf{affine\ pieces}
\end{array}
$ \\ \hline\hline
$\mathbf{A}_{1},\mathbf{A}_{2}$ & $-$ & $-$ \\ \hline
$\mathbf{A}_{n},\,n\geq 3$ & $\mathbf{A}_{n-2}$ & $U_{1}$ \\ \hline
$\mathbf{D}_{4}$ & $\mathbf{A}_{1},\ \mathbf{A}_{1},\ \mathbf{A}_{1}$ & $%
U_{2},\,U_{1}\cap U_{2},\,U_{1}\cap U_{2}$ \\ \hline
$\mathbf{D}_{5}$ & $\mathbf{A}_{3},\ \mathbf{A}_{1}$ & $U_{1},U_{2}$ \\
\hline
$\mathbf{D}_{n},\,n\geq 6$ & $\mathbf{D}_{n-2},\ \mathbf{A}_{1}$ & $%
U_{1},U_{2}$ \\ \hline
$\mathbf{E}_{6}$ & $\mathbf{A}_{5}$ & $U_{2}$ \\ \hline
$\mathbf{E}_{7}$ & $\mathbf{D}_{6}$ & $U_{2}$ \\ \hline
$\mathbf{E}_{8}$ & $\mathbf{E}_{7}$ & $U_{2}$ \\ \hline
\end{tabular}
\end{equation*}
}
\end{lemma}

\noindent\textit{Proof.} The affine pieces in which the singularities of $%
\mathbf{Bl}_{\mathbf{0}}(X_{f})$ are located are obviously those of the
above table (simply by partial derivative checking). Let us now examine the
types of the appearing singularities in each case separately.\smallskip

\noindent$\triangleright$ Blowing up singularity $\mathbf{A}_{n},\,n\geq3,$
we obtain an $\mathbf{A}_{n-2}$-singularity in its normal form $\ \widetilde
{f_{1}}.\smallskip$

\noindent$\triangleright$ Blowing up $\mathbf{D}_{n}$'s, and working first
with the patch $U_{1}$, we get a $\mathbf{D}_{n-2}$-singularity in its
normal form $\ \widetilde{f_{1}}$ whenever $n\geq6,$ no singularity for $%
n=4, $ and an $\mathbf{A}_{3}$-singularity for $n=5,$ just by utilizing the
analytic coordinate change
\begin{equation*}
y_{1,i}=\left\{
\begin{array}{ll}
y_{1,i}^{\prime}, & i\in\{2,3,4\} \\
\, & \, \\
y_{1,1}^{\prime}-\dfrac{1}{2}(y_{1,2}^{\prime})^{2}, & i=1
\end{array}
\right.
\end{equation*}
and writing the corresponding defining polynomial as:
\begin{equation*}
y_{1,1}^{2}+y_{1,1}y_{1,2}^{2}+y_{1,3}^{2}+y_{1,4}^{2}=-\frac{1}{4}%
(y_{1,2}^{\prime})^{4}+(y_{1,1}^{\prime})^{2}+(y_{1,3}^{%
\prime})^{2}+(y_{1,4}^{\prime})^{2}\ .
\end{equation*}
Passing to $U_{2},$ we have\ {\small $\mathbf{Bl}_{\mathbf{0}}(X_{f})\left|
_{U_{2}}\right. =$
\begin{equation*}
=\left\{ \left( y_{2,1},...,y_{2,4}\right) \in\mathbb{C}^{4}\left|
\theta(y_{2,1},...,y_{2,4}):=y_{2,1}^{n-1}y_{2,2}^{n-3}+y_{2,1}y_{2,2}+y_{2,3}^{2}+y_{2,4}^{2}=0\right. \right\}
\end{equation*}
} with partial derivatives w.r.t. $\theta=\theta(y_{2,1},...,y_{2,4})$:
{\small
\begin{equation*}
\left\{
\begin{array}{l}
\dfrac{\partial\theta}{\partial y_{2,1}}=(n-1)\
y_{2,1}^{n-2}y_{2,2}^{n-3}+y_{2,2}=y_{2,2}\,\left( (n-1)\
y_{2,1}^{n-2}y_{2,2}^{n-4}+1\right) \\
\  \\
\dfrac{\partial\theta}{\partial y_{2,2}}=(n-3)\
y_{2,1}^{n-1}y_{2,2}^{n-4}+y_{2,1}=y_{2,1}\,\left( (n-3)\
y_{2,1}^{n-2}y_{2,2}^{n-4}+1\right) \\
\  \\
\dfrac{\partial\theta}{\partial y_{2,3}}=2\,y_{2,3}\text{ \ \ \ \ and \ \ \
\ }\dfrac{\partial\theta}{\partial y_{2,4}}=2\,y_{2,4}\ .
\end{array}
\right.
\end{equation*}
} Clearly, for $n=4,$ the singular locus of $\mathbf{Bl}_{\mathbf{0}%
}(X_{f})\left| _{U_{2}}\right. $ consists of the points
\begin{equation*}
\left( 0,0,0,0\right) ,\ \ \ \ (\sqrt{-1},0,0,0)\ \ \ \ \text{ and \ \ \ \ }%
(-\sqrt{-1},0,0,0)
\end{equation*}
which can be expressed as the singularities \textit{at the origin} $\mathbf{0%
}$ of $\mathbb{C}^{4}$ for {\small \
\begin{equation}
\left\{
\begin{array}{l}
y_{2,1}^{3}y_{2,2}+y_{2,1}y_{2,2}+y_{2,3}^{2}+y_{2,4}^{2}=0 \\
\  \\
\ y_{2,2}^{\prime}\,(y_{2,1}^{\prime})^{3}\pm3\,\sqrt{-1}\,y_{2,2}^{\prime
}\,(y_{2,1}^{\prime})^{2}-2\,y_{2,2}^{\prime}\,y_{2,1}^{\prime}+(y_{2,3}^{%
\prime})^{2}+(y_{2,4}^{\prime})^{2}=0
\end{array}
\right. \allowbreak\ \   \label{GLEICHUNGEN}
\end{equation}
} (just by setting $y_{2,1}=y_{2,1}^{\prime}\pm\sqrt{-1}$ and $%
y_{2,i}=y_{2,i}^{\prime},$ for $i\in\{2,3,4\}$). Next, applying a result of
B\u {a}descu (in a very special case of it, \cite[Thm. 1, p. 209]{BADESCU}),
we see that \textit{all} normal isolated singularities which can be fully
resolved after a single blow-up and have exceptional divisor $E\cong \mathbb{%
P}_{\mathbb{C}}^{1}\times\mathbb{P}_{\mathbb{C}}^{1}$ with conormal bundle $%
\mathcal{N}_{E}^{\vee}$ isomorphic to $\mathcal{O}_{E}\left( 1,1\right) $
are analytically isomorphic to each other. It is easy to verify that this is
valid for all singularities (\ref{GLEICHUNGEN}). Hence, they are all
analytically isomorphic to an $\mathbf{A}_{1}$-singularity (which has the
same property). Alternatively, one can show that these are analytically
isomorphic to $\mathbf{A}_{1}$-singularities by exploiting the fact that
they are semiquasihomogeneous of weight $(\frac{1}{2},\frac{1}{2},\frac{1}{2}%
,\frac{1}{2})$ and by using \cite[Corollary 3.3]{ROCZEN3}. (The completions
are isomorphic to the singularities defined by that polynomial part
consisting of all terms of weight $1$, which is obviously equal to $%
y_{2,1}y_{2,2}+y_{2,3}^{2}+y_{2,4}^{2}$ and $-2\,y_{2,2}^{\prime}\,y_{2,1}^{%
\prime }+(y_{2,3}^{\prime})^{2}+(y_{2,4}^{\prime})^{2},$ respectively). On
the other hand, for $n\geq5$, the only singular point of $\mathbf{Bl}_{%
\mathbf{0}}(X_{f})\left| _{U_{2}}\right. $ is $\left( 0,0,0,0\right) ,$
which again turns out to be an $\mathbf{A}_{1}$-singularity (by the same
reasoning).\smallskip

\noindent$\triangleright$ Now the singularity $\mathbf{E}_{6}$ passes after
blowing up to an $\mathbf{A}_{5}$-singularity, because using the analytic
coordinate change
\begin{equation*}
y_{2,i}=\left\{
\begin{array}{ll}
y_{2,i}^{\prime}, & i\in\{1,3,4\} \\
\, & \, \\
y_{2,2}^{\prime}-\dfrac{1}{2}(y_{2,1}^{\prime})^{3}, & i=2
\end{array}
\right.
\end{equation*}
we get
\begin{equation*}
y_{2,1}^{3}y_{2,2}+y_{2,2}^{2}+y_{2,3}^{2}+y_{2,4}^{2}=-\frac{1}{4}%
(y_{2,1}^{\prime})^{6}+(y_{2,2}^{\prime})^{2}+(y_{2,3}^{%
\prime})^{2}+(y_{2,4}^{\prime})^{2}\ .
\end{equation*}
$\triangleright$ Starting with $\mathbf{E}_{7}$ we obtain a $\mathbf{D}_{6}$%
-singularity, because the analytic coordinate change
\begin{equation*}
y_{2,i}=\left\{
\begin{array}{ll}
y_{2,i}^{\prime}, & i\in\{1,3,4\} \\
\, & \, \\
y_{2,2}^{\prime}-\dfrac{1}{2}(y_{2,1}^{\prime})^{2}, & i=2
\end{array}
\right.
\end{equation*}
implies
\begin{equation*}
y_{2,1}^{3}y_{2,2}+y_{2,1}y_{2,2}^{2}+y_{2,3}^{2}+y_{2,4}^{2}=-\frac{1}{4}%
(y_{2,1}^{\prime})^{5}+y_{2,1}^{\prime}\,(y_{2,2}^{\prime})^{2}+(y_{2,3}^{%
\prime})^{2}+(y_{2,4}^{\prime})^{2}\ .
\end{equation*}
$\triangleright$ Finally, blowing up singularity $\mathbf{E}_{8},$ we
acquire an $\mathbf{E}_{7}$-singularity in its normal form $\ \widetilde{%
f_{2}}.$\hfill$\square$\bigskip

\noindent $\blacktriangleright $ \textbf{Global description of }$\mathbf{Bl}%
_{\mathbf{0}}(X_{f})$ \textbf{and} $\mathcal{E}_{f}.$ This can be realized
after coming back to our global coordinates:{\small
\begin{equation*}
\begin{tabular}{|c|c|}
\hline
\textbf{Types} & $\mathbf{Bl}_{\mathbf{0}}(X_{f})=$ all $\ \left( \left(
x_{1},..,x_{4}\right) ,\left( t_{1}:t_{2}:t_{3}:t_{4}\right) \right) \in \
\mathbf{Bl}_{\mathbf{0}}(\mathbb{C}^{4})\ \ $with:$\ $ \\ \hline\hline
$
\begin{array}{c}
\, \\
\mathbf{A}_{n} \\
\,
\end{array}
$ & $x_{1}^{n-1}\,t_{1}^{2}+t_{2}^{2}+t_{3}^{2}+t_{4}^{2}=0$ \\ \hline
$
\begin{array}{c}
\, \\
\mathbf{D}_{n} \\
\,
\end{array}
$ & $x_{1}^{n-3}\,t_{1}^{2}+x_{1}\,t_{2}^{2}+t_{3}^{2}+t_{4}^{2}=0$ \\ \hline
$
\begin{array}{c}
\, \\
\mathbf{E}_{6} \\
\,
\end{array}
$ & $x_{1}\,t_{1}^{2}+x_{2}^{2}\,t_{2}^{2}+t_{3}^{2}+t_{4}^{2}=0$ \\ \hline
$
\begin{array}{c}
\, \\
\mathbf{E}_{7} \\
\,
\end{array}
$ & $x_{1}\,t_{1}^{2}+x_{1}\,x_{2}\,t_{2}^{2}+t_{3}^{2}+t_{4}^{2}=0$ \\
\hline
$
\begin{array}{c}
\, \\
\mathbf{E}_{8} \\
\,
\end{array}
$ & $x_{1}\,t_{1}^{2}+x_{2}^{3}\,t_{2}^{2}+t_{3}^{2}+t_{4}^{2}=0$ \\ \hline
\end{tabular}
\medskip
\end{equation*}
\noindent }In particular, this means that the exceptional locus $\mathcal{E}%
_{f}$ is given globally\textit{\ }by{\small
\begin{equation*}
\begin{tabular}{|c|c|}
\hline
\textbf{Types \ of \ }$X_{f}$\textbf{'s} & $\mathcal{E}_{f}$ $=$ all $\
\left( \mathbf{0},\left( t_{1}:t_{2}:t_{3}:t_{4}\right) \right) \in \
\left\{ \mathbf{0}\right\} \times \mathbb{P}_{\mathbb{C}}^{3}\ \ $with:$\ $
\\ \hline\hline
$
\begin{array}{c}
\  \\
\mathbf{A}_{1} \\
\
\end{array}
$ & $t_{1}^{2}+t_{2}^{2}+t_{3}^{2}+t_{4}^{2}=0$ \\ \hline
$
\begin{array}{c}
\  \\
\mathbf{A}_{n},n\geq 2 \\
\
\end{array}
$ & $t_{2}^{2}+t_{3}^{2}+t_{4}^{2}=0$ \\ \hline
$
\begin{array}{c}
\  \\
\mathbf{D}_{n},\mathbf{E}_{6},\mathbf{E}_{7},\mathbf{E}_{8} \\
\
\end{array}
$ & $t_{3}^{2}+t_{4}^{2}=\left( t_{3}+\sqrt{-1}\ t_{4}\right) \,\left( t_{3}-%
\sqrt{-1}\ t_{4}\right) =0$ \\ \hline
\end{tabular}
\end{equation*}
}In the latter four cases $\mathcal{E}_{f}$ $\ $consists of two exceptional
prime divisors, say $\mathcal{E}_{f}^{\prime }$ $\ $and $\mathcal{E}%
_{f}^{\prime \prime }$ (which are $\cong $ $\mathbb{P}_{\mathbb{C}}^{2}$).
Moreover, taking into account the above local description of singularities
of $\mathbf{Bl}_{\mathbf{0}}(X_{f}),$ we may rewrite them in homogeneous
coordinates on $\left\{ \mathbf{0}\right\} \times \mathbb{P}_{\mathbb{C}%
}^{3} $ as follows:{\small
\begin{equation*}
\begin{tabular}{|c|c|}
\hline
\textbf{Types \ of \ }$X_{f}$\textbf{'s} & \textbf{Singular points of\ } $%
\mathbf{Bl}_{\mathbf{0}}(X_{f})$ \\ \hline\hline
$\mathbf{A}_{1},\ \mathbf{A}_{2}$ & $-$ \\ \hline
$
\begin{array}{c}
\  \\
\mathbf{A}_{n},\,n\geq 3 \\
\
\end{array}
$ & $\ \left( \mathbf{0},\left( 1:0:0:0\right) \right) \in \mathcal{E}_{f}$
\\ \hline
$\mathbf{D}_{4}$ & $
\begin{array}{c}
\  \\
\left( \mathbf{0},\left( 0:1:0:0\right) \right) ,\ (\mathbf{0},(\pm \sqrt{-1}%
:1:0:0))\in \mathcal{E}_{f}^{\prime }\cap \mathcal{E}_{f}^{^{\prime \prime }}
\\
\
\end{array}
$ \\ \hline
$
\begin{array}{c}
\  \\
\mathbf{D}_{n},\,n\geq 5 \\
\
\end{array}
$ & $\left( \mathbf{0},\left( 1:0:0:0\right) \right) ,\ \left( \mathbf{0}%
,\left( 0:1:0:0\right) \right) \in \mathcal{E}_{f}^{\prime }\cap \mathcal{E}%
_{f}^{^{\prime \prime }}$ \\ \hline
$
\begin{array}{c}
\  \\
\mathbf{E}_{6},\ \mathbf{E}_{7},\ \mathbf{E}_{8} \\
\
\end{array}
$ & $\left( \mathbf{0},\left( 0:1:0:0\right) \right) \in \mathcal{E}%
_{f}^{\prime }\cap \mathcal{E}_{f}^{^{\prime \prime }}$ \\ \hline
\end{tabular}
\medskip
\end{equation*}
}

\noindent{}\underline{$\bullet$ \textbf{Step 2: The next blow-ups.}} The
desired snc-desingularizations of $X_{f}$'s, say $\varphi:\widetilde
{X}\rightarrow X_{f}$, will be constructed by blowing up the possibly new
singular points again and again until we reach to a smooth threefold $%
\widetilde{X}$ with exceptional locus $\mathfrak{Ex}\left( \varphi\right) $
consisting of smooth prime divisors with normal crossings. We give a
complete characterization of $\varphi$'s by the following data:\medskip

\noindent$\triangleright$ the \textit{local resolution diagrams}
(abbreviated \textit{LR-diagrams}) which are constructed after repeated
applications of Lemma \ref{LRL} (with each arrow indicating a \textit{local}
blow-up at a single closed point),\medskip

\noindent$\triangleright$ the \textit{intersection }(\textit{plane})\textit{%
\ graphs} whose vertices represent the exceptional prime divisors w.r.t. the
$\varphi$'s and their edges insinuate that the corresponding vertices are
divisors which have non-empty intersection,\medskip

\noindent$\triangleright$ the \textit{structure} of the exceptional prime
divisors up to biregular isomorphism (which turn out to be certain compact
rational surfaces of Picard number either $2$ or $4$), and finally\medskip

\noindent$\triangleright$ the \textit{intersection cycles} of all
intersecting pairs of exceptional prime divisors $(D_{i}\cdot D_{j})\left|
_{D_{k}}\right. ,\,k\in\{i,j\},$ as \textit{divisors} on $D_{k}$ (cf. \cite
{ROCZEN1}, \cite{ROCZEN2}), though we are primarily interested in their
underlying topological spaces (see below lemma \ref{LEMMA1}).\medskip\

\noindent The interplay of local and global data (simultaneous blow-ups,
strict transfoms after each step etc.) will be explained explicitly only for
types $\mathbf{A}_{n},\mathbf{D}_{4},\mathbf{E}_{6}$. (For reasons of
economy, further details -in this connection- about the other types will be
omitted. The not so difficult verification of the way one builds the
corresponding intersection graphs step by step is left to the
reader).\bigskip

\noindent\textbf{(i) Type }$\mathbf{A}_{1}$\textbf{. }Blowing up the origin
once$,$ we achieve immediately the required desingularization. The
exceptional prime divisor
\begin{equation*}
\mathcal{E}_{f}\cong\left\{ \left( t_{1}:t_{2}:t_{3}:t_{4}\right) \in\
\mathbb{P}_{\mathbb{C}}^{3}\ \left| \ \right.
t_{1}^{2}+t_{2}^{2}+t_{3}^{2}+t_{4}^{2}=0\right\}
\end{equation*}
is biregularly isomorphic to $\left\{ \left( t_{1}^{\prime}:t_{2}^{\prime
}:t_{3}^{\prime}:t_{4}^{\prime}\right) \in\mathbb{P}_{\mathbb{C}}^{3}\
\left| \ \right.
t_{1}^{\prime}\,t_{2}^{\prime}-t_{3}^{\prime}\,t_{4}^{\prime}=0\right\} =$ Im%
$\left( \gamma\right) ,$ where $\gamma$ denotes the Segre embedding
\begin{equation*}
\mathbb{P}_{\mathbb{C}}^{1}\times\mathbb{P}_{\mathbb{C}}^{1}\ni(\left(
\varpi_{1}:\varpi_{2}\right) ,\left( \varpi_{1}^{\prime}:\varpi_{2}^{\prime
}\right) )\overset{\gamma}{\longmapsto}(z_{1}:z_{2}:z_{3}:z_{4})\in \mathbb{P%
}_{\mathbb{C}}^{3}
\end{equation*}
with
\begin{equation*}
\left\{
\begin{array}{l}
z_{1}=\varpi_{1}\,\varpi_{1}^{\prime},\,\
z_{2}=\varpi_{1}\,\varpi_{2}^{\prime},\,\
z_{3}=\varpi_{2}\,\varpi_{1}^{\prime},\,\ z_{4}=\varpi
_{2}\,\varpi_{2}^{\prime}, \\
\  \\
\ t_{1}^{\prime}=z_{1},\,t_{2}^{\prime}=z_{4},\,\ t_{3}^{\prime}=z_{2},\,\
t_{4}^{\prime}=z_{3}.
\end{array}
\right.
\end{equation*}
\ Indeed, defining $\delta$ to be the biregular isomorphism {\small
\begin{equation*}
\left( t_{1}^{\prime}:t_{2}^{\prime}:t_{3}^{\prime}:t_{4}^{\prime}\right)
\overset{\delta}{\longmapsto}\left( t_{1}-\sqrt{-1}\,t_{2}:t_{1}+\sqrt
{-1}\,t_{2}:t_{3}-\sqrt{-1}\,t_{4}:-(t_{3}+\sqrt{-1}\,t_{4})\right) ,
\end{equation*}
} we obtain $\delta\left( \text{Im}\left( \gamma\right) \right) =\mathcal{E}%
_{f}.$ Consequently, $\mathcal{E}_{f}\cong\mathbb{P}_{\mathbb{C}}^{1}\times%
\mathbb{P}_{\mathbb{C}}^{1}$ and has conormal bundle $\mathcal{O}_{\mathcal{E%
}_{f}}(1,1).\bigskip$

\noindent\textbf{(ii) Type }$\mathbf{A}_{2}$\textbf{. }Blowing up the origin
once, $\mathbf{Bl}_{\mathbf{0}}(X_{f})$ is smooth (as threefold), though
\begin{equation*}
\mathcal{E}_{f}=\left\{ \left( \mathbf{0},\left(
t_{1}:t_{2}:t_{3}:t_{4}\right) \right) \in\ \left\{ \mathbf{0}\right\}
\times \mathbb{P}_{\mathbb{C}}^{3}\ \left| \ \right.
t_{2}^{2}+t_{3}^{2}+t_{4}^{2}=0\right\} \subset\mathbf{Bl}_{\mathbf{0}%
}(X_{f})
\end{equation*}
(as surface on the threefold $\mathbf{Bl}_{\mathbf{0}}(X_{f})$) has a
singular, ordinary double point at $q=\left( \mathbf{0},(1:0:0:0\right) )$
in $\mathcal{E}_{f}\left| _{U_{1}}\right. .$ For this reason, in order to
form an snc-resolution of the original singularity, we have to blow-up once
more our threefold at $q$ and consider
\begin{equation*}
\varphi:\widetilde{X}=\mathbf{Bl}_{\,q}(\mathbf{Bl}_{\mathbf{0}%
}(X_{f}))\longrightarrow X_{f}\ .
\end{equation*}
\ The new exceptional prime divisor is obviously a $\mathbb{P}_{\mathbb{C}%
}^{2}$, while the strict transform of the old one is nothing but the ($2$%
-dimensional) blow-up of $\mathcal{E}_{f}$ at $q.$ Since $\mathcal{E}_{f}$
can be viewed as the projective cone $\subset\mathbb{P}_{\mathbb{C}}^{3}$
over the smooth quadratic hypersurface $V=\left\{ \left(
t_{2}:t_{3}:t_{4}\right) \in\mathbb{P}_{\mathbb{C}}^{2}\ \left| \ \right.
t_{2}^{2}+t_{3}^{2}+t_{4}^{2}=0\right\} $ with $(1:0:0:0)$ as its vertex,
blowing up $(1:0:0:0),$ we obtain a ruled (compact) surface over $V\cong
\mathbb{P}_{\mathbb{C}}^{1}$ having the inverse image of $(1:0:0:0)$ as a
section $\mathsf{C}_{0}$ with self-intersection $\mathsf{C}_{0}^{2}=-2$ (see
Hartshorne \cite[V.2.11.4, pp. 374-375]{HARTSHORNE}). Hence, the strict
transform of $\mathcal{E}_{f}$ under $\varphi$ has to be the rational ruled
surface $\mathbb{F}_{2}:=\mathbb{P}(\mathcal{O}_{\mathbb{P}_{\mathbb{C}%
}^{1}}\oplus\mathcal{O}_{\mathbb{P}_{\mathbb{C}}^{1}}\left( -2\right) )$
(because $\mathbb{F}_{2}$ is the unique $\mathbb{P}_{\mathbb{C}}^{1}$-bundle
over $\mathbb{P}_{\mathbb{C}}^{1}$ having an irreducible curve of
self-intersection $-2,$ cf. \cite[p. 519]{GR-H}).

\begin{remark}
\label{REMARK-S}Among the three-dimensional $\mathbf{A}$-$\mathbf{D}$-$%
\mathbf{E}$'s, type $\mathbf{A}_{2},$ and, in general, type $\mathbf{A}%
_{n},\,n$ even, constitutes the only exception in which one has to blow up a
\textit{smooth} threefold point at the last step to ensure an
snc-resolution. In all the other cases the snc-condition will be present
immediately after the last blow-ups of singular points (becoming clear from
the LR-diagrams which have only $\mathbf{A}_{1}$'s at their last but one
ends).
\end{remark}

\noindent \textbf{(iii) Types }$\mathbf{A}_{n},\,n\geq 3$\textbf{. }The
LR-diagram for these types depends on the $($mod $2$)-behaviour of $n,$ and
the number of the required blow-ups equals $m:=\left\lfloor \frac{n+2}{2}%
\right\rfloor .$%
\begin{equation*}
\fbox{$
\begin{array}{cc}
\mathbf{A}_{n}\rightarrow \mathbf{A}_{n-2}\rightarrow \mathbf{A}%
_{n-4}\rightarrow \cdots \rightarrow \mathbf{A}_{3}\rightarrow \mathbf{A}%
_{1}\rightarrow \mathbf{A}_{0} & (\text{if\ \ }n\equiv 1(\text{mod 2)}) \\
\  & \, \\
\mathbf{A}_{n}\rightarrow \mathbf{A}_{n-2}\rightarrow \mathbf{A}%
_{n-4}\rightarrow \cdots \rightarrow \mathbf{A}_{2}\rightarrow \mathbf{A}%
_{0}\rightarrow \mathbf{A}_{0} & (\text{if\ \ }n\equiv 0(\text{mod 2)})
\end{array}
$}
\end{equation*}
($\mathbf{A}_{0}$ stands for a ``smooth chart'' on the threefold). But $%
\varphi :\widetilde{X}\rightarrow X_{f}$ \ is decomposed also globally into $%
m$ blow-ups {\small
\begin{equation*}
\begin{array}{r}
\widetilde{X}=\mathbf{Bl}_{\,q_{m}}(\mathbf{Bl}_{\,q_{m-1}}(\,\cdots \,(%
\mathbf{Bl}_{\,q_{1}}(X_{f}))))\overset{\pi _{m}}{\longrightarrow }\,\cdots
\,\overset{\pi _{3}}{\longrightarrow }\mathbf{Bl}_{\,q_{2}}(\mathbf{Bl}%
_{\,q_{1}}(X_{f}))\overset{\pi _{2}}{\longrightarrow }\mathbf{Bl}%
_{\,q_{1}}(X_{f}) \\
\,_{\pi _{1}=\pi }\downarrow \\
X_{f}
\end{array}
\,
\end{equation*}
} of $m$ points $q_{1}=\mathbf{0},\ q_{2}=\left( \mathbf{0},\left(
1:0:0:0\right) \right) ,\ldots ,q_{m}$, and is endowed with the ``separation
property''. By this we mean that, if $E_{1}=\mathcal{E}_{f},E_{2},\ldots
,E_{m}$ are the exceptional loci of $\pi _{1},\pi _{2},\ldots ,\pi _{m},$
respectively, then for $i\geq 2$ a singular point $q_{i}$ is resolved by $%
\pi _{i}$ and the (possibly existing) new singular point $q_{i+1}$ is
\textit{not} contained in the strict transforms of $E_{1},E_{2},\ldots
,E_{i-1}$ under $\pi _{i}$. Thus, defining $D_{i}$ to be the strict
transform of $E_{i}$ under $\pi _{i+1}\circ \pi _{i+2}\circ \cdots \circ \pi
_{m-1}\circ \pi _{m}$ on $\widetilde{X}$, we obtain an intersection graph of
the form:
\begin{figure}[h]
\begin{center}
\includegraphics[width=11cm,height=0.9cm]{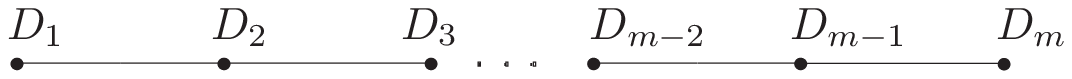} \vspace{0.3cm} \\[0pt]
\text{Case} $\mathbf{A}_{n}$.
\end{center}
\end{figure}

\noindent It is clear by \textbf{(i)} and \textbf{(ii)} that $D_{m}\cong%
\mathbb{P}_{\mathbb{C}}^{1}\times\mathbb{P}_{\mathbb{C}}^{1},$ for $%
n\equiv1( $mod $2$), and $D_{m}\cong\mathbb{P}_{\mathbb{C}}^{2},$ for $%
n\equiv0($mod $2 $), while $D_{j}\cong\mathbb{F}_{2}$ for all $j,\,1\leq
j\leq m-1.$ The Picard group Pic$\left( \mathbb{F}_{2}\right) \cong \mathbb{Z%
}^{2}$ of each $\mathbb{F}_{2}$ is generated by two projective lines: a
fiber $\mathsf{f}$ and a section $\mathsf{C}_{0}$ with $\mathsf{C}%
_{0}^{2}=-2.$ The intersection cycles read as follows:
\begin{equation*}
\left( D_{j}\cdot D_{j+1}\right) \left| _{D_{j}}\right. =\mathsf{C}_{0},\
\left( D_{j}\cdot D_{j+1}\right) \left| _{D_{j+1}}\right. \,\sim\,\mathsf{C}%
_{0}+2\,\mathsf{f},\ \ \forall j,\ 1\leq j\leq m-2\ ,
\end{equation*}
and {\small
\begin{equation*}
\left( D_{m-1}\cdot D_{m}\right) \left| _{D_{m-1}}\right. =\mathsf{C}_{0},\
\left( D_{m-1}\cdot D_{m}\right) \left| _{D_{m}}\right. \,\sim\,\left\{
\begin{array}{ll}
\mathsf{H}_{1}+\mathsf{H}_{2}, & \text{if\ \ }n\equiv1(\text{mod 2) } \\
\  & \, \\
2\,\mathsf{H}, & \text{if\ \ }n\equiv0(\text{mod 2)}
\end{array}
\right.
\end{equation*}
} where $\mathcal{O}_{\mathbb{P}_{\mathbb{C}}^{2}}(\mathsf{H})=\mathcal{O}_{%
\mathbb{P}_{\mathbb{C}}^{2}}(1)$ in Pic$(\mathbb{P}_{\mathbb{C}}^{2}),\ $and
\begin{equation*}
\mathcal{O}_{\mathbb{P}_{\mathbb{C}}^{1}\times\mathbb{P}_{\mathbb{C}}^{1}}(%
\mathsf{H}_{1})=\mathcal{O}_{\mathbb{P}_{\mathbb{C}}^{1}\times \mathbb{P}_{%
\mathbb{C}}^{1}}(1,0),\ \ \ \ \mathcal{O}_{\mathbb{P}_{\mathbb{C}}^{1}\times%
\mathbb{P}_{\mathbb{C}}^{1}}(\mathsf{H}_{2})=\mathcal{O}_{\mathbb{P}_{%
\mathbb{C}}^{1}\times\mathbb{P}_{\mathbb{C}}^{1}}(0,1)
\end{equation*}
in Pic$(\mathbb{P}_{\mathbb{C}}^{1}\times\mathbb{P}_{\mathbb{C}}^{1}).$ (We
shall keep the notation below whenever the arising exceptional prime
divisors are biregularly isomorphic to $\mathbb{F}_{2}$ or to $\mathbb{P}_{%
\mathbb{C}}^{1}\times\mathbb{P}_{\mathbb{C}}^{1}$). Obviously, $\left(
\mathsf{H}\cdot\mathsf{H}\right) \left| _{\mathbb{P}_{\mathbb{C}%
}^{2}}\right. =\left( \mathsf{H}_{1}\cdot\mathsf{H}_{2}\right) \left| _{%
\mathbb{P}_{\mathbb{C}}^{1}\times\mathbb{P}_{\mathbb{C}}^{1}}\right. =1$ and
$\left( \mathsf{H}_{1}\cdot\mathsf{H}_{1}\right) \left| _{\mathbb{P}_{%
\mathbb{C}}^{1}\times\mathbb{P}_{\mathbb{C}}^{1}}\right. =\left( \mathsf{H}%
_{2}\cdot\mathsf{H}_{2}\right) \left| _{\mathbb{P}_{\mathbb{C}}^{1}\times%
\mathbb{P}_{\mathbb{C}}^{1}}\right. =0.\bigskip$

\noindent$\blacktriangleright$ \textbf{Three characteristic rational
surfaces. }The remaining types $\mathbf{D}$-$\mathbf{E}$ of singularities $%
(X_{f},\mathbf{0})$ are more complicated as the $\varphi$'s under
construction will not fulfil the above ``separation property''. Furthermore,
since the exceptional locus after the first blow-up consists of two
irreducible components $\mathcal{E}_{f}^{\prime}$ and $\mathcal{E}%
_{f}^{\prime\prime},$ and the appearing new singular points ($3$ in case $%
\mathbf{D}_{4},2$ in case $\mathbf{D}_{n},\,n\geq5,$ and $1$ in cases $%
\mathbf{E}_{6},\mathbf{E}_{7},\mathbf{E}_{8}$) lie on the line $\mathcal{G}=%
\mathcal{E}_{f}^{\prime }\cap\mathcal{E}_{f}^{^{\prime\prime}},$ the strict
transforms of $\mathcal{G}$ together with their intersections with other
components (due to the next desingularization steps) will accompany us until
we arrive at $\widetilde{X}.$ In addition, to ensure a uniform resolution
procedure from the ``global'' point of view, one has to blow up the new
singularities \textit{simultaneously} (in each step) and take into account
the related intrinsic geometry. That's why, before proceeding to the
examination of the remaining cases, we define three rational compact complex
surfaces which will appear in a natural way as exceptional prime divisors of
our $\varphi$'s. (In fact, they will be inherited from the strict transforms
of the original $\mathcal{E}_{f}^{\prime}$ and $\mathcal{E}%
_{f}^{\prime\prime}$ as well as from the other intermediate components which
arise on one's way on the ``surface level''.)\medskip

\noindent$\triangleright$ Let $\mathbb{P}_{\mathbb{C}}^{2}\mathbf{[3]}$ be
the surface resulting after the blow-up $\mathbf{Bl}_{\{q_{0},q_{1},q_{2}\}}(%
\mathbb{P}_{\mathbb{C}}^{2})$ of $\mathbb{P}_{\mathbb{C}}^{2}$
simultaneously at three different points $q_{0},q_{1},q_{2}$ of a line $%
\mathcal{G}\subset\mathbb{P}_{\mathbb{C}}^{2}.$ (This surface is unique up
to biregular isomorphism, because for any other triple $q_{0}^{%
\prime},q_{1}^{\prime},q_{2}^{\prime}$ of different points of a line $%
\mathcal{G}^{\prime}\subset\mathbb{P}_{\mathbb{C}}^{2}$ the linear
isomorphism $\mathcal{G}\overset{\cong}{\rightarrow}\mathcal{G}^{\prime}$
mapping $q_{i}$ to $q_{i}^{\prime},i=1,2,3,$ can be extended to an
isomorphism $\mathbb{P}_{\mathbb{C}}^{2}\overset{\cong}{\rightarrow}\mathbb{P%
}_{\mathbb{C}}^{2}$). If we denote by $\mathsf{C}_{i}$ the inverse image of $%
q_{i}$ in $\mathbb{P}_{\mathbb{C}}^{2}\mathbf{[3],}$ then Pic$(\mathbb{P}_{%
\mathbb{C}}^{2}\mathbf{[3]})\cong\mathbb{Z}^{4}$ with $\{\mathsf{C}_{0},%
\mathsf{C}_{1},\mathsf{C}_{2},\mathsf{G}\}$ as generating system, where $%
\mathsf{G}$ is the strict transform of the original line $\mathcal{G}.$
Topologically $\{\mathsf{C}_{0},\mathsf{C}_{1},\mathsf{C}_{2},\mathsf{G}\}$
looks like:

\begin{figure}[h]
\begin{center}
\includegraphics[width=3.3in,height=1.2in]{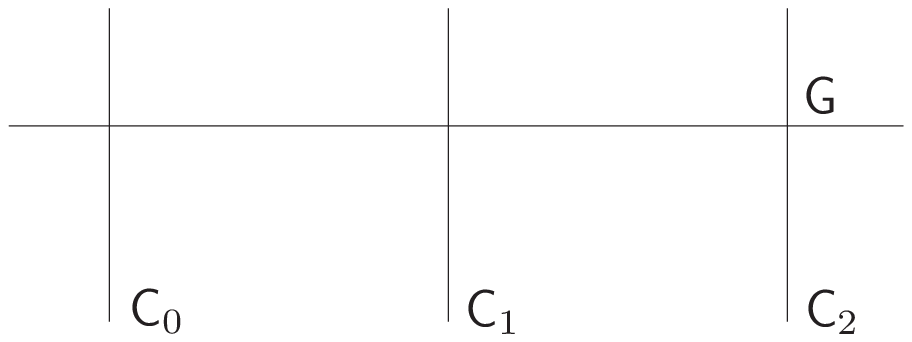}
\end{center}
\end{figure}

\noindent The intersection numbers of these generators on $\mathbb{P}_{%
\mathbb{C}}^{2}\mathbf{[3]}$ are the following:
\begin{equation*}
\left\{
\begin{array}{l}
\mathsf{C}_{0}^{2}=\mathsf{C}_{1}^{2}=\mathsf{C}_{2}^{2}=-1,\mathsf{G}%
^{2}=-2, \\
\left( \mathsf{G}\cdot\mathsf{C}_{0}\right) =\left( \mathsf{G}\cdot\mathsf{C}%
_{1}\right) =\left( \mathsf{G}\cdot\mathsf{C}_{2}\right) =1 \\
(\text{and }=0\text{ \ otherwise})
\end{array}
\right.
\end{equation*}

\noindent$\triangleright$ Let now $\mathbb{P}_{\mathbb{C}}^{2}\mathbf{[\,}%
\overline{\mathbf{3}}\,\mathbf{]}$ be the surface $\mathbf{Bl}_{\{q_{2}\}}(%
\mathbf{Bl}_{\{q_{0},q_{1}\}}(\mathbb{P}_{\mathbb{C}}^{2}))$ being
constructed by simultaneously blowing-up of $\mathbb{P}_{\mathbb{C}}^{2}$ at
two different points $q_{0},q_{1},$ followed by the blow-up at the
intersection point $q_{2}$ of the strict transform of $\overline{q_{0}\,q_{1}%
}$ and the blow-up of $q_{1}$ on $\mathbf{Bl}_{\{q_{0},q_{1}\}}(\mathbb{P}_{%
\mathbb{C}}^{2}).$ (The isomorphism type of $\mathbb{P}_{\mathbb{C}}^{2}%
\mathbf{[\,}\overline{\mathbf{3}}\,\mathbf{]}$ is unique, and one can use
arbitrary points $q_{0}\neq q_{1}$ for the construction). If we denote by $%
\mathsf{G}$ the strict transform of $\overline{q_{0}\,q_{1}},$ by $\mathsf{C}%
_{i}$ the strict transform of $q_{i},i\in\{0,1\},$ and by $\mathsf{C}_{2}$
the blow-up of $q_{2}$ within $\mathbb{P}_{\mathbb{C}}^{2}\mathbf{[\,}%
\overline{\mathbf{3}}\,\mathbf{],}$ then Pic$(\mathbb{P}_{\mathbb{C}}^{2}%
\mathbf{[\,}\overline{\mathbf{3}}\,\mathbf{]})\cong \mathbb{Z}^{4}$ with $\{%
\mathsf{C}_{0},\mathsf{C}_{1},\mathsf{C}_{2},\mathsf{G}\}$ as generating
system:

\begin{figure}[h]
\begin{center}
\includegraphics[width=6.2cm,height=2.95cm]{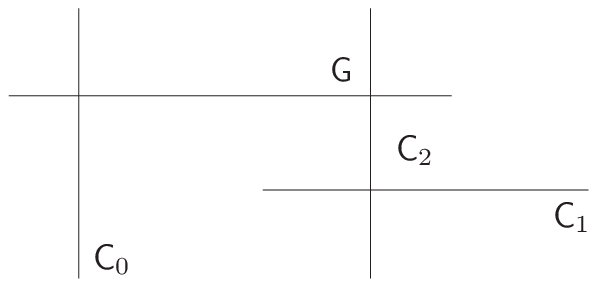}
\end{center}
\end{figure}

\noindent and intersection numbers:

\begin{equation*}
\left\{
\begin{array}{l}
\mathsf{C}_{0}^{2}=\mathsf{C}_{2}^{2}=-1,\mathsf{C}_{1}^{2}=\mathsf{G}%
^{2}=-2, \\
\left( \mathsf{G}\cdot\mathsf{C}_{0}\right) =\left( \mathsf{G}\cdot\mathsf{C}%
_{2}\right) =\left( \mathsf{C}_{1}\cdot\mathsf{C}_{2}\right) =1 \\
(\text{and }=0\text{ \ otherwise})
\end{array}
\right.
\end{equation*}

\noindent$\triangleright$ Finally, let $\mathbb{P}_{\mathbb{C}}^{2}\mathbf{%
[\,}\overline{\overline{\mathbf{3}}}\,\mathbf{]}$ denote the surface $%
\mathbf{Bl}_{\{q_{2}\}}(\mathbf{Bl}_{\{q_{1}\}}(\mathbf{Bl}_{\{q_{0}\}}(%
\mathbb{P}_{\mathbb{C}}^{2})))$ determined by blowing up a point $q_{0}$ of $%
\mathbb{P}_{\mathbb{C}}^{2},$ taking a line $\mathcal{G}\subset \mathbb{P}_{%
\mathbb{C}}^{2},$ with $q_{0}\in\mathcal{G},$ such that (strict transform of
$\mathcal{G}$)$\,\cap\,\mathbf{Bl}_{\{q_{0}\}}(\mathbb{P}_{\mathbb{C}%
}^{2})=\{q_{1}\},$ blowing up in turn $q_{1},$ and blowing up (at the last
step) $q_{2},$ where (strict transform of $\mathcal{G}$)$\,\cap\,\mathbf{Bl}%
_{\{q_{1}\}}(\mathbf{Bl}_{\{q_{0}\}}(\mathbb{P}_{\mathbb{C}%
}^{2}))=\{q_{2}\}. $ The isomorphism type of $\mathbb{P}_{\mathbb{C}}^{2}%
\mathbf{[\,}\overline{\overline{\mathbf{3}}}\,\mathbf{]}$ is again unique,
Pic$(\mathbb{P}_{\mathbb{C}}^{2}\mathbf{[\,}\overline {\overline{\mathbf{3}}%
}\,\mathbf{]})\cong\mathbb{Z}^{4}$ is generated by $\{\mathsf{C}_{0},\mathsf{%
C}_{1},\mathsf{C}_{2},\mathsf{G}\},$ where $\mathsf{G}$ is the (final)
strict transform of $\mathcal{G}$, $\mathsf{C}_{i}$ the strict transform of $%
q_{i},i\in\{0,1\},$ and $\mathsf{C}_{2}$ the blow-up of $q_{2}$ within $%
\mathbb{P}_{\mathbb{C}}^{2}\mathbf{[\,}\overline{\overline{\mathbf{3}}}\,%
\mathbf{].}$ Topologically $\{\mathsf{C}_{0},\mathsf{C}_{1},\mathsf{C}_{2},%
\mathsf{G}\}$ looks like:

\begin{figure}[h]
\begin{center}
\includegraphics[width=2.7069in,height=1.2773in]{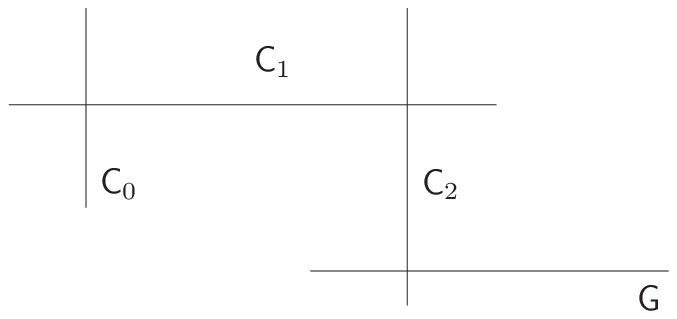}
\end{center}
\end{figure}

\noindent and the corresponding intersection numbers equal:
\begin{equation*}
\left\{
\begin{array}{l}
\mathsf{C}_{2}^{2}=-1,\mathsf{C}_{0}^{2}=\mathsf{C}_{1}^{2}=\mathsf{G}%
^{2}=-2, \\
\left( \mathsf{G}\cdot\mathsf{C}_{2}\right) =\left( \mathsf{C}_{0}\cdot%
\mathsf{C}_{1}\right) =\left( \mathsf{C}_{1}\cdot\mathsf{C}_{2}\right) =1 \\
(\text{and }=0\text{ \ otherwise})
\end{array}
\right.
\end{equation*}
\bigskip

\noindent \textbf{(iv) Types} $\mathbf{D}_{n}$\textbf{\ for }$n=2k,\ k\geq 2$%
. Let us first explain what happens in the $\mathbf{D}_{4}$-case. Blowing up
the origin $\mathbf{0}\in X_{f}$ we get {\small
\begin{equation*}
\mathbf{Bl}_{\,\mathbf{0}}(X_{f})=\left\{ \left( \left(
x_{1},..,x_{4}\right) ,\left( t_{1}:..:t_{4}\right) \right) \in \ \mathbf{Bl}%
_{\mathbf{0}}(\mathbb{C}^{4})\left|
x_{1}\,t_{1}^{2}+x_{1}\,t_{2}^{2}+t_{3}^{2}+t_{4}^{2}=0\right. \right\}
\end{equation*}
} with $\mathcal{E}_{f}=\mathcal{E}_{f}^{\prime }\cup \mathcal{E}%
_{f}^{\prime \prime }$ as exceptional locus. As we have already mentioned
above, $\mathbf{Bl}_{\,\mathbf{0}}(X_{f})$ possess the three $\mathbf{A}_{1}$%
-singularities {\small
\begin{equation*}
q_{0}=\left( \mathbf{0},\left( 0:1:0:0\right) \right) ,\,\text{\ }q_{1}=(%
\mathbf{0},(\sqrt{-1}:1:0:0)),\,\text{\ }q_{2}=(\mathbf{0},(-\sqrt{-1}%
:1:0:0)),
\end{equation*}
} which belong to the line $\mathcal{G}=\mathcal{E}_{f}^{\prime }\cap
\mathcal{E}_{f}^{^{\prime \prime }}.$ To obtain our global desingularization
$\varphi :\widetilde{X}\longrightarrow X_{f}$ it is enough to blow up once
more all three points $q_{0},q_{1},q_{2}$ simultaneously: {\small
\begin{equation*}
\widetilde{X}=\mathbf{Bl}_{\,\left\{ q_{0},q_{1},q_{2}\right\} }(\mathbf{Bl}%
_{\,\mathbf{0}}(X_{f}))\overset{\pi _{2}}{\longrightarrow }\mathbf{Bl}_{\,%
\mathbf{0}}(X_{f})\,\overset{\pi _{1}=\pi }{\longrightarrow }\,X_{f}\ .
\end{equation*}
} Let us denote by $D_{1}^{\prime }$ (resp. $D_{2}^{\prime \prime }$) the
strict transform of $\mathcal{E}_{f}^{\prime }$ (resp. $\mathcal{E}%
_{f}^{\prime \prime }$ $)$ under $\pi _{2}$, $D_{3}=\pi _{2}^{-1}\left(
q_{0}\right) ,$ $D_{j}=\pi _{2}^{-1}\left( q_{j}\right) ,$ for $j\in \left\{
1,2\right\} ,$ and define
\begin{equation*}
\mathsf{C}_{i}:=\pi _{2}^{-1}\left| _{D_{1}^{\prime }\text{ (resp. }%
D_{1}^{\prime \prime }\text{)}}\right. \ (q_{i}),\text{ \ }i\in \{0,1,2\}\ .
\end{equation*}
Then obviously $D_{1}\cong $ $D_{2}\cong D_{3}\cong \mathbb{P}_{\mathbb{C}%
}^{1}\times \mathbb{P}_{\mathbb{C}}^{1}$ and $D_{1}^{\prime }\cong
D_{1}^{\prime \prime }\cong \mathbb{P}_{\mathbb{C}}^{2}\mathbf{[3]}$ with
Picard group generated by $\mathsf{C}_{0},\mathsf{C}_{1},\mathsf{C}_{2}$ and
$\mathsf{G}$, where $\mathsf{G}$ is the strict transform of $\mathcal{G}$
under $\pi _{2}$. The intersection graph of these five exceptional divisors
is illustrated as follows:

\begin{figure}[h]
\begin{center}
\includegraphics[width=1.6328in,height=1.6501in]{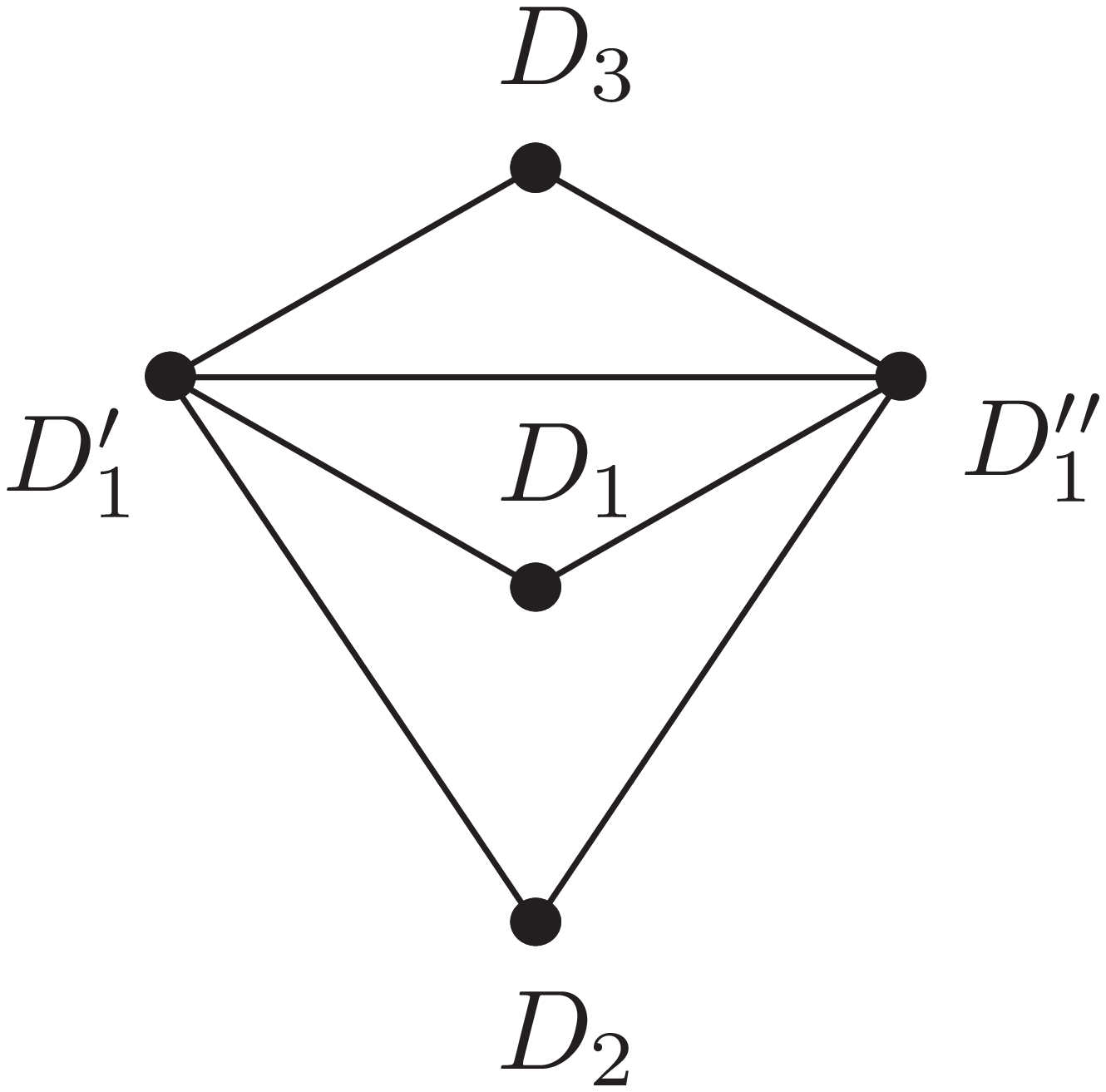}
\end{center}
\end{figure}
\noindent {}Generalizing to $\mathbf{D}_{2k\text{ }},$ the LR-diagram has
the form:
\begin{equation*}
\fbox{$
\begin{array}[t]{ccccccccccccc}
&  &  &  &  &  &  &  & \mathbf{A}_{0} &  &  &  &  \\
&  &  &  &  &  &  &  & \uparrow &  &  &  &  \\
&  &  &  &  &  &  &  & \mathbf{A}_{1} &  &  &  &  \\
&  &  &  &  &  &  &  & \uparrow &  &  &  &  \\
\mathbf{D}_{2k} & \rightarrow & \mathbf{D}_{2\left( k-1\right) } &
\rightarrow & \cdots & \rightarrow & \mathbf{D}_{6} & \rightarrow & \mathbf{D%
}_{4} & \rightarrow & \mathbf{A}_{1} & \rightarrow & \mathbf{A}_{0} \\
\downarrow &  & \downarrow &  &  &  & \downarrow &  & \downarrow &  &  &  &
\\
\mathbf{A}_{1} &  & \mathbf{A}_{1} &  &  &  & \mathbf{A}_{1} &  & \mathbf{A}%
_{1} &  &  &  &  \\
\downarrow &  & \downarrow &  &  &  & \downarrow &  & \downarrow &  &  &  &
\\
\mathbf{A}_{0} &  & \mathbf{A}_{0} &  &  &  & \mathbf{A}_{0} &  & \mathbf{A}%
_{0} &  &  &  &
\end{array}
$}
\end{equation*}
\medskip with a $D_{4}$ at its right-hand side and the intersection graph
looks like:
\begin{figure}[h]
\begin{center}
\includegraphics[width=10.9633cm,height=4.6327cm]{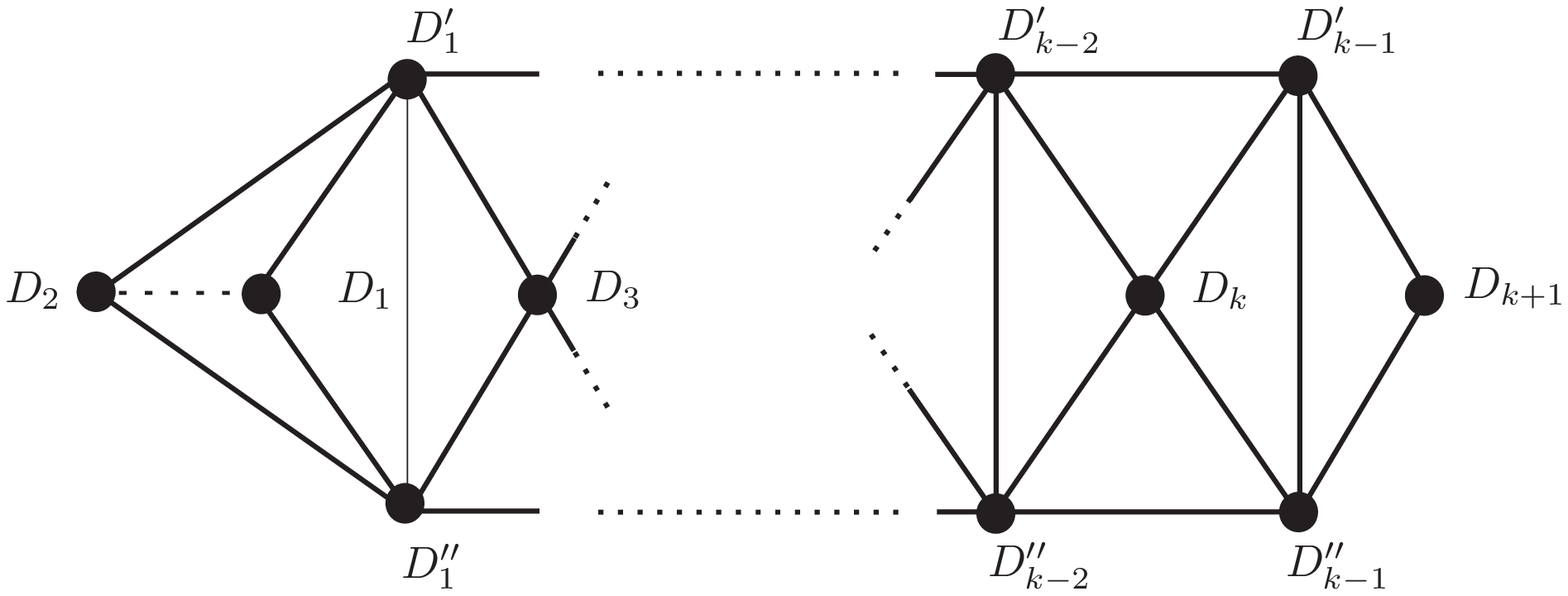} \vspace{0.3cm} \\%
[0pt]
\text{Case} $\mathbf{D}_{n}$ .
\end{center}
\end{figure}
\noindent {}

\noindent {}(The dotted line from $D_{2}$ to $D_{1}$ will be used only for
the case of odd $n$'s and it should be ignored for the time being). The
ordering of the subscripts of the divisors of the top and the bottom row is $%
1,2,...,k-2,k-1,$ whereas that of the divisors of the middle row is $%
2,1,3,4,...,k,k+1$. In this general case one needs altogether $k+1$ global ($%
=$ simultaneous) blow-ups to construct $\varphi :\widetilde{X}%
\longrightarrow X_{f}$. The exceptional prime divisors which occur are $%
D_{j}\cong \mathbb{P}_{\mathbb{C}}^{1}\times \mathbb{P}_{\mathbb{C}}^{1},\
\forall j,\ 1\leq j\leq k+1,$ and
\begin{equation*}
D_{1}^{\prime }\cong D_{1}^{\prime \prime }\cong \mathbb{P}_{\mathbb{C}}^{2}%
\mathbf{[3]},\ \ \ \ \ D_{j}^{\prime }\cong D_{j}^{\prime \prime }\cong
\mathbb{P}_{\mathbb{C}}^{2}\mathbf{[\,}\overline{\mathbf{3}}\,\mathbf{]},\ \
\ \forall j,\ \ 2\leq j\leq k-1\ ,
\end{equation*}
with the $k+1$ $\mathbb{P}_{\mathbb{C}}^{1}\times \mathbb{P}_{\mathbb{C}%
}^{1} $'s coming from the $\mathbf{A}_{1}$'s of the LR-diagram, and the $k-2$
pairs of $\mathbb{P}_{\mathbb{C}}^{2}\mathbf{[\,}\overline{\mathbf{3}}\,%
\mathbf{]}$'s inherited from the strict transforms of the $\mathcal{E}%
_{f}^{\prime }$ and $\mathcal{E}_{f}^{^{\prime \prime }}$ with respect to
the first $k-2$ global blow-ups (where in each step the singularities appear
pairwise). The corresponding intersection cycles are: {\footnotesize
\begin{equation*}
\begin{array}{ll}
\begin{array}{l}
\  \\
(D_{1}\cdot D_{1}^{\prime })\left| _{D_{1}}\right. =\mathsf{H}_{2},
\end{array}
&
\begin{array}{l}
\  \\
(D_{1}\cdot D_{1}^{\prime })\left| _{D_{1}^{\prime }}\right. =\mathsf{C}_{1},
\end{array}
\\
\begin{array}{l}
\  \\
(D_{1}\cdot D_{1}^{\prime \prime })\left| _{D_{1}}\right. =\mathsf{H}_{1},
\end{array}
&
\begin{array}{l}
\  \\
(D_{1}\cdot D_{1}^{\prime \prime })\left| _{D_{1}^{\prime \prime }}\right. =%
\mathsf{C}_{1},
\end{array}
\\
\begin{array}{l}
\  \\
(D_{2}\cdot D_{1}^{\prime })\left| _{D_{2}}\right. =\mathsf{H}_{2},
\end{array}
&
\begin{array}{l}
\  \\
(D_{2}\cdot D_{1}^{\prime })\left| _{D_{1}^{\prime }}\right. =\mathsf{C}_{2},
\end{array}
\\
\begin{array}{l}
\  \\
(D_{2}\cdot D_{1}^{\prime \prime })\left| _{D_{2}}\right. =\mathsf{H}_{1},
\end{array}
&
\begin{array}{l}
\  \\
(D_{2}\cdot D_{1}^{\prime \prime })\left| _{D_{1}^{\prime \prime }}\right. =%
\mathsf{C}_{2},
\end{array}
\\
\begin{array}{l}
\  \\
(D_{k+1}\cdot D_{k-1}^{\prime })\left| _{D_{k+1}}\right. =\mathsf{H}_{1},
\end{array}
&
\begin{array}{l}
\  \\
(D_{k+1}\cdot D_{k-1}^{\prime })\left| _{D_{k-1}^{\prime }}\right. =\mathsf{C%
}_{0},
\end{array}
\\
\begin{array}{l}
\  \\
(D_{k+1}\cdot D_{k-1}^{\prime \prime })\left| _{D_{k+1}}\right. =\mathsf{H}%
_{2}, \\
\,
\end{array}
&
\begin{array}{l}
\  \\
(D_{k+1}\cdot D_{k-1}^{\prime \prime })\left| _{D_{k-1}^{\prime \prime
}}\right. =\mathsf{C}_{0}, \\
\,
\end{array}
\end{array}
\end{equation*}
} \noindent while for $k\geq 3,$ and all $j$, $3\leq j\leq k,\medskip $
{\footnotesize
\begin{equation*}
\begin{array}{cc}
(D_{j}\cdot D_{j-1}^{\prime })\left| _{D_{j}}\right. =\mathsf{H}_{2}, &
(D_{j}\cdot D_{j-1}^{\prime })\left| _{D_{j-1}^{\prime }}\right. =\mathsf{C}%
_{2}, \\
\, & \, \\
(D_{j}\cdot D_{j-2}^{\prime })\left| _{D_{j}}\right. =\mathsf{H}_{1}, &
(D_{j}\cdot D_{j-2}^{\prime })\left| _{D_{j-2}^{\prime }}\right. =\mathsf{C}%
_{0},
\end{array}
\end{equation*}
\begin{equation*}
\begin{array}{cc}
(D_{j}\cdot D_{j-1}^{\prime \prime })\left| _{D_{j}}\right. =\mathsf{H}_{1},
& (D_{j}\cdot D_{j-1}^{\prime \prime })\left| _{D_{j-1}^{\prime \prime
}}\right. =\mathsf{C}_{2}, \\
\, & \, \\
(D_{j}\cdot D_{j-2}^{\prime \prime })\left| _{D_{j}}\right. =\mathsf{H}_{2},
& (D_{j}\cdot D_{j-2}^{\prime \prime })\left| _{D_{j-2}^{\prime \prime
}}\right. =\mathsf{C}_{0},
\end{array}
\end{equation*}
\begin{equation*}
\begin{array}{cc}
(D_{1}^{\prime }\cdot D_{2}^{\prime })\left| _{D_{1}^{\prime }}\right. \sim
\mathsf{G}+\mathsf{C}_{1}+\mathsf{C}_{2}, & (D_{1}^{\prime }\cdot
D_{2}^{\prime })\left| _{D_{2}^{\prime }}\right. =\mathsf{C}_{1}, \\
\, & \, \\
\noindent (D_{1}^{\prime \prime }\cdot D_{2}^{\prime \prime })\left|
_{D_{1}^{\prime \prime }}\right. \sim \mathsf{G}+\mathsf{C}_{1}+\mathsf{C}%
_{2}, & (D_{1}^{\prime \prime }\cdot D_{2}^{\prime \prime })\left|
_{D_{2}^{\prime \prime }}\right. =\mathsf{C}_{1},
\end{array}
\end{equation*}
} and for all $j$, $\ 2\leq j\leq k-2,$ {\small
\begin{equation*}
\begin{array}{ll}
(D_{j}^{\prime \,(\prime \prime )}\cdot D_{j+1}^{\prime \,(\prime \prime
)})\left| _{D_{j}^{\prime \,(\prime \prime )}}\right. \sim \mathsf{G}+%
\mathsf{C}_{1}+2\mathsf{C}_{2}, & (D_{j}^{\prime \,(\prime \prime )}\cdot
D_{j+1}^{\prime \,(\prime \prime )})\left| _{D_{j+1}^{\prime \,(\prime
\prime )}}\right. =\mathsf{C}_{1}. \\
\, & \,
\end{array}
\end{equation*}
} and finally, for all $j$, $\ 1\leq j\leq k-1,$ {\small
\begin{equation*}
\begin{array}{ll}
(D_{j}^{\prime }\cdot D_{j}^{\prime \prime })\left| _{D_{j}^{\prime
}}\right. =\mathsf{G},\ \ \  & (D_{j}^{\prime }\cdot D_{j}^{\prime \prime
})\left| _{D_{j}^{\prime \prime }}\right. =\mathsf{G}. \\
\, & \,
\end{array}
\end{equation*}
} \noindent \textbf{(v) Types} $\mathbf{D}_{n}$\textbf{\ for }$n=2k+1$. The
LR-diagram in this case reads as follows: {\small \newline
\begin{equation*}
\fbox{$
\begin{array}{ccccccccccccc}
\mathbf{D}_{2k+1} & \rightarrow & \mathbf{D}_{2(k-1)+1} & \rightarrow &
\cdots & \rightarrow & \mathbf{D}_{5} & \rightarrow & \mathbf{A}_{3} &
\rightarrow & \mathbf{A}_{1} & \rightarrow & \mathbf{A}_{0} \\
\downarrow &  & \downarrow &  &  &  & \downarrow &  &  &  &  &  &  \\
\mathbf{A}_{1} &  & \mathbf{A}_{1} &  &  &  & \mathbf{A}_{1} &  &  &  &  &
&  \\
\downarrow &  & \downarrow &  &  &  & \downarrow &  &  &  &  &  &  \\
\mathbf{A}_{0} &  & \mathbf{A}_{0} &  &  &  & \mathbf{A}_{0} &  &  &  &  &
&
\end{array}
$}
\end{equation*}
} Up to the introduction of the extra dotted edge into the game, the
intersection diagram remains the same, and the exceptional prime divisors
are
\begin{equation*}
D_{1}\cong \mathbb{F}_{2},\,\ \ D_{j}\cong \mathbb{P}_{\mathbb{C}}^{1}\times
\mathbb{P}_{\mathbb{C}}^{1},\ \ \forall j,\ \ 2\leq j\leq k+1,
\end{equation*}
and
\begin{equation*}
D_{j}^{\prime }\cong D_{j}^{\prime \prime }\cong \mathbb{P}_{\mathbb{C}}^{2}%
\mathbf{[\,}\overline{\mathbf{3}}\,\mathbf{]},\ \ \ \forall j,\ \ 1\leq
j\leq k-1\ .
\end{equation*}
Moreover, the intersection cycles are identical with those we have
encountered before in \textbf{(iv)}, up to the following ones:
{\footnotesize
\begin{equation*}
\begin{array}{ll}
(D_{1}\cdot D_{2})\left| _{D_{1}}\right. =\mathsf{C}_{0},\ \ \ (D_{1}\cdot
D_{2})\left| _{D_{2}}\right. \sim \,\mathsf{H}_{1}+\mathsf{H}_{2}, &
(D_{1}\cdot D_{1}^{\prime })\left| _{D_{1}}\right. =\mathsf{f},\ \ \
(D_{1}\cdot D_{1}^{\prime })\left| _{D_{1}^{\prime }}\right. =\mathsf{C}_{1},
\\
\  & \  \\
(D_{1}\cdot D_{1}^{\prime \prime })\left| _{D_{1}}\right. =\mathsf{f}%
^{\prime },\ \ \ \ \ \ (D_{1}\cdot D_{1}^{\prime \prime })\left|
_{D_{1}^{\prime \prime }}\right. =\mathsf{C}_{1}, & \text{(}\mathsf{f}%
\underset{\text{set th.}}{\neq }\mathsf{f}^{\prime }\text{ fibers of }%
\mathbb{F}_{2}\text{)} \\
\  & \  \\
(D_{1}^{\prime \,(\prime \prime )}\cdot D_{2}^{\prime \,(\prime \prime
)})\left| _{D_{1}^{\prime \,(\prime \prime )}}\right. \sim \mathsf{G}+%
\mathsf{C}_{1}+2\mathsf{C}_{2}, & (D_{1}^{\prime \,(\prime \prime )}\cdot
D_{2}^{\prime \,(\prime \prime )})\left| _{D_{2}^{\prime \,(\prime \prime
)}}\right. =\mathsf{C}_{1}, \\
\, & \,
\end{array}
\end{equation*}
} \noindent \textbf{(vi) Type} $\mathbf{E}_{6}$\textbf{. }The LR-diagram in
this case reads as:
\begin{equation*}
\fbox{$
\begin{array}{c}
\mathbf{\ } \\
\mathbf{E}_{6}\rightarrow \mathbf{A}_{5}\rightarrow \mathbf{A}%
_{3}\rightarrow \mathbf{A}_{1}\rightarrow \mathbf{A}_{0}\  \\
\mathbf{\ }
\end{array}
$}
\end{equation*}
Globally, the desingularization procedure is described as follows. To obtain
the morphism $\varphi :\widetilde{X}\longrightarrow X_{f},$ we need $3$
additional blow-ups at three points $q_{0},q_{1},q_{2}$ \textit{after} $%
\mathbf{Bl}_{\,\mathbf{0}}(X_{f})\,\overset{\pi }{\longrightarrow }\,X_{f},$
i.e.,
\begin{equation*}
\begin{array}{l}
\mathbf{Bl}_{\,q_{1}}(\mathbf{Bl}_{\,q_{0}}(\mathbf{Bl}_{\,\mathbf{0}%
}(X_{f})))\overset{\pi _{2}}{\longrightarrow }\mathbf{Bl}_{\,q_{0}}(\mathbf{%
Bl}_{\,\mathbf{0}}(X_{f}))\overset{\pi _{1}}{\longrightarrow }\mathbf{Bl}_{\,%
\mathbf{0}}(X_{f})\overset{\pi _{0}=\pi }{\longrightarrow }X_{f} \\
\  \\
\uparrow _{\pi _{3}} \\
\  \\
\widetilde{X}=\mathbf{Bl}_{\,q_{2}}(\mathbf{Bl}_{\,q_{1}}(\mathbf{Bl}%
_{\,q_{0}}(\mathbf{Bl}_{\,\mathbf{0}}(X_{f}))))
\end{array}
\end{equation*}
where$\ q_{0}=\left( \mathbf{0},\left( 0:1:0:0\right) \right) \in U_{2}$ on%
{\small
\begin{equation*}
\mathbf{Bl}_{\,\mathbf{0}}(X_{f})=\left\{ \left( \left(
x_{1},..,x_{4}\right) ,\left( t_{1}:t_{2}:t_{3}:t_{4}\right) \right) \in \
\mathbf{Bl}_{\mathbf{0}}(\mathbb{C}^{4})\ \ \left| \ \right.
x_{1}\,t_{1}^{2}+x_{2}^{2}\,t_{2}^{2}+t_{3}^{2}+t_{4}^{2}=0\right\} .
\end{equation*}
} Analogously, one gets $q_{1}=\left( \mathbf{0},\left( 0:1:0:0\right)
\right) $ on $\mathbf{Bl}_{\,q_{0}}(\mathbf{Bl}_{\,\mathbf{0}}(X_{f})\left|
_{U_{2}}\right. ),$ which equals {\small
\begin{equation*}
\left\{ \left( \left( y_{2,1},..,y_{2,4}\right) ,\left( \lambda
_{1}:..:\lambda _{4}\right) \right) \in \ U_{2}\times \mathbb{P}_{\mathbb{C}%
}^{3}\ \ \left| \ \right. \left( y_{2,1}\right) ^{2}\,\lambda _{1}\lambda
_{2}+\lambda _{2}^{2}+\lambda _{3}^{2}+\lambda _{4}^{2}=0\right\}
\end{equation*}
} (and similarly for $q_{2}\in \mathbf{Bl}_{\,q_{1}}(\mathbf{Bl}_{\,q_{0}}(%
\mathbf{Bl}_{\,\mathbf{0}}(X_{f})\left| _{U_{2}}\right. ))$ in the last
step). The point $q_{0}$ belongs to the line $\mathcal{G}\mathsf{\,}=\,%
\mathcal{E}_{f}^{\prime }\cap \mathcal{E}_{f}^{\prime \prime }$ (where, as
usual, $\pi ^{-1}(\mathbf{0})=\mathcal{E}_{f}^{\prime }\cup \mathcal{E}%
_{f}^{\prime \prime }$) and $\left( \mathbf{Bl}_{\,\mathbf{0}%
}(X_{f}),q_{0}\right) $ is an $\mathbf{A}_{5}$-singularity. According to
\textbf{(iii),} this will be resolved by $\pi _{1}\circ \pi _{2}\circ \pi
_{3}$ to give two $\mathbb{F}_{2}$'s and one $\mathbb{P}_{\mathbb{C}%
}^{1}\times \mathbb{P}_{\mathbb{C}}^{1}$ as exceptional divisors. More
precisely,
\begin{equation*}
q_{1}\in \left( \text{strict transform of\thinspace\ }\mathcal{G}\text{
\thinspace under\thinspace\ }\pi _{1}\right) \cap \left( \text{exceptional
locus of \thinspace }\pi _{1}\right)
\end{equation*}
is the new, $\mathbf{A}_{3}$-singularity, while
\begin{equation*}
q_{2}\in \left( \text{strict transform of\thinspace\ }\mathcal{G}\text{
\thinspace under\thinspace\ }\pi _{1}\circ \pi _{2}\right) \mathbb{r}\left(
\begin{array}{c}
\text{strict transform } \\
\text{of the exceptional} \\
\text{ locus of \thinspace }\pi _{1}\text{ under }\pi _{2}
\end{array}
\right)
\end{equation*}
is the final $\mathbf{A}_{1}$-singularity. Let us denote by $D_{1}$ the
strict transform of the exceptional locus of $\pi _{1}$ under $\pi _{2}\circ
\pi _{3},$ by $D_{2}$ the strict transform of the exceptional locus of $\pi
_{2}$ under $\pi _{3},$ by $D_{3}$ the exceptional locus of $\pi _{3},$ and
finally by $D_{4}$ (resp. $D_{4}^{\prime },\mathsf{G}$) the strict transform
of the original $\mathcal{E}_{f}^{\prime }$ (resp. $\mathcal{E}_{f}^{\prime
\prime },\mathcal{G}$) under $\pi _{1}\circ \pi _{2}\circ \pi _{3},$ and
define
\begin{equation*}
\left\{
\begin{array}{l}
\mathsf{C}_{0}:=\left( \text{strict transform of\thinspace\ }q_{0}\text{
\thinspace under\thinspace\ }\pi _{1}\circ \pi _{2}\circ \pi _{3}\text{ on }%
D_{4}\ \ (\text{resp. }D_{4}^{\prime })\right) \smallskip \\
\mathsf{C}_{1}:=\left( \text{strict transform of\thinspace\ }q_{1}\text{
\thinspace under\thinspace\ }\pi _{2}\circ \pi _{3}\text{ on }D_{4}\ \ (%
\text{resp. }D_{4}^{\prime })\right) \smallskip \\
\mathsf{C}_{2}:=\left( \text{the blow-up of\thinspace\ }q_{2}\text{
\thinspace by\thinspace\ }\pi _{3}\text{ on }D_{4}\ \ (\text{resp. }%
D_{4}^{\prime })\right) \ .
\end{array}
\right.
\end{equation*}
Then
\begin{equation*}
D_{1}\cong D_{2}\cong \mathbb{F}_{2},\ \ \ D_{3}\cong \mathbb{P}_{\mathbb{C}%
}^{1}\times \mathbb{P}_{\mathbb{C}}^{1},\ \ \ D_{4}\cong D_{4}^{\prime
}\cong \mathbb{P}_{\mathbb{C}}^{2}\mathbf{[\,}\overline{\overline{\mathbf{3}}%
}\,\mathbf{],}
\end{equation*}
with Pic($D_{4}$) (resp. Pic($D_{4}^{\prime }$)$)$ generated by $\mathsf{C}%
_{0},\mathsf{C}_{1},\mathsf{C}_{2},\mathsf{G},$ intersection graph

\begin{figure}[h]
\begin{center}
\includegraphics[width=3.6cm,height=4.6cm]{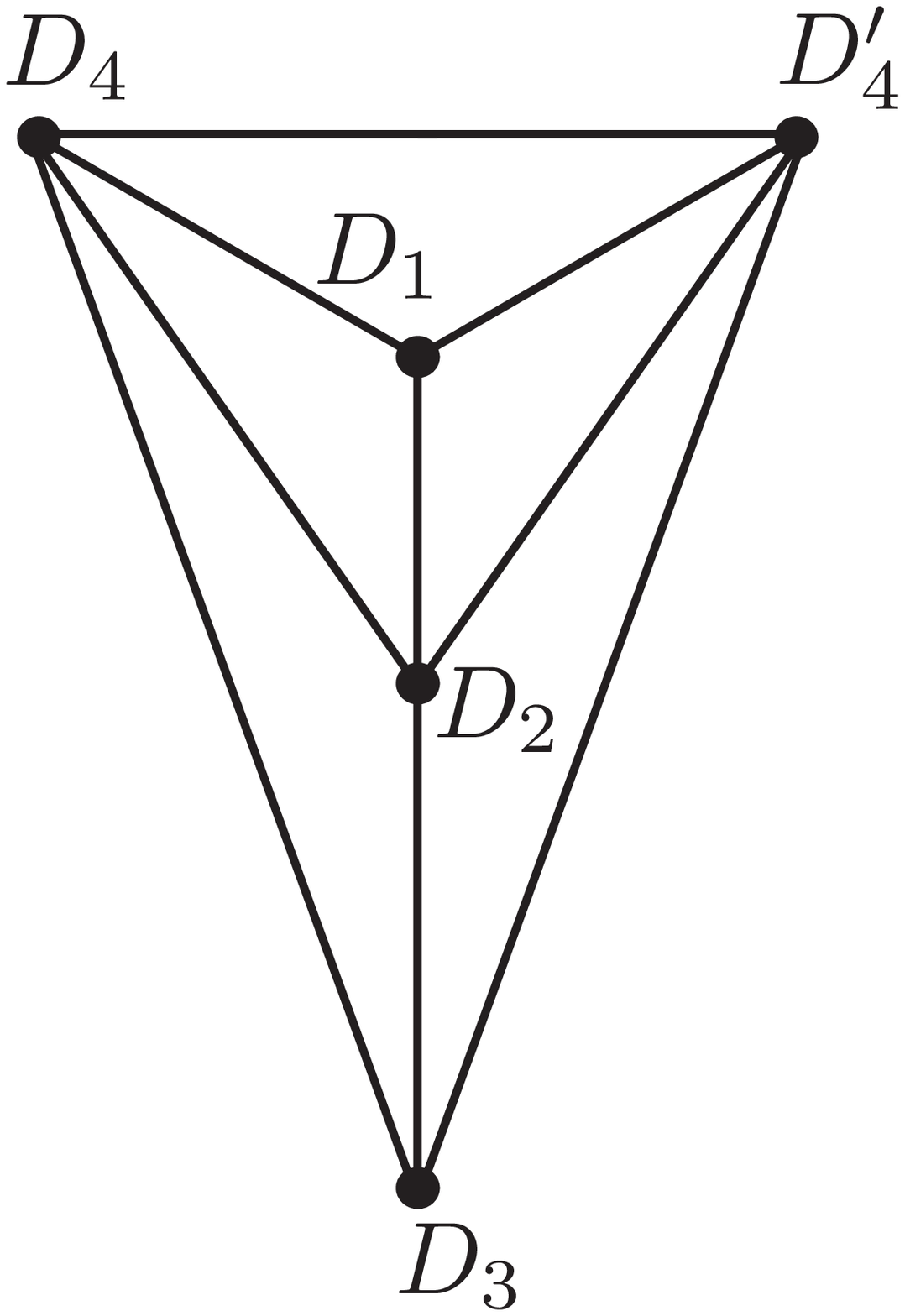} \vspace{0.3cm} \\[0pt]
\text{Case} $\mathbf{E}_{6}$.
\end{center}
\end{figure}

\noindent {}and intersections cycles: {\footnotesize
\begin{equation*}
\begin{array}{cc}
(D_{1}\cdot D_{2})\left| _{D_{1}}\right. =\mathsf{C}_{0}, & (D_{1}\cdot
D_{2})\left| _{D_{2}}\right. \sim \,\mathsf{C}_{0}+2\mathsf{f}, \\
\  & \  \\
(D_{1}\cdot D_{4})\left| _{D_{1}}\right. =\mathsf{f}, & (D_{1}\cdot
D_{4})\left| _{D_{4}}\right. =\mathsf{C}_{0},
\end{array}
\end{equation*}
\begin{equation*}
\begin{array}{cc}
(D_{1}\cdot D_{4}^{\prime })\left| _{D_{1}}\right. =\mathsf{f}^{\prime }, &
(D_{1}\cdot D_{4}^{\prime })\left| _{D_{4}^{\prime }}\right. =\mathsf{C}_{0},
\\
\  & \  \\
(D_{2}\cdot D_{3})\left| _{D_{2}}\right. =\mathsf{C}_{0}, & (D_{2}\cdot
D_{3})\left| _{D_{3}}\right. \sim \mathsf{H}_{1}+\mathsf{H}_{2},
\end{array}
\end{equation*}
\begin{equation*}
\begin{array}{cc}
(D_{2}\cdot D_{4})\left| _{D_{2}}\right. =\mathsf{f}, & (D_{2}\cdot
D_{4})\left| _{D_{4}}\right. =\mathsf{C}_{1}, \\
\  & \  \\
(D_{2}\cdot D_{4}^{\prime })\left| _{D_{2}}\right. =\mathsf{f}^{\prime }, &
(D_{2}\cdot D_{4}^{\prime })\left| _{D_{4}^{\prime }}\right. =\mathsf{C}_{1},
\end{array}
\end{equation*}
\begin{equation*}
\begin{array}{cc}
(D_{3}\cdot D_{4})\left| _{D_{3}}\right. =\mathsf{H}_{1}, & (D_{3}\cdot
D_{4})\left| _{D_{4}}\right. =\mathsf{C}_{2}, \\
\  & \  \\
(D_{3}\cdot D_{4}^{\prime })\left| _{D_{3}}\right. =\mathsf{H}_{2}, &
(D_{3}\cdot D_{4}^{\prime })\left| _{D_{4}^{\prime }}\right. =\mathsf{C}_{2},
\end{array}
\end{equation*}
\begin{equation*}
\begin{array}{cc}
(D_{4}\cdot D_{4}^{\prime })\left| _{D_{4}}\right. =\mathsf{G}, &
(D_{4}\cdot D_{4}^{\prime })\left| _{D_{4}^{\prime }}\right. =\mathsf{G}.
\end{array}
\end{equation*}
} (where $\mathsf{f}\underset{\text{set th.}}{\neq }\mathsf{f}^{\prime }$
fibers of $\mathbb{F}_{2}$).\medskip

\noindent\textbf{(vii) The cases} $\mathbf{E}_{7}$ \textbf{and} $\mathbf{E}%
_{8}$\textbf{. }Since $\mathbf{E}_{8}$ passes to an $\mathbf{E}_{7}$ after
the first blow-up, the LR-diagram looks like: {\small
\begin{equation*}
\fbox{$
\begin{array}[t]{ccccccccccc}
&  &  &  &  &  & \mathbf{A}_{0} &  &  &  &  \\
&  &  &  &  &  & \uparrow &  &  &  &  \\
&  &  &  &  &  & \mathbf{A}_{1} &  &  &  &  \\
&  &  &  &  &  & \uparrow &  &  &  &  \\
\mathbf{E}_{8} & \dashrightarrow & \mathbf{E}_{7} & \longrightarrow &
\mathbf{D}_{6} & \longrightarrow & \mathbf{D}_{4} & \longrightarrow &
\mathbf{A}_{1} & \longrightarrow & \mathbf{A}_{0} \\
&  &  &  & \downarrow &  & \downarrow &  &  &  &  \\
&  &  &  & \mathbf{A}_{1} &  & \mathbf{A}_{1} &  &  &  &  \\
&  &  &  & \downarrow &  & \downarrow &  &  &  &  \\
&  &  &  & \mathbf{A}_{0} &  & \mathbf{A}_{0} &  &  &  &
\end{array}
$}
\end{equation*}
} \smallskip\noindent{}Globally, for the resolution of $\mathbf{E}_{7}$-
(resp. $\mathbf{E}_{8}$-) singularity, we need $4$ (resp. $5$) blow-ups. The
intersection graph contains $10$ (resp. $12$) vertices (with the dotted
edges only in the $\mathbf{E}_{8}$-case)

\begin{figure}[h]
\begin{center}
\includegraphics[width=5.3312cm,height=7.2cm]{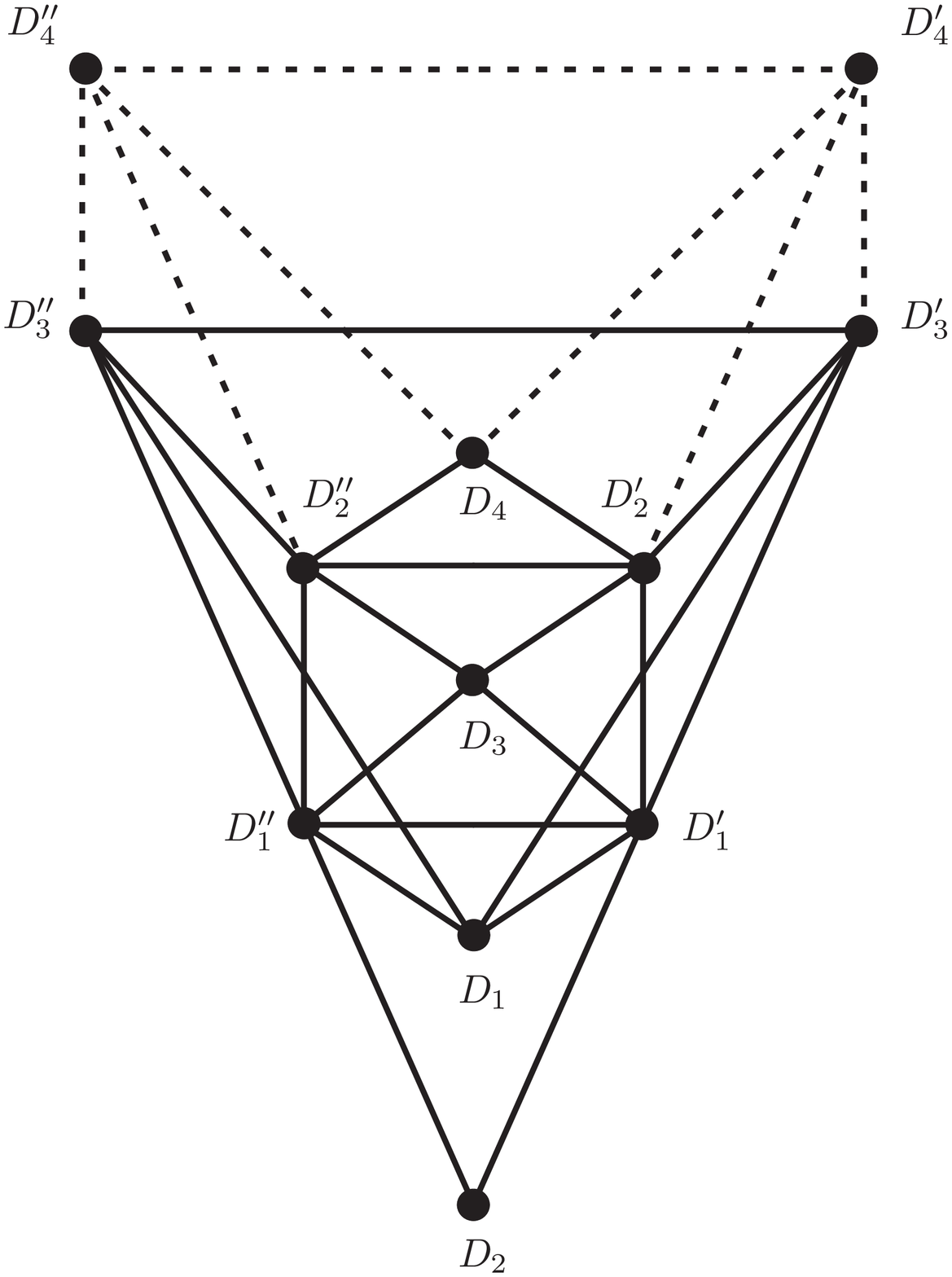} \vspace{0.3cm} \\[0pt%
]
\text{Cases }$\mathbf{E}_{7}$\text{ and }$\mathbf{E}_{8}$.
\end{center}
\end{figure}

\noindent corresponding to the $12$ exceptional prime divisors
\begin{equation*}
\begin{array}{l}
D_{1}\cong D_{2}\cong D_{3}\cong D_{4}\cong\mathbb{P}_{\mathbb{C}}^{1}\times%
\mathbb{P}_{\mathbb{C}}^{1}, \\
\, \\
D_{1}^{\prime}\cong D_{1}^{\prime\prime}\cong\mathbb{P}_{\mathbb{C}}^{2}%
\mathbf{[3],\ \ }D_{2}^{\prime}\cong D_{2}^{\prime\prime}\cong \mathbb{P}_{%
\mathbb{C}}^{2}\mathbf{[\,}\overline{\mathbf{3}}\,\mathbf{],} \\
\, \\
D_{3}^{\prime}\cong D_{3}^{\prime\prime}\cong D_{4}^{\prime}\cong
D_{4}^{\prime\prime}\cong\mathbb{P}_{\mathbb{C}}^{2}\mathbf{[\,}\overline{%
\overline{\mathbf{3}}}\,\mathbf{].}
\end{array}
\end{equation*}
The ``central'' four $\mathbb{P}_{\mathbb{C}}^{1}\times\mathbb{P}_{\mathbb{C}%
}^{1}$'s come from the four lastly appearing $\mathbf{A}_{1}$'s, and the
four top $\mathbb{P}_{\mathbb{C}}^{2}\mathbf{[\,}\overline{\overline{\mathbf{%
3}}}\,\mathbf{]}$'s are due to the last three sucessive blow-ups of $%
\mathcal{E}_{f}^{\prime}$ $\ $and $\mathcal{E}_{f}^{\prime\prime}.$ The two $%
\mathbb{P}_{\mathbb{C}}^{2}\mathbf{[3]}$'s (resp. the two $\mathbb{P}_{%
\mathbb{C}}^{2}\mathbf{[\,}\overline{\mathbf{3}}\,\mathbf{]}$'s) are in turn
inherited from the strict transforms of $\mathcal{E}_{f}^{\prime}$ $\ $and $%
\mathcal{E}_{f}^{\prime\prime}$ after passing from $\mathbf{D}_{4}$ to the
three $\mathbf{A}_{1}$'s (resp. from $\mathbf{D}_{6}$ to $\mathbf{D}_{4}$).
Making use of the previously introduced notation, the intersection cycles
read as follows:{\footnotesize
\begin{equation*}
\begin{array}{cc}
(D_{1}\cdot D_{1}^{\prime})\left| _{D_{1}}\right. =\mathsf{H}_{1}, &
(D_{1}\cdot D_{1}^{\prime})\left| _{D_{1}^{\prime}}\right. =\,\mathsf{C}_{1},
\\
\  & \  \\
(D_{1}\cdot D_{1}^{\prime\prime})\left| _{D_{1}}\right. =\mathsf{H}_{2}, &
(D_{1}\cdot D_{1}^{\prime\prime})\left| _{D_{1}^{\prime\prime}}\right. =%
\mathsf{C}_{2},
\end{array}
\end{equation*}
\begin{equation*}
\begin{array}{cc}
(D_{1}\cdot D_{3}^{\prime})\left| _{D_{1}}\right. =\mathsf{H}_{2}, &
(D_{1}\cdot D_{3}^{\prime})\left| _{D_{3}^{\prime}}\right. =\,\mathsf{C}_{2},
\\
\  & \  \\
(D_{1}\cdot D_{3}^{\prime\prime})\left| _{D_{1}}\right. =\mathsf{H}_{1}, &
(D_{1}\cdot D_{3}^{\prime\prime})\left| _{D_{3}^{\prime\prime}}\right. =%
\mathsf{C}_{2},
\end{array}
\end{equation*}
\begin{equation*}
\begin{array}{cc}
(D_{1}^{\prime}\cdot D_{1}^{\prime\prime})\left| _{D_{1}^{\prime}}\right. =%
\mathsf{G}, & (D_{1}^{\prime}\cdot D_{1}^{\prime\prime})\left|
_{D_{1}^{\prime\prime}}\right. =\mathsf{G}, \\
\  & \  \\
(D_{1}^{\prime}\cdot D_{2})\left| _{D_{1}^{\prime}}\right. =\mathsf{C}_{2},
& (D_{1}^{\prime}\cdot D_{2})\left| _{D_{2}}\right. =\mathsf{H}_{1},
\end{array}
\end{equation*}
\begin{equation*}
\begin{array}{cc}
(D_{1}^{\prime}\cdot D_{2}^{\prime})\left| _{D_{1}^{\prime}}\right. \sim%
\mathsf{G}+\mathsf{C}_{1}+2\mathsf{C}_{2}, & (D_{1}^{\prime}\cdot
D_{2}^{\prime})\left| _{D_{2}^{\prime}}\right. =\mathsf{C}_{1}, \\
\  & \  \\
(D_{1}^{\prime}\cdot D_{3})\left| _{D_{1}^{\prime}}\right. =\mathsf{C}_{0},
& (D_{1}^{\prime}\cdot D_{3})\left| _{D_{3}}\right. =\mathsf{H}_{2},
\end{array}
\end{equation*}
\begin{equation*}
\begin{array}{cc}
(D_{1}^{\prime}\cdot D_{3}^{\prime})\left| _{D_{1}^{\prime}}\right. \sim%
\mathsf{G}+\mathsf{C}_{0}+\mathsf{C}_{2}, & (D_{1}^{\prime}\cdot
D_{3}^{\prime})\left| _{D_{3}^{\prime}}\right. =\mathsf{C}_{1}, \\
\  & \  \\
(D_{1}^{\prime\prime}\cdot D_{2})\left| _{D_{1}^{\prime\prime}}\right. =%
\mathsf{C}_{2}, & (D_{1}^{\prime\prime}\cdot D_{2})\left| _{D_{2}}\right. =%
\mathsf{H}_{2},
\end{array}
\end{equation*}
\begin{equation*}
\begin{array}{cc}
(D_{1}^{\prime\prime}\cdot D_{2}^{\prime\prime})\left| _{D_{1}^{\prime\prime
}}\right. \sim\mathsf{G}+\mathsf{C}_{1}+2\mathsf{C}_{2}, &
(D_{1}^{\prime\prime}\cdot D_{2}^{\prime\prime})\left|
_{D_{2}^{\prime\prime}}\right. =\mathsf{C}_{1}, \\
\  & \  \\
(D_{1}^{\prime\prime}\cdot D_{3})\left| _{D_{1}^{\prime\prime}}\right. =%
\mathsf{C}_{0}, & (D_{1}^{\prime\prime}\cdot D_{3})\left| _{D_{3}}\right. =%
\mathsf{H}_{1},
\end{array}
\end{equation*}
\begin{equation*}
\begin{array}{cc}
(D_{1}^{\prime\prime}\cdot D_{3}^{\prime\prime})\left| _{D_{1}^{\prime\prime
}}\right. \sim\mathsf{G}+\mathsf{C}_{0}+\mathsf{C}_{2}, &
(D_{1}^{\prime\prime}\cdot D_{3}^{\prime\prime})\left|
_{D_{3}^{\prime\prime}}\right. =\mathsf{C}_{1}, \\
\  & \  \\
(D_{2}^{\prime}\cdot D_{2}^{\prime\prime})\left| _{D_{2}^{\prime}}\right. =%
\mathsf{G}, & (D_{2}^{\prime}\cdot D_{2}^{\prime\prime})\left|
_{D_{2}^{\prime\prime}}\right. =\mathsf{G},
\end{array}
\end{equation*}
\begin{equation*}
\begin{array}{cc}
(D_{2}^{\prime}\cdot D_{3}^{\prime})\left| _{D_{2}^{\prime}}\right. \sim%
\mathsf{G}+\mathsf{C}_{0}+\mathsf{C}_{2}, & (D_{2}^{\prime}\cdot
D_{3}^{\prime})\left| _{D_{3}^{\prime}}\right. =\mathsf{C}_{0}, \\
\  & \  \\
(D_{2}^{\prime}\cdot D_{3})\left| _{D_{2}^{\prime}}\right. =\mathsf{C}_{2},
& (D_{2}^{\prime}\cdot D_{3})\left| _{D_{3}}\right. =\mathsf{H}_{1},
\end{array}
\end{equation*}
\begin{equation*}
\begin{array}{cc}
(D_{2}^{\prime\prime}\cdot D_{3}^{\prime\prime})\left| _{D_{2}^{\prime\prime
}}\right. \sim\mathsf{G}+\mathsf{C}_{0}+\mathsf{C}_{2}, &
(D_{2}^{\prime\prime}\cdot D_{3}^{\prime\prime})\left|
_{D_{3}^{\prime\prime}}\right. =\mathsf{C}_{0}, \\
\  & \  \\
(D_{2}^{\prime\prime}\cdot D_{3})\left| _{D_{2}^{\prime\prime}}\right. =%
\mathsf{C}_{2}, & (D_{2}^{\prime\prime}\cdot D_{3})\left| _{D_{3}}\right. =%
\mathsf{H}_{1},
\end{array}
\end{equation*}
\begin{equation*}
\begin{array}{cc}
(D_{2}^{\prime}\cdot D_{4})\left| _{D_{2}^{\prime}}\right. =\mathsf{H}_{1},
& (D_{2}^{\prime}\cdot D_{4})\left| _{D_{4}}\right. =\mathsf{C}_{0}, \\
\  & \  \\
(D_{2}^{\prime\prime}\cdot D_{4})\left| _{D_{2}^{\prime\prime}}\right. =%
\mathsf{C}_{0}, & (D_{2}^{\prime\prime}\cdot D_{4})\left| _{D_{4}}\right. =%
\mathsf{H}_{1}, \\
\  & \
\end{array}
\end{equation*}
} with {\small $(D_{3}^{\prime}\cdot D_{3}^{\prime\prime})\left|
_{D_{3}^{\prime}}\right. =\mathsf{G},\ (D_{3}^{\prime}\cdot
D_{3}^{\prime\prime})\left| _{D_{3}^{\prime\prime}}\right. =\mathsf{G},$ } \
and {\footnotesize
\begin{equation*}
\begin{array}{cc}
(D_{2}^{\prime}\cdot D_{4}^{\prime})\left| _{D_{2}^{\prime}}\right. \sim%
\mathsf{G}+\mathsf{C}_{1}+2\mathsf{C}_{2}, & (D_{2}^{\prime}\cdot
D_{4}^{\prime})\left| _{D_{4}^{\prime}}\right. =\mathsf{C}_{1}, \\
\  & \  \\
(D_{4}\cdot D_{4}^{\prime\prime})\left| _{D_{4}}\right. =\mathsf{H}_{1}, &
(D_{4}\cdot D_{4}^{\prime\prime})\left| _{D_{4}^{\prime\prime}}\right. =%
\mathsf{C}_{2},
\end{array}
\end{equation*}
\begin{equation*}
\begin{array}{cc}
(D_{2}^{\prime\prime}\cdot D_{4}^{\prime\prime})\left| _{D_{2}^{\prime\prime
}}\right. \sim\mathsf{G}+\mathsf{C}_{1}+2\mathsf{C}_{2}, &
(D_{2}^{\prime\prime}\cdot D_{4}^{\prime\prime})\left|
_{D_{4}^{\prime\prime}}\right. =\mathsf{C}_{1}, \\
\  & \  \\
(D_{4}\cdot D_{4}^{\prime})\left| _{D_{4}}\right. =\mathsf{H}_{2}, &
(D_{4}\cdot D_{4}^{\prime})\left| _{D_{4}^{\prime}}\right. =\,\mathsf{C}_{2},
\end{array}
\end{equation*}
\begin{equation*}
\begin{array}{cc}
(D_{3}^{\prime}\cdot D_{4}^{\prime})\left| _{D_{3}^{\prime}}\right. \sim%
\mathsf{G}+\mathsf{C}_{1}+2\mathsf{C}_{2}, & (D_{3}^{\prime}\cdot
D_{4}^{\prime})\left| _{D_{4}^{\prime}}\right. =\mathsf{C}_{0}, \\
\  & \  \\
(D_{4}^{\prime}\cdot D_{4}^{\prime\prime})\left| _{D_{4}^{\prime}}\right. =%
\mathsf{G}, & (D_{4}^{\prime}\cdot D_{4}^{\prime\prime})\left|
_{D_{4}^{\prime\prime}}\right. =\mathsf{G},
\end{array}
\end{equation*}
\begin{equation*}
\begin{array}{cc}
(D_{3}^{\prime\prime}\cdot D_{4}^{\prime\prime})\left| _{D_{3}^{\prime\prime
}}\right. \sim\mathsf{G}+\mathsf{C}_{1}+2\mathsf{C}_{2}, &
(D_{3}^{\prime\prime}\cdot D_{4}^{\prime\prime})\left|
_{D_{4}^{\prime\prime}}\right. =\mathsf{C}_{0},
\end{array}
\end{equation*}
} \noindent where these last $2\cdot7$ intersections concern only the
snc-resolution of the $\mathbf{E}_{8}$-type singularity.

\begin{lemma}
\label{LEMMA1}\emph{(i)} All the edges of the intersection graphs represent
smooth, irreducible, rational compact complex curves.\smallskip \newline
\emph{(ii) }Let $\mathbf{b}\left( X\right) $ denote the total number of the
edges of the intersection graph associated to the desingularization $\varphi
:\widetilde{X}\rightarrow X_{f}=X,$ and let $\mathbf{t}\left( X\right) $ be
the number of those triangles of the graph for which the corresponding three
exceptional prime divisors have non-empty intersection in common. Then each
of the $\mathbf{t}\left( X\right) $ triple non-empty intersections consists
topologically of exactly one point. In addition, $\mathbf{b}\left( X\right) $
and $\mathbf{t}\left( X\right) $ take the following values\emph{:}
\begin{equation*}
\begin{tabular}{|c|c|c|}
\hline
$\mathbf{Types}$ & $\mathbf{b}\left( X\right) $ & $\mathbf{t}\left( X\right)
$ \\ \hline\hline
$\mathbf{A}_{n}$ \ \emph{(}$n$\emph{\ \ odd)} & $
\begin{array}{c}
\  \\
m-1\ (=\dfrac{n-1}{2}) \\
\
\end{array}
$ & $0$ \\ \hline
$\mathbf{A}_{n}$ \ \emph{(}$n$\emph{\ \ even)} & $
\begin{array}{c}
\  \\
m-1\ (=\dfrac{n}{2}) \\
\
\end{array}
$ & $0$ \\ \hline
$\mathbf{D}_{2k}$ & $7(k-1)$ & $3+4\left( k-2\right) $ \\ \hline
$\mathbf{D}_{2k+1}$ & $7k-6$ & $4+4\left( k-2\right) $ \\ \hline
$\mathbf{E}_{6}$ & $9$ & $5$ \\ \hline
$\mathbf{E}_{7}$ & $21$ & $12$ \\ \hline
$\mathbf{E}_{8}$ & $28$ & $17$ \\ \hline
\end{tabular}
\end{equation*}
\newline
\emph{(iii)} In all the cases, there are no four exceptional prime divisors
having non-empty intersection in common.
\end{lemma}

\noindent\textit{Proof. }(i) The underlying topological spaces of all
divisors $\mathsf{H}$, $\mathsf{H}_{1}$, $\mathsf{H}_{2}$, $\mathsf{f}$, $%
\mathsf{f}^{\prime}$, $\mathsf{C}_{0}$, $\mathsf{C}_{1}$, $\mathsf{C}_{2}$, $%
\mathsf{G}$ are in all the cases homeomorphic to $\mathbb{P}_{\mathbb{C}%
}^{1}.$ But also all the other divisors $(D_{i}\cdot D_{j})\left|
_{D_{k}}\right. ,\,k\in\{i,j\},$ for which we gave (just for geometric
reasons and completeness' sake) certain expressions in terms of the
generators of Pic$\left( D_{k}\right) $ up to linear equivalence `$\sim$'$,$
are actually lines (living on $D_{k}$ and being strict transforms of other
lines which are intersections of the exceptional divisors with affine
patches in the previous steps). Therefore they have underlying topological
spaces homeomorphic to $\mathbb{P}_{\mathbb{C}}^{1}.$ (It is better to
compare with the corresponding intersections $(D_{i}\cdot D_{j})\left|
_{D_{\{i,j\}\mathbb{r}\{k\}}}\right. $ for a quick check!$)\smallskip$

\noindent(ii) We find $\mathbf{b}\left( X\right) $ by simply counting all
the edges of each of our graphs. The graph for type $\mathbf{A}_{n}$
contains no triangles. For the remaining types $\mathbf{D}_{2k}$, $\mathbf{D}%
_{2k+1}$, $\mathbf{E}_{6}$, $\mathbf{E}_{7}$, $\mathbf{E}_{8}$, the
intersection graphs contain $3+4\left( k-2\right) $, $\ 5+4\left( k-2\right)
$,\ $7$,$\ 12\ $ and $\ 17$ triangles, respectively, whose vertices are the
only graph-vertices lying on their boundaries. Using the explicitly just
described behaviour of the intersections between the corresponding
exceptional prime divisors, one verifies easily that the number $\mathbf{t}%
\left( X\right) $ equals $3+4\left( k-2\right) ,\ 4+4\left( k-2\right) ,\
5,\ 12$ \ and\ $17,$ respectively. The only triangles which have to be
excluded are those associated to $D_{1}\cap D_{1}^{\prime}\cap
D_{1}^{\prime\prime}=\varnothing$ (for type $\mathbf{D}_{2k+1}$) and to $%
D_{1}\cap D_{4}\cap D_{4}^{\prime }=D_{2}\cap D_{4}\cap
D_{4}^{\prime}=\varnothing$ (for type $\mathbf{E}_{6}$), and each triple
non-empty intersection consists obviously of exactly one point.\smallskip\

\noindent (iii) Examining each (not necessarily convex or non-degenerate)
quadrilateral of the intersection graphs (with no interior points in its
edges), we obtain by the above given data: $D_{i}\cap D_{j}\cap D_{k}\cap
D_{l}=\varnothing $, for all possible pairwise distinct indices $i,j,k,l.$%
\hfill $\square $

\begin{lemma}
\label{LEMMA2}\emph{(i)} The $E$-polynomials of \ $\mathbb{F}_{2}$ and $%
\mathbb{P}_{\mathbb{C}}^{1}\times \mathbb{P}_{\mathbb{C}}^{1}$ are equal%
\emph{:}
\begin{equation}
E\left( \mathbb{F}_{2};u,v\right) =E(\mathbb{P}_{\mathbb{C}}^{1}\times
\mathbb{P}_{\mathbb{C}}^{1};u,v)=1+2\,uv+(uv)^{2}=\left( 1+uv\right) ^{2}
\label{E-F}
\end{equation}
\emph{(ii)} \ $\mathbb{P}_{\mathbb{C}}^{2}\mathbf{[3],\,}\mathbb{P}_{\mathbb{%
C}}^{2}\mathbf{[\,}\overline{\mathbf{3}}\,\mathbf{]}$ and $\mathbb{P}_{%
\mathbb{C}}^{2}\mathbf{[\,}\overline{\overline{\mathbf{3}}}\,\mathbf{]}$
have identical $E$-polynomials, with
\begin{equation}
E(\mathbb{P}_{\mathbb{C}}^{2}\mathbf{[3]};u,v)=E(\mathbb{P}_{\mathbb{C}}^{2}%
\mathbf{[\,}\overline{\mathbf{3}}\,\mathbf{]};u,v)=E(\mathbb{P}_{\mathbb{C}%
}^{2}\mathbf{[\,}\overline{\overline{\mathbf{3}}}\,\mathbf{]}%
;u,v)=1+4\,uv+(uv)^{2}  \label{E-Ps}
\end{equation}
\end{lemma}

\noindent\textit{Proof.} (i) is obvious. (For the fibration $\mathbb{F}%
_{2}\rightarrow$ $\mathbb{P}_{\mathbb{C}}^{1}$ one may use directly (\ref
{E-FIB})). (ii) follows easily from the fact that the $E$-polynomial of a
non-singular surface increases by $uv$ after a blow-up (cf. (\ref{BLOWUP}%
)).\hfill$\square$

\newpage
%
%
%

\section{Computing the discrepancy coefficients\label{DISCREPANCY}}

\noindent\noindent This section is devoted to the exact computation of the
discrepancy coefficients w.r.t. the above snc-desingularizations $\varphi:%
\widetilde{X}\longrightarrow X=X_{f}^{\left( 3\right) }$ $\ $of $3$%
-dimensional $\mathbf{A}$-$\mathbf{D}$-$\mathbf{E}$'s and to a subsequent \
simplification of \ applying formula (\ref{E-STR}).

\begin{proposition}
\label{DISCR}The discrepancies of the snc-desingularizations
\begin{equation*}
\varphi :\widetilde{X}\longrightarrow X
\end{equation*}
of the underlying spaces $X=X_{f}^{\left( 3\right) }$ of the
three-dimensional $\mathbf{A}$-$\mathbf{D}$-$\mathbf{E}$ singularities \emph{%
(}discussed in \emph{\S \ref{CANDES}) }are given by the following table\emph{%
:}
\begin{equation*}
\begin{tabular}{|c|c|}
\hline
$\mathbf{Types}$ & $
\begin{array}{c}
\, \\
\mathbf{Discrepancy}\ \ \ \ K_{\widetilde{X}}-\varphi ^{\ast }\left(
K_{X}\right)  \\
\,
\end{array}
$ \\ \hline\hline
$\mathbf{A}_{n},$ $n$ \emph{even} & $
\begin{array}{c}
\  \\
\sum\limits_{i=1}^{\frac{n}{2}}i\,D_{i}+(n+2)\,D_{\frac{n}{2}+1} \\
\
\end{array}
$ \\ \hline
$\mathbf{A}_{n},$ $n$ \emph{odd} & $
\begin{array}{c}
\  \\
\sum\limits_{i=1}^{\frac{n+1}{2}}i\,D_{i} \\
\
\end{array}
$ \\ \hline
$\mathbf{D}_{n},$ $n$ \emph{even} & $
\begin{array}{c}
\  \\
\left( n-1\right) D_{1}+\left( n-1\right) D_{2}+\sum\limits_{i=3}^{\frac{n}{2%
}+1}\left( 2(n-2i)+7\right) D_{i} \\
\  \\
+\sum\limits_{i=1}^{\frac{n}{2}-1}\left( \frac{n}{2}-i\right) \left(
D_{i}^{\prime }+D_{i}^{\prime \prime }\right) \!\!\! \\
\
\end{array}
$ \\ \hline
$\mathbf{D}_{n},$ $n$ \emph{odd} & $
\begin{array}{c}
\  \\
\left( n-2\right) D_{1}+\left( n-1\right) D_{2}+\sum\limits_{i=3}^{\frac{n+1%
}{2}}\left( 2(n-2i-1)+7\right) D_{i} \\
\  \\
+\sum\limits_{i=1}^{\frac{n-3}{2}}\left( \frac{n-1}{2}-i\right) \left(
D_{i}^{\prime }+D_{i}^{\prime \prime }\right) \!\!\! \\
\
\end{array}
$ \\ \hline
$\mathbf{E}_{6}$ & $
\begin{array}{c}
\  \\
3D_{1}+6D_{2}+9D_{3}+D_{4}+D_{4}^{\prime } \\
\
\end{array}
$ \\ \hline
$\mathbf{E}_{7}$ & $
\begin{array}{c}
\, \\
\allowbreak 11D_{1}+9D_{2}+13D_{3}+5D_{4}+4D_{1}^{\prime }+4D_{1}^{\prime
\prime } \\
\  \\
+2D_{2}^{\prime }+2D_{2}^{\prime \prime }+D_{3}^{\prime }+D_{3}^{\prime
\prime } \\
\,
\end{array}
$ \\ \hline
$\mathbf{E}_{8}$ & $
\begin{array}{c}
\, \\
\allowbreak 19D_{1}+15D_{2}+23D_{3}+11D_{4}+7D_{1}^{\prime }+7D_{1}^{\prime
\prime } \\
\  \\
+4D_{2}^{\prime }+4D_{2}^{\prime \prime }+2D_{3}^{\prime }+2D_{3}^{\prime
\prime }+D_{4}^{\prime }+D_{4}^{\prime \prime } \\
\,
\end{array}
$ \\ \hline
\end{tabular}
\end{equation*}
\end{proposition}

\noindent \textit{Proof. }\ By construction, $\varphi :\widetilde{X}%
\longrightarrow X$ is composed of ``partial'' resolution morphisms. To use a
uniform notation (from a \textit{global} point of view) in what follows, we
shall write $\varphi =\varphi _{1}\circ \varphi _{2}\circ \cdots \circ
\varphi _{\nu }$ and
\begin{equation}
\widetilde{X}=X_{\nu }\overset{\varphi _{\nu }}{\longrightarrow }X_{\nu -1}%
\overset{\varphi _{\nu -1}}{\longrightarrow }\cdots \overset{\varphi _{3}}{%
\longrightarrow }X_{2}\overset{\varphi _{2}}{\longrightarrow }X_{1}\overset{%
\varphi _{1}}{\longrightarrow }X_{0}=X  \label{COMPOSITES}
\end{equation}
for these partial resolutions (where $\nu =\left\lfloor \frac{n+2}{2}%
\right\rfloor ,\left\lfloor \frac{n+1}{2}\right\rfloor ,4,4,5$ for types $%
\mathbf{A}_{n}$, $\mathbf{D}_{n}$, $\mathbf{E}_{6}$, $\mathbf{E}_{7}$, and $%
\mathbf{E}_{8}$, respectively, as one deduces from \S \ref{CANDES}). The
discrepancy w.r.t. $\varphi $ equals:
\begin{equation}
\begin{array}{l}
K_{\widetilde{X}}-\varphi ^{\ast }\left( K_{X}\right) = \\
\  \\
\sum\limits_{i=1}^{\nu -1}\left( \varphi _{i+1}\circ \varphi _{i+2}\circ
\cdots \circ \varphi _{\nu }\right) ^{\ast }\left( K_{X_{i}}-\varphi
_{i}^{\ast }\left( K_{X_{i-1}}\right) \right) +K_{X_{\nu }}-\varphi _{\nu
}^{\ast }\left( K_{X_{\nu -1}}\right)
\end{array}
\label{PARDIS}
\end{equation}
Therefore, for its computation, it suffices to determine the discrepancies
w.r.t. each of the $\varphi _{i}$'s, and then to specify the pull-backs
which are involved in (\ref{PARDIS}).\bigskip

\noindent {}\textbf{I) Computation of the intermediate discrepancies}. Since
the arising singularities are isolated, we may investigate the zeros of
canonical differentials locally around them.\medskip

\noindent \textbf{(i) Type }$\mathbf{A}_{n}$\textbf{. }The defining
polynomial of the singularity is
\begin{equation}
f(x_{1},\ldots ,x_{4})=x_{1}^{n+1}+x_{2}^{2}+x_{3}^{2}+x_{4}^{2}\ .
\label{EQAN}
\end{equation}
Let $n\geq 2$, and consider the rational canonical differential
\begin{equation*}
\mathfrak{s:}=\text{Res}_{X}\left( \frac{dx_{1}\wedge dx_{2}\wedge
dx_{3}\wedge dx_{4}}{f}\right) =\frac{dx_{2}\wedge dx_{3}\wedge dx_{4}}{%
\left( \partial f\,/\,\partial x_{1}\right) }\in \Omega _{\mathbb{C}\left(
X\right) /\mathbb{C}}^{3}\ .
\end{equation*}
$\mathfrak{s}$ is a basis of the dualizing sheaf $\omega _{X}=\mathcal{O}%
_{X}(K_{X})=(\Omega _{X}^{3})^{\vee \,\vee }$ \ whose sections are defined
by
\begin{equation*}
\left\{
\begin{array}{c}
\text{open sets } \\
\text{of }X
\end{array}
\right\} \ni V\longmapsto \Gamma \left( V,\omega _{X}\right) :=\left\{ %
\mathfrak{y}\in \Omega _{\mathbb{C}\left( X\right) /\mathbb{C}}^{3}\left|
\begin{array}{l}
\mathfrak{y}\text{ \ is a } \\
\text{regular canonical} \\
\text{differential } \\
\text{on \ }V\cap (X\mathbb{r}\{\mathbf{0}\})
\end{array}
\right. \right\} .
\end{equation*}
Blow up $X$ at $\mathbf{0}$ and consider the affine piece $U_{1}\cap $ $%
\mathbf{Bl}_{\mathbf{0}}(X),$ with
\begin{equation*}
U_{1}=\text{Spec}\left( \mathbb{C}\left[ y_{1,1},y_{1,2},y_{1,3},y_{1,4}%
\right] \right) .
\end{equation*}
The restriction of the exceptional locus $\mathcal{E}_{f}$ on $U_{1}$ is
nothing but
\begin{equation*}
\mathbf{Bl}_{\mathbf{0}}(X)\cap \mathcal{E}_{1}=\mathcal{E}_{f}\left|
_{U_{1}}\right. =\left\{ \left( y_{1,1},..,y_{1,4}\right) \in \mathbb{C}%
^{4}\left| y_{1,1}=\widetilde{f_{1}}\left( y_{1,1},..,y_{1,4}\right)
=0\right. \right\}
\end{equation*}
where
\begin{equation*}
\widetilde{f_{1}}\left( y_{1,1},y_{1,2},y_{1,3},y_{1,4}\right)
=y_{1,1}^{n-1}+y_{1,2}^{2}+y_{1,3}^{2}+y_{1,4}^{2}\ .
\end{equation*}
(As we explained before, the possibly existing new ($\mathbf{A}_{n-2}$%
\thinspace -) singularity on $\mathbf{Bl}_{\mathbf{0}}(X)$ lies in $\mathcal{%
E}_{f}\left| _{U_{1}}\right. $). To find the discrepancy coefficient w.r.t. $%
\mathbf{Bl}_{\mathbf{0}}(X)\longrightarrow X,$ it suffices to compare $%
\mathfrak{s}$ with the rational canonical differential
\begin{equation*}
\overline{\mathfrak{s}}:=\frac{dy_{1,2}\wedge dy_{1,3}\wedge dy_{1,4}}{%
(\partial \widetilde{f_{1}}\,/\,\partial y_{1,1})}\in \Omega _{\mathbb{C}%
\left( U_{1}\right) /\mathbb{C}}^{3}\ .
\end{equation*}
($U_{1}$ is non-singular with local coordinates $y_{1,2},y_{1,3},y_{1,4}$ at
any point $q$ for which $\partial \widetilde{f_{1}}(q)\,/\,\partial
y_{1,1}\neq 0$). In $U_{1}$ we have $x_{1}=y_{1,1}$ and $x_{j}=x_{1}\,\xi
_{j}=y_{1,1}\,y_{1,j},$ for all $j\in \{2,3,4\}.$ Hence,
\begin{align}
&
\begin{array}{l}
dx_{2}\wedge dx_{3}\wedge dx_{4} \\
=\left( y_{1,1}dy_{1,2}+y_{1,2}\,dy_{1,1}\right) \wedge \left(
y_{1,1}dy_{1,3}+y_{1,3}\,dy_{1,1}\right) \wedge \left(
y_{1,1}dy_{1,4}+y_{1,4}\,dy_{1,1}\right)
\end{array}
\notag \\
& =y_{1,1}^{2}\ (y_{1,2}\,dy_{1,1}\wedge dy_{1,3}\wedge
dy_{1,4}-y_{1,3}\,dy_{1,1}\wedge dy_{1,2}\wedge dy_{1,4}  \notag \\
& +y_{1,4}\,dy_{1,1}\wedge dy_{1,2}\wedge dy_{1,3}+y_{1,1}\,dy_{1,2}\wedge
dy_{1,3}\wedge dy_{1,4})  \label{DIF1}
\end{align}
and
\begin{equation}
\partial f\,/\,\partial x_{1}=\left( n+1\right) \,x_{1}^{n}=\left(
n+1\right) \,\,y_{1,1}^{n}=\left( \frac{n+1}{n-1}\right) \ y_{1,1}^{2}\
(\partial \widetilde{f_{1}}\,/\,\partial y_{1,1})  \label{DIF2}
\end{equation}
On the other hand, note that
\begin{equation*}
d\widetilde{f_{1}}=\left( n-1\right) \ y_{1,1}^{n-2}\,dy_{1,1}+2\,\left(
y_{1,2}\,dy_{1,2}+y_{1,3}\,dy_{1,3}+y_{1,4}\,dy_{1,4}\right) =0
\end{equation*}
if and only if
\begin{equation}
dy_{1,1}=-\frac{2}{n-1}\ y_{1,1}^{2-n}\ \left(
y_{1,2}\,dy_{1,2}+y_{1,3}\,dy_{1,3}+y_{1,4}\,dy_{1,4}\right)  \label{DIF3}
\end{equation}
Substituting the expression (\ref{DIF3}) for $dy_{1,1}$ into the right-hand
side of (\ref{DIF1}), we obtain easily
\begin{equation}
\begin{array}{l}
dx_{2}\wedge dx_{3}\wedge dx_{4} \\
=\left( -\frac{2}{n-1}\ y_{1,1}^{4-n}\
(y_{1,2}^{2}+y_{1,3}^{2}+y_{1,4}^{2})+y_{1,1}^{3}\right) \ dy_{1,2}\wedge
dy_{1,3}\wedge dy_{1,4}
\end{array}
\   \label{DIF4}
\end{equation}
Combining now (\ref{DIF4}) with $%
y_{1,2}^{2}+y_{1,3}^{2}+y_{1,4}^{2}=-y_{1,1}^{n-1}$ and (\ref{DIF2}), we get
\begin{equation}
\mathfrak{s}=\frac{\left( \frac{n+1}{n-1}\ y_{1,1}^{3}\right) \
dy_{1,2}\wedge dy_{1,3}\wedge dy_{1,4}}{\left( \frac{n+1}{n-1}\right) \
y_{1,1}^{2}\ (\partial \widetilde{f_{1}}\,/\,\partial y_{1,1})}=y_{1,1}\
\overline{\mathfrak{s}}  \label{DIF5}
\end{equation}
The equality (\ref{DIF5}) shows that the discrepancy coefficient of the
exceptional prime divisor $\mathcal{E}_{f}$ w.r.t. $\mathbf{Bl}_{\mathbf{0}%
}(X)\longrightarrow X$ equals $1$.\medskip

\noindent {}If $n=1,$ then we compare
\begin{equation*}
\mathfrak{s}=\frac{dx_{1}\wedge dx_{2}\wedge dx_{3}}{\left( \partial
f\,/\,\partial x_{4}\right) }\text{ \ \ \ \ with \ \ \ \ }\overline{%
\mathfrak{s}}=\frac{dy_{1,1}\wedge dy_{1,2}\wedge dy_{1,3}}{(\partial
\widetilde{f_{1}}\,/\,\partial y_{1,4})}.
\end{equation*}
Since $\frac{\partial f}{\partial x_{4}}=2x_{4}=2y_{1,1}y_{1,4}$, $\frac{%
\partial \widetilde{f_{1}}}{\partial y_{1,4}}=2y_{1,4}$, and
\begin{equation*}
dx_{1}\wedge dx_{2}\wedge dx_{3}=y_{1,1}^{2}dy_{1,1}\wedge dy_{1,2}\wedge
dy_{1,3},
\end{equation*}
we conclude again $\mathfrak{s}=y_{1,1}\overline{\mathfrak{s}}.$ In fact,
this kind of argumentation covers all but one steps of the resolution
procedure for $\mathbf{A}_{n}$'s. The indicated ``special'' case occurs only
in the last step and only for $n$ even, where we blow-up once more to get
rid of the singularity of the exceptional locus for the purpose of ensuring
the snc-condition for $\varphi :\widetilde{X}\longrightarrow X$ ($``n=0$%
''-case). But since we blow-up a point which is \textit{smooth} on the
\textit{3-fold}, the discrepancy coefficient of the lastly created
exceptional prime divisor $D_{\frac{n}{2}+1}$ equals $2$ (see remark \ref
{REMARK-S} and Griffiths \& Harris \cite[Lemma of p. 187]{GR-H})$.\medskip $

\noindent \textbf{(ii) Type }$\mathbf{D}_{n}$\textbf{. }For this type we
proceed analogously by making use of the affine piece $U_{1}$. The only
difference here is that the exceptional divisor $\mathcal{E}_{f}$ under the
first blow-up has two irreducible components $\mathcal{E}_{f}^{\prime }$ and
$\mathcal{E}_{f}^{\prime \prime }$. Nevertheless, the corresponding local
computation with rational canonical differentials gives again
\begin{equation*}
\frac{dx_{2}\wedge dx_{3}\wedge dx_{4}}{\left( \partial f\,/\,\partial
x_{1}\right) }=y_{1,1}\ \frac{dy_{1,2}\wedge dy_{1,3}\wedge dy_{1,4}}{%
(\partial \widetilde{f_{1}}\,/\,\partial y_{1,1})}
\end{equation*}
and the discrepancy coefficient for both of them equals $1$. As it is clear
from Lemma \ref{LRL} and \textbf{(i)}, the discrepancy coefficients in all
resolution steps will be again $1$.\medskip\

\noindent \textbf{(iii) Types }$\mathbf{E}_{6},\mathbf{E}_{7},\mathbf{E}_{8}$%
\textbf{. \ }For these types one may work along the same lines with respect
to the affine piece $U_{2}=$ Spec$\left( \mathbb{C}\left[
y_{2,1},y_{2,2},y_{2,3},y_{2,4}\right] \right) .$ The exceptional divisor $%
\mathcal{E}_{f}$ w.r.t. $\mathbf{Bl}_{\mathbf{0}}(X)\longrightarrow X$
consists again of two prime ones. Each of them has discrepancy coefficient
equal to $1.$ This property remains also valid for all other composites (\ref
{COMPOSITES}) of $\varphi ,$ exactly as in the case of type $\mathbf{D}_{n}.$
Further details will be omitted$.\medskip $

\noindent {}Recapitulating, we should stress that in \textbf{(i)}, \textbf{%
(ii)}, \textbf{(iii)}, the discrepancy coefficient for \textit{each} of the
prime divisors of the exceptional locus of the $\varphi _{i}$'s in (\ref
{COMPOSITES}) equals $1$, up to the last resolution morphism for type $%
\mathbf{A}_{n}$, $n$ even, which has discrepancy $2$. This fact will be used
below in an essential way.\bigskip

\noindent {}\textbf{II) Computation of the pull-backs}. To determine the
required pullbacks of our discrepancies (see (\ref{COMPOSITES}), (\ref
{PARDIS})), we shall denote by $E_{j}$ (resp., $E_{j}^{^{\prime (\prime
\prime )}}$) those exceptional prime divisors which are created (for the
first time) after the application of a $\varphi _{i}$ (i.e., actually the
members of $\mathfrak{Ex}\left( \varphi _{i}\right) $), so that their strict
transforms (on $\widetilde{X}$) are exactly the exceptional prime divisors
(w.r.t. $\varphi $) which are denoted by $D_{j}$ (resp., $D_{j}^{^{\prime
(\prime \prime )}}$) in \S \ref{CANDES}.\medskip

\noindent {}\textbf{(i) Type }$\mathbf{A}_{n}.$ Defining $m=\left\lfloor
\frac{n+2}{2}\right\rfloor ,$ as in \S \ref{CANDES}, $\varphi $ is
decomposed into $m$ birational morphisms:
\begin{equation*}
\widetilde{X}=X_{m}\overset{\varphi _{m}}{\longrightarrow }X_{m-1}\overset{%
\varphi _{m-1}}{\longrightarrow }\cdots \overset{\varphi _{3}}{%
\longrightarrow }X_{2}\overset{\varphi _{2}}{\longrightarrow }X_{1}\overset{%
\varphi _{1}}{\longrightarrow }X_{0}=X.
\end{equation*}
Each $\varphi _{i}(=\pi _{i}$ of \S \ref{CANDES}$)$ gives rise to an
exceptional prime divisor $E_{i}.$ By \textbf{I)} we get
\begin{equation}
K_{X_{i}}-\varphi _{i}^{\ast }(K_{X_{i-1}})=E_{i},\ \ \forall i,\ 1\leq
i\leq m-1,  \label{FORMULAA}
\end{equation}
and
\begin{equation}
K_{X_{m}}-\varphi _{m}^{\ast }(K_{X_{m-1}})=\left\{
\begin{array}{ll}
D_{m}, & \text{if \ }n\text{ \ is odd,} \\
\  & \  \\
2D_{m}, & \text{if \ }n\text{ \ is even.}
\end{array}
\right.   \label{FORMULAB}
\end{equation}
We claim that for all $i$, $1\leq i\leq m-1,$
\begin{equation}
\left( \varphi _{i+1}\circ \varphi _{i+2}\circ \cdots \circ \varphi
_{m}\right) ^{\ast }\left( E_{i}\right) =\left\{
\begin{array}{ll}
\sum\limits_{j=i}^{m}D_{j}, & \text{if \ }n\text{ \ is odd,} \\
\  & \  \\
\sum\limits_{j=i}^{m}D_{j}+2D_{m}, & \text{if \ }n\text{ \ is even.}
\end{array}
\right.   \label{FORMULAC}
\end{equation}
To prove (\ref{FORMULAC}) we shall work with local equations for the
corresponding divisors. Consider two successive blow-ups
\begin{equation*}
X_{j+1}\overset{\varphi _{j+1}}{\longrightarrow }X_{j}\overset{\varphi _{j}}{%
\longrightarrow }X_{j-1}
\end{equation*}
and assume that $X_{j}$ has a singularity of type $\mathbf{A}_{n},$ $n\geq 1$%
, (with equation (\ref{EQAN})), where $\varphi _{j}$ denotes the blow-up of
the $\mathbf{A}_{n+2}$-singularity of $X_{j-1}.$ The local equation $(%
\widetilde{f_{2}}=0)$ is the equation of $X_{j+1}$ on the affine chart $%
U_{2}=$ Spec$\left( \mathbb{C}\left[ y_{2,1},..,y_{2,4}\right] \right) $,
where
\begin{equation*}
\widetilde{f_{2}}\left( y_{2,1},y_{2,2},y_{2,3},y_{2,4}\right)
=y_{2,1}^{n+1}y_{2,2}^{n-1}+1+y_{2,3}^{2}+y_{2,4}^{2}
\end{equation*}
(cf. \S \ref{CANDES}). The new exceptional locus $E_{j+1}$ of $\varphi _{j+1}
$ on $U_{2}\cap X_{j+1}$ is given by the local equation $(y_{2,1}=0).$ On
the other hand, $(x_{1}=0)$ and $(y_{2,2}=0)$ express the local equations
for $E_{j}$ on $X_{j}$ and for its strict transform $E_{j\text{, \textsc{st}}%
}$ on $U_{2}\cap X_{j+1},$ respectively. Since the preimage of $(x_{1}=0)$
under $\varphi _{j+1}$ equals $(y_{2,1}\cdot y_{2,2}=0),$ we have:
\begin{equation}
\varphi _{j+1}^{\ast }(E_{j})=E_{j+1}+E_{j\text{, \textsc{st}}}.
\label{EQ-STRT}
\end{equation}
It remains to see what happens in the case in which $\varphi _{j+1}$ is the
blow up of a (regular) $\mathbf{A}_{0}$-point, i.e., whenever $j=m-1=k$ and $%
X_{k+1}$ is the last step of the resolution process for a singularity of
type $\mathbf{A}_{2k}$. For $n=0$, we get equations
\begin{equation*}
x_{1}+x_{2}^{2}+x_{3}^{2}+x_{4}^{2}=0\ \ \text{and}\ \ \
z_{2,1}+z_{2,2}(1+z_{2,3}^{2}+z_{2,4}^{2})=0,
\end{equation*}
on $X_{k}$ and $U_{2}\cap X_{k+1}$, respectively. The divisors $D_{k+1}$, $%
E_{k}$, $E_{k\text{, \textsc{st}}}$ have local equations $(z_{2,2}=0)$, $%
(x_{1}=0)$ and $(z_{2,1}=0)$, respectively. Since
\begin{equation*}
x_{1}=z_{2,1}\,z_{2,2}=z_{2,2}^{2}(1+z_{2,3}^{2}+z_{2,4}^{2}),
\end{equation*}
we deduce
\begin{equation}
\varphi _{k+1}^{\ast }(E_{k})=2D_{k+1}+E_{k\text{, \textsc{st}}%
}=2D_{m}+E_{m-1\text{, \textsc{st}}}.  \label{EQ-STRT2}
\end{equation}
(\ref{FORMULAC}) follows after repeated application of equations like (\ref
{EQ-STRT}) and (\ref{EQ-STRT2}).\ Now inserting the data of (\ref{FORMULAA}%
), (\ref{FORMULAB}), (\ref{FORMULAC}) into (\ref{PARDIS}) we obtain:
\begin{equation*}
K_{\widetilde{X}}-\varphi ^{\ast }\left( K_{X}\right) =\left\{
\begin{array}{ll}
\sum\limits_{i=1}^{\frac{n+1}{2}}i\,D_{i}, & \text{if \ }n\text{ \ is odd,}
\\
\  & \  \\
\sum\limits_{i=1}^{\frac{n}{2}}i\,D_{i}+(n+2)\,D_{\frac{n}{2}+1}, & \text{if
\ }n\text{ \ is even.}
\end{array}
\right.
\end{equation*}

\noindent {}\textbf{(ii) Type }$\mathbf{D}_{n}$, $n=2k.$ In this case $%
\varphi $ is decomposed into $k$ birational morphisms:
\begin{equation*}
\widetilde{X}=X_{k}\overset{\varphi _{k}}{\longrightarrow }X_{k-1}\overset{%
\varphi _{k-1}}{\longrightarrow }\cdots \overset{\varphi _{3}}{%
\longrightarrow }X_{2}\overset{\varphi _{2}}{\longrightarrow }X_{1}\overset{%
\varphi _{1}}{\longrightarrow }X_{0}=X.
\end{equation*}
By construction, $\mathfrak{Ex}\left( \varphi _{1}\right) =\left\{
E_{k-1}^{\prime },E_{k-1}^{\prime \prime }\right\} ,$
\begin{equation*}
\mathfrak{Ex}\left( \varphi _{i+1}\right) =\left\{ E_{k-i-1}^{\prime
},E_{k-i-1}^{\prime \prime },E_{k-i+2}\right\} ,\ \forall i,\ 1\leq i\leq
k-2,
\end{equation*}
and $\mathfrak{Ex}\left( \varphi _{k}\right) =\left\{
D_{1},D_{2},D_{3}\right\} $. By \textbf{I)} we have
\begin{equation*}
\begin{array}{l}
K_{X_{1}}-\varphi _{1}^{\ast }(K_{X_{0}})=E_{k-1}^{\prime }+E_{k-1}^{\prime
\prime }, \\
\  \\
K_{X_{i+1}}-\varphi _{i+1}^{\ast }(K_{X_{i}})=E_{k-i-1}^{\prime
}+E_{k-i-1}^{\prime \prime }+E_{k-i+2},\ \forall i,\ 1\leq i\leq k-2, \\
\  \\
K_{X_{k}}-\varphi _{k}^{\ast }(K_{X_{k-1}})=D_{1}+D_{2}+D_{3}.
\end{array}
\end{equation*}
We shall prove that
\begin{equation}
K_{\widetilde{X}}-\varphi ^{\ast }\left( K_{X}\right) =\left( 2k-1\right)
(D_{1}+D_{2})+\sum\limits_{i=1}^{k-1}i(D_{k-i}^{\prime }+D_{k-i}^{\prime
\prime })+\sum\limits_{j=1}^{k-1}(4j-1)D_{k-j+2}.  \label{KXDACH}
\end{equation}
For $k=2$ this can be shown easily. Suppose that $\ k\geq 3$. Then
\begin{equation*}
\begin{array}{l}
\varphi _{i+1}^{\ast }(E_{k-i}^{^{\prime (\prime \prime
)}})=E_{k-i-1}^{^{\prime (\prime \prime )}}+E_{k-i+2}+E_{k-i,\text{\textsc{st%
}}}^{^{\prime (\prime \prime )}},\ \forall i,\ 1\leq i\leq k-2, \\
\  \\
\varphi _{k}^{\ast }(E_{1}^{^{\prime (\prime \prime
)}})=D_{1}+D_{2}+D_{3}+D_{1}^{^{\prime (\prime \prime )}},
\end{array}
\end{equation*}
and for all $i$, $1\leq i\leq k-2,$
\begin{equation*}
\begin{array}{l}
\left( \varphi _{i+1}\circ \varphi _{i+2}\right) ^{\ast }(E_{k-i}^{^{\prime
(\prime \prime )}})=E_{k-i+1}+E_{k-i+2}+E_{k-i-1,\text{\textsc{st}}%
}^{^{\prime (\prime \prime )}}+E_{k-i,\text{\textsc{st}, \textsc{st}}%
}^{^{\prime (\prime \prime )}}.
\end{array}
\end{equation*}
This means that
\begin{equation*}
\begin{array}{l}
\left( \varphi _{2}\circ \varphi _{3}\circ \cdots \circ \varphi _{k}\right)
^{\ast }\left( E_{k-1}^{\prime }+E_{k-1}^{\prime \prime }\right) \\
\  \\
=\sum\limits_{j=1}^{k-1}\left( D_{k-j}^{\prime }+D_{k-j}^{\prime \prime
}\right) +2(D_{1}+D_{2}+D_{3})+2\left(
D_{3}+2\sum\limits_{j=2}^{k-2}D_{k-j+2}+D_{k+1}\right) ,
\end{array}
\end{equation*}
and that for all $i$, $2\leq i\leq k-2,$%
\begin{equation*}
\begin{array}{l}
\left( \varphi _{i+1}\circ \varphi _{i+2}\circ \cdots \circ \varphi
_{k}\right) ^{\ast }(E_{k-i}^{\prime }+E_{k-i}^{\prime \prime }+E_{k-i+3})
\\
\  \\
=\sum\limits_{j=i}^{k-1}\left( D_{k-j}^{\prime }+D_{k-j}^{\prime \prime
}\right) +2(D_{1}+D_{2}+D_{3}) \\
\  \\
+2\left( D_{3}+2\sum\limits_{j=i+1}^{k-2}D_{k-j+2}+D_{k-i+2}\right)
+D_{k-i+3}
\end{array}
\end{equation*}
and
\begin{equation*}
\varphi _{k}^{\ast }(E_{1}^{\prime }+E_{1}^{\prime \prime }+E_{4})=\left(
D_{1}^{\prime }+D_{1}^{\prime \prime }\right) +2(D_{1}+D_{2}+D_{3})+D_{4}.
\end{equation*}
Thus, (\ref{PARDIS}) implies (\ref{KXDACH}).\medskip

\noindent {}\textbf{(iii) Type }$\mathbf{D}_{n}$, $n=2k+1.$ Here $\varphi $
is decomposed into $k+1$ birational morphisms:
\begin{equation*}
\widetilde{X}=X_{k+1}\overset{\varphi _{k+1}}{\longrightarrow }X_{k}\overset{%
\varphi _{k}}{\longrightarrow }\cdots \overset{\varphi _{3}}{\longrightarrow
}X_{2}\overset{\varphi _{2}}{\longrightarrow }X_{1}\overset{\varphi _{1}}{%
\longrightarrow }X_{0}=X.
\end{equation*}
Computing the total discrepancy, we find analogously:
\begin{equation*}
K_{\widetilde{X}}-\varphi ^{\ast }\left( K_{X}\right) =\left( 2k-1\right)
D_{1}+2kD_{2}+\sum\limits_{i=1}^{k-1}i(D_{k-i}^{\prime }+D_{k-i}^{\prime
\prime })+\sum\limits_{j=1}^{k-1}(4j-1)D_{k-j+2}.
\end{equation*}

\noindent {}\textbf{(iv) Type }$\mathbf{E}_{6}$\textbf{. }In this case $%
\varphi $ is decomposed into $4$ birational morphisms:
\begin{equation*}
\widetilde{X}=X_{4}\overset{\varphi _{4}}{\longrightarrow }X_{3}\overset{%
\varphi _{3}}{\longrightarrow }X_{2}\overset{\varphi _{2}}{\longrightarrow }%
X_{1}=\mathbf{Bl}_{\mathbf{0}}(X)\overset{\varphi _{1}}{\longrightarrow }%
X_{0}=X
\end{equation*}
By construction,
\begin{equation*}
\mathfrak{Ex}\left( \varphi _{1}\right) =\left\{ E_{4},E_{4}^{\prime
}\right\} ,\ \ \mathfrak{Ex}\left( \varphi _{2}\right) =\left\{
E_{1}\right\} ,\ \ \mathfrak{Ex}\left( \varphi _{3}\right) =\left\{
E_{2}\right\} ,
\end{equation*}
and $\ \mathfrak{Ex}\left( \varphi _{4}\right) =\left\{ D_{3}\right\} $
(where $\varphi _{i}=\pi _{i-1}$ of \S \ref{CANDES}). By \textbf{I)} we have
\begin{equation*}
\begin{array}{ll}
K_{X_{1}}-\varphi _{1}^{\ast }\left( K_{X_{0}}\right) =E_{4}+E_{4}^{\prime },
& K_{X_{2}}-\varphi _{2}^{\ast }\left( K_{X_{1}}\right) =E_{1}, \\
\  & \  \\
K_{X_{3}}-\varphi _{3}^{\ast }\left( K_{X_{2}}\right) =E_{2}, &
K_{X_{4}}-\varphi _{4}^{\ast }\left( K_{X_{3}}\right) =D_{3}.
\end{array}
\end{equation*}
The intersection diagrams imply
\begin{equation*}
\begin{array}{l}
\left( \varphi _{2}\circ \varphi _{3}\circ \varphi _{4}\right) ^{\ast
}\left( E_{4}+E_{4}^{\prime }\right)
=2D_{1}+4D_{2}+6D_{3}+D_{4}+D_{4}^{\prime }, \\
\  \\
\left( \varphi _{3}\circ \varphi _{4}\right) ^{\ast }\left( E_{1}\right)
=D_{1}+D_{2}+D_{3} , \\
\  \\
\varphi _{4}^{\ast }\left( E_{2}\right) =D_{2}+D_{3}.
\end{array}
\end{equation*}
Hence, by (\ref{PARDIS}), the discrepancy w.r.t. $\varphi $ equals $%
3D_{1}+6D_{2}+9D_{3}+D_{4}+D_{4}^{\prime }.\medskip $

\noindent {}\textbf{(v) Type }$\mathbf{E}_{7}$\textbf{. }Here $\varphi $ is
decomposed into $4$ birational morphisms:\textbf{\ }
\begin{equation*}
\widetilde{X}=X_{4}\overset{\varphi _{4}}{\longrightarrow }X_{3}\overset{%
\varphi _{3}}{\longrightarrow }X_{2}\overset{\varphi _{2}}{\longrightarrow }%
X_{1}=\mathbf{Bl}_{\mathbf{0}}(X)\overset{\varphi _{1}}{\longrightarrow }%
X_{0}=X
\end{equation*}
By construction,
\begin{equation*}
\mathfrak{Ex}\left( \varphi _{1}\right) =\left\{ E_{3}^{\prime
},E_{3}^{\prime \prime }\right\} ,\ \ \mathfrak{Ex}\left( \varphi
_{2}\right) =\left\{ E_{2}^{\prime },E_{2}^{\prime \prime }\right\} ,\ \ %
\mathfrak{Ex}\left( \varphi _{3}\right) =\left\{ E_{1}^{\prime
},E_{1}^{\prime \prime },E_{4}\right\} ,
\end{equation*}
and $\ \mathfrak{Ex}\left( \varphi _{4}\right) =\left\{
D_{1},D_{2},D_{3}\right\} $. By \textbf{I)} we obtain
\begin{equation*}
\begin{array}{ll}
K_{X_{1}}-\varphi _{1}^{\ast }\left( K_{X_{0}}\right) =E_{3}^{\prime
}+E_{3}^{\prime \prime }, & K_{X_{2}}-\varphi _{2}^{\ast }\left(
K_{X_{1}}\right) =E_{2}^{\prime }+E_{2}^{\prime \prime }, \\
\  & \  \\
K_{X_{3}}-\varphi _{3}^{\ast }\left( K_{X_{2}}\right) =E_{1}^{\prime
}+E_{1}^{\prime \prime }+E_{4}, & K_{X_{4}}-\varphi _{4}^{\ast }\left(
K_{X_{3}}\right) =D_{1}+D_{2}+D_{3}.
\end{array}
\end{equation*}
The computation of the pull-backs gives
\begin{equation*}
\begin{array}{l}
\left( \varphi _{2}\circ \varphi _{3}\circ \varphi _{4}\right) ^{\ast}
\left( E_{3}^{\prime }+E_{3}^{\prime \prime }\right) = \\
6D_1+4D_2+6D_3+2D_4+2(D_{1}^{\prime }+D_{1}^{\prime \prime })
+D_{2}^{\prime }+D_{2}^{\prime \prime }
+D_{3}^{\prime }+D_{3}^{\prime \prime }, \\
\  \\
\left( \varphi _{3}\circ \varphi _{4}\right) ^{\ast }
\left( E_{2}^{\prime}+E_{2}^{\prime \prime }\right) =
2D_1+2D_2+4D_3+2D_4
+D_{1}^{\prime }+D_{1}^{\prime \prime }
+D_{2}^{\prime }+D_{2}^{\prime \prime }, \\
\  \\
\varphi _{4}^{\ast }\left( E_{1}^{\prime }+E_{1}^{\prime \prime}+E_{4}\right) =
2D_1+2D_2+2D_3+D_4
+D_{1}^{\prime }+D_{1}^{\prime \prime }.
\end{array}
\end{equation*}
Now apply (\ref{PARDIS}).\medskip

\noindent {}\textbf{(vi) Type }$\mathbf{E}_{8}$\textbf{. }In this case $%
\varphi $ is decomposed into $5$ birational morphisms:\textbf{\ }
\begin{equation*}
\widetilde{X}=X_{5}\overset{\varphi _{5}}{\longrightarrow }X_{4}\overset{%
\varphi _{4}}{\longrightarrow }X_{3}\overset{\varphi _{3}}{\longrightarrow }%
X_{2}\overset{\varphi _{2}}{\longrightarrow }X_{1}\overset{\varphi _{1}}{%
\longrightarrow }X_{0}=X
\end{equation*}
By construction, $\mathfrak{Ex}\left( \varphi _{1}\right) =\left\{
E_{4}^{\prime },E_{4}^{\prime \prime }\right\} ,$
\begin{equation*}
\mathfrak{Ex}\left( \varphi _{2}\right) =\left\{ E_{3}^{\prime
},E_{3}^{\prime \prime }\right\} ,\ \ \mathfrak{Ex}\left( \varphi
_{3}\right) =\left\{ E_{2}^{\prime },E_{2}^{\prime \prime }\right\} ,\ \ %
\mathfrak{Ex}\left( \varphi _{4}\right) =\left\{ E_{1}^{\prime
},E_{1}^{\prime \prime },E_{4}\right\} ,
\end{equation*}
and $\ \mathfrak{Ex}\left( \varphi _{5}\right) =\left\{
D_{1},D_{2},D_{3}\right\} $. By \textbf{I)} we have
\begin{equation*}
\begin{array}{ll}
K_{X_{1}}-\varphi _{1}^{\ast }\left( K_{X_{0}}\right) =E_{4}^{\prime
}+E_{4}^{\prime \prime }, & K_{X_{2}}-\varphi _{2}^{\ast }\left(
K_{X_{1}}\right) =E_{3}^{\prime }+E_{3}^{\prime \prime }, \\
\  & \  \\
K_{X_{3}}-\varphi _{3}^{\ast }\left( K_{X_{2}}\right) =E_{2}^{\prime
}+E_{2}^{\prime \prime }, & K_{X_{4}}-\varphi _{4}^{\ast }\left(
K_{X_{3}}\right) =E_{1}^{\prime }+E_{1}^{\prime \prime }+E_{4},
\end{array}
\end{equation*}
and $K_{X_{5}}-\varphi _{5}^{\ast }\left( K_{X_{4}}\right)
=D_{1}+D_{2}+D_{3}.$ We obtain
\begin{equation*}
\begin{array}{l}
\left( \varphi _{2}\circ \varphi _{3}\circ \varphi _{4}\circ \varphi
_{5}\right) ^{\ast }\left( E_{4}^{\prime }+E_{4}^{\prime \prime }\right)
= \\
 8D_1+6D_2+10D_3+6D_4
+3(D_{1}^{\prime }+D_{1}^{\prime \prime })
+2(D_{2}^{\prime }+D_{2}^{\prime \prime })
+D_{3}^{\prime }+D_{3}^{\prime \prime }
+D_{4}^{\prime }+D_{4}^{\prime \prime } .
\end{array}
\end{equation*}
The remaining inverse images
$\left( \varphi _{3}\circ \varphi _{4}\circ \varphi _{5}\right) ^{\ast
}\left( E_{3}^{\prime }+E_{3}^{\prime \prime }\right)$,
$\left( \varphi _{4}\circ \varphi _{5}\right) ^{\ast }\left( E_{2}^{\prime
}+E_{2}^{\prime \prime }\right) $
and
$\varphi _{5}^{\ast }\left( E_{1}^{\prime }+E_{1}^{\prime \prime
}+E_{4}\right) $
coincide with (v), where in each case
$\varphi _{i}\circ \dots \circ \varphi _{4}$
has to be replaced by
$\varphi _{i+1}\circ \dots \circ \varphi _{5}$.
Finally, apply again (\ref{PARDIS}).\hfill $\square $

\begin{proposition}
Suppose that \textit{\ }$X=X_{f}^{\left( 3\right) }$ is the underlying space
of an $\mathbf{A}$-$\mathbf{D}$-$\mathbf{E}$-singularity, $\varphi :%
\widetilde{X}\longrightarrow X$ its snc-desingularization, $\mathfrak{Ex}%
\left( \varphi \right) =\left\{ D_{1},..,D_{r}\right\} $ the corresponding
exceptional set with discrepacy coefficients $a_{1},\ldots ,a_{r}$, $%
I:=\{1,2,...,r\},$ and
\begin{equation*}
\frak{R}_{\varphi }:=\{\left. (i,j)\in I^{2}\right| D_{\{i,j\}}\neq
\varnothing \},\ \ \frak{Q}_{\varphi }:=\{\left. (i,j,k)\in I^{3}\right|
D_{\{i,j,k\}}\neq \varnothing \}.
\end{equation*}
Then the string-theoretic $E$-function of $X$ satisfies the following
equality\emph{:}
\begin{equation}
\begin{array}{l}
E_{\text{\emph{str}}}\left( X;u,v\right) =E\left( D_{\varnothing }^{\circ
};u,v\right) +\sum\limits_{i=1}^{r}\frac{E\left( D_{i};u,v\right) \left(
uv-1\right) }{\left( uv\right) ^{a_{i}+1}-1} \\
\, \\
+\left( 1+uv\right) \left[ \sum\limits_{(i,j)\in \frak{R}_{\varphi }}\left(
\frac{uv-(uv)^{a_{i}+1}}{(uv)^{a_{i}+1}-1}\right) \left( \frac{%
uv-(uv)^{a_{j}+1}}{(uv)^{a_{j}+1}-1}\right) -\mathbf{b}(X)\right] \\
\, \\
+\sum\limits_{(i,j,k)\in \frak{Q}_{\varphi }}\left( \frac{uv-(uv)^{a_{i}+1}}{%
(uv)^{a_{i}+1}-1}\right) \left( \frac{uv-(uv)^{a_{j}+1}}{(uv)^{a_{j}+1}-1}%
\right) \left( \frac{uv-(uv)^{a_{k}+1}}{(uv)^{a_{k}+1}-1}\right) +\mathbf{t}%
(X)\smallskip
\end{array}
\label{ESTRING}
\end{equation}
with $\mathbf{b}\left( X\right) ,\,\mathbf{t}\left( X\right) $ as defined in
\emph{\ref{LEMMA1} (ii)}. In particular,
\begin{equation}
\begin{array}{l}
e_{\text{\emph{str}}}\left( X\right) -e\left( D_{\varnothing }^{\circ
}\right) =\sum\limits_{i=1}^{r}\frac{e\left( D_{i}\right) }{a_{i}+1}+2\left[
\sum\limits_{(i,j)\in \frak{R}_{\varphi }}\left( \frac{a_{i}}{a_{i}+1}%
\right) \left( \frac{a_{j}}{a_{j}+1}\right) -\mathbf{b}(X)\right] \\
\, \\
-\sum\limits_{(i,j,k)\in \frak{Q}_{\varphi }}\left( \frac{a_{i}}{a_{i}+1}%
\right) \left( \frac{a_{j}}{a_{j}+1}\right) \left( \frac{a_{k}}{a_{k}+1}%
\right) +\mathbf{t}(X)
\end{array}
\,\smallskip  \label{eSTRING}
\end{equation}
\ \emph{(As we shall see below in \ref{H11-COR}, }$e\left( D_{\varnothing
}^{\circ }\right) =0$\emph{)}.
\end{proposition}

\noindent \textit{Proof.} Using inclusion-exclusion principle (\ref{INC-EX})
for the $E$-polynomial of $D_{J}^{\circ }$, we obtain
\begin{equation}
E\left( D_{J}^{\circ };u,v\right) =E\left( D_{J};u,v\right)
-\sum_{\varnothing \neq J^{\prime }\subseteq I\mathbb{r}J}\,\left( -1\right)
^{\left| J^{\prime }\right| -1}\,E\left( D_{J^{\prime }};u,v\right)
\label{STRAT}
\end{equation}
Formula (\ref{E-STR}) can be rewritten via (\ref{STRAT}) as follows:
\begin{equation*}
\begin{array}{l}
E_{\text{str}}\left( X;u,v\right) \\
\, \\
=\sum\limits_{J\subseteq I}\left( E\left( D_{J};u,v\right)
-\sum\limits_{\varnothing \neq J^{\prime }\subseteq I\mathbb{r}J}\,\left(
-1\right) ^{\left| J^{\prime }\right| -1}\,E\left( D_{J^{\prime }\cup
J};u,v\right) \right) \prod\limits_{j\in J}\left( \frac{uv-1}{\left(
uv\right) ^{a_{j}+1}-1}\right) \\
\, \\
=\sum\limits_{J\subseteq I}E\left( D_{J};u,v\right) \,\prod\limits_{j\in
J}\left( \frac{uv-1}{\left( uv\right) ^{a_{j}+1}-1}-1\right) \\
\, \\
=\sum\limits_{J\subseteq I}E\left( D_{J};u,v\right) \,\prod\limits_{j\in
J}\left( \frac{uv-\left( uv\right) ^{a_{j}+1}}{\left( uv\right) ^{a_{j}+1}-1}%
\right) .
\end{array}
\end{equation*}
Hence,
\begin{equation}
\begin{array}{l}
E_{\text{str}}\left( X;u,v\right) -E\left( D_{\varnothing }^{\circ
};u,v\right) \\
\, \\
=E(\bigcup\limits_{i\in I}D_{i};u,v)+\sum\limits_{\varnothing \neq
J\subseteq I}E\left( D_{J};u,v\right) \prod\limits_{j\in J}\left( \frac{%
uv-\left( uv\right) ^{a_{j}+1}}{\left( uv\right) ^{a_{j}+1}-1}\right) \\
\, \\
=\sum\limits_{i=1}^{r}E(D_{i};u,v)-\sum\limits_{(i,j)\in \frak{R}_{\varphi
}}E(D_{\{i,j\}};u,v)+\sum\limits_{(i,j,k)\in \frak{Q}_{\varphi
}}E(D_{\{i,j,k\}};u,v) \\
\, \\
+\sum\limits_{i=1}^{r}E\left( D_{j};u,v\right) \left( \frac{uv-\left(
uv\right) ^{a_{j}+1}}{\left( uv\right) ^{a_{j}+1}-1}\right) +\sum\limits
_{\substack{ J\subseteq I  \\ \left| J\right| \in \{2,3\}}}E\left(
D_{J};u,v\right) \prod_{j\in J}\left( \frac{uv-\left( uv\right) ^{a_{j}+1}}{%
\left( uv\right) ^{a_{j}+1}-1}\right)
\end{array}
\label{SXESH}
\end{equation}
Since $\left| \frak{R}_{\varphi }\right| =\,\mathbf{b}(X)$, $\left| \frak{Q}%
_{\varphi }\right| =\mathbf{t}(X)$, and
\begin{equation*}
E(D_{\{i,j\}};u,v)=1+uv,\ \forall (i,j)\in \frak{R}_{\varphi },\ \ \
E(D_{\{i,j,k\}};u,v)=1,\ \forall (i,j,k)\in \frak{Q}_{\varphi },\text{ }
\end{equation*}
Formula (\ref{ESTRING}) follows from (\ref{SXESH}), and (\ref{eSTRING}) from
(\ref{ESTRING}) by passing to the limit $u,v\longrightarrow 1.$\hfill $%
\square $

%
%
%

\section{Proof of the Theorem}

\noindent \noindent Theorem \ref{MAIN} will be proved by direct evaluation
of formula (\ref{ESTRING}). For this it is obviously enough to determine the
coefficients of the $E$-polynomials of all exceptional prime divisors, on
the one hand, and those of $E\left( D_{\varnothing }^{\circ };u,v\right) ,$
on the other. Hence, in view of lemma \ref{LEMMA2} and of our explicit
description of a canonical desingularization, what remains to be done is the
study of the coefficients of this ``first summand'' $E\left( D_{\varnothing
}^{\circ };u,v\right) $ which depend exclusively on the intrinsic geometry
around the singularities. We begin with a general proposition being valid in
all dimensions.

\begin{proposition}
Let $\left( X,x\right) $ be an isolated complete intersection singularity of
pure dimension $d\geq 2$ and $(\widetilde{X},\frak{Ex}\left( \varphi \right)
)\overset{\varphi }{\longrightarrow }\left( X,x\right) $ a resolution with
exceptional locus $\mathfrak{Ex}\left( \varphi \right) =\cup
_{i=1}^{r}D_{i}. $ Then the coefficients of the $E$-polynomial
\begin{equation}
\begin{array}{l}
E(\widetilde{X}\mathbb{r}\mathfrak{Ex}\left( \varphi \right) ;u,v)=E\left(
D_{\varnothing }^{\circ };u,v\right) = \\
\  \\
=E\left( X\mathbb{r}\left\{ x\right\} ;u,v\right) =\left( uv\right)
^{d}\,E\left( L;u^{-1},v^{-1}\right)
\end{array}
\label{1FORMEL}
\end{equation}
of $\widetilde{X}\mathbb{r}\mathfrak{Ex}\left( \varphi \right) $ depend on
those of the $E$-polynomial of its link $L,$ and, in fact, only on the Hodge
numbers of the $(d-1)$-cohomology group of $L.$\medskip \newline
If $\left( X,x\right) $ is, in addition, a rational singularity, then
\begin{equation}
\left\{
\begin{array}{l}
E(\widetilde{X}\mathbb{r}\mathfrak{Ex}\left( \varphi \right) ;u,v)=E\left( X%
\mathbb{r}\left\{ x\right\} ;u,v\right) =\left( uv\right) ^{d}-1\,+ \\
\! \\
+\left( -1\right) ^{d}\,\ \left[ \sum\limits_{\substack{ 1\leq p,q\leq d-1
\\ 2\leq p+q\leq d-1}}\ h^{p,q}(H^{d-1}(L,\mathbb{C}))\,\,u^{p}\,v^{q}\right]
+ \\
\! \\
+\left( -1\right) ^{d-1}\,\ \left[ \sum\limits_{\substack{ 1\leq p,q\leq d-1
\\ d+1\leq p+q\leq 2d-2}}h^{d-p,d-q}(H^{d-1}(L,\mathbb{C}))\,\,u^{p}\,v^{q}\,%
\right]
\end{array}
\right.  \label{2FORMEL}
\end{equation}
\end{proposition}

\noindent\textit{Proof. }Let $L=L\left( X,x\right) $ denote the link of the
singularity $\left( X,x\right) $, i.e., the intersection of a closed
neighbourhood of $x$ containing it with a small sphere. $L$ is a
differentiable, compact, oriented manifold of dimension $2d-1$, and there
are isomorphisms:
\begin{equation*}
H^{i+1}\left( X,X\mathbb{r}\left\{ x\right\} ,\mathbb{Q}\right) \cong
H^{i}\left( X\mathbb{r}\left\{ x\right\} ,\mathbb{Q}\right) \cong
H^{i}\left( L,\mathbb{Q}\right) \ .
\end{equation*}
For this reason it is sufficient to consider the natural MHS on the
cohomologies of $L$. Note that
\begin{equation}
h^{p,q}(H^{i}\left( L,\mathbb{C}\right) )=h^{q,p}(H^{i}\left( L,\mathbb{C}%
\right) )  \label{PQP}
\end{equation}
while Poincar\'{e} duality implies (\ref{1FORMEL}) because
\begin{equation*}
h^{p,q}(H^{i}\left( L,\mathbb{C}\right) )=h^{d-p,d-q}(H^{2d-i-1}\left( L,%
\mathbb{C}\right) )
\end{equation*}
equals
\begin{equation}
h^{p,q}(H^{i}\left( L,\mathbb{C}\right) )=h^{p,q}(H^{i}\left( X\mathbb{r}%
\left\{ x\right\} ,\mathbb{C}\right) )=h^{d-p,d-q}(H_{c}^{2d-i}\left( X%
\mathbb{r}\left\{ x\right\} ,\mathbb{C}\right) )  \label{PDUAL}
\end{equation}
For the computation of these dimensions it is therefore enough to assume,
from now on, that $i\leq d$. According to \cite[Cor. (15.9)]
{STEENBRINK-SCRIPT}, the restriction map
\begin{equation*}
H^{i}(\widetilde{X},\mathbb{Q)}\longrightarrow H^{i}(\widetilde{X}\mathbb{r}%
\mathfrak{Ex}\left( \varphi\right) ,\mathbb{Q)}\cong H^{i}\left( L,\mathbb{Q}%
\right)
\end{equation*}
is surjective for $i<d$ and equals the zero-map for $i=d.$ From the induced
exact MHS-sequences
\begin{equation*}
\begin{array}{ll}
0\longrightarrow H_{\mathfrak{Ex}\left( \varphi\right) }^{i}(\widetilde {X},%
\mathbb{Q)}\longrightarrow H^{i}\left( \mathfrak{Ex}\left( \varphi\right) ,%
\mathbb{Q}\right) \longrightarrow H^{i}\left( L,\mathbb{Q}\right)
\longrightarrow0\text{ } & \left( i<d\right) \\
\  & \  \\
0\longrightarrow H_{\mathfrak{Ex}\left( \varphi\right) }^{d}(\widetilde {X},%
\mathbb{Q)}\longrightarrow H^{d}\left( \mathfrak{Ex}\left( \varphi\right) ,%
\mathbb{Q}\right) \longrightarrow0\text{ } & \left( i=d\right)
\end{array}
\end{equation*}
one gets the vanishing of $Gr_{j}^{\mathcal{W}_{\bullet}}(H_{\mathfrak
{Ex}\left( \varphi\right) }^{i}(\widetilde{X},\mathbb{Q)}),$ $j\neq i$, and
of $Gr_{j}^{\mathcal{W}_{\bullet}}(H^{i}\left( L,\mathbb{Q}\right) ),$ for $%
j\geq i-1$ (cf. \cite[Cor. 1.12]{STEENBRINK1}), and consequently, for $i<d,$
$h^{p,q}\left( H^{i}\left( L,\mathbb{C}\right) \right) $ equals
\begin{equation}
\left\{
\begin{array}{ll}
h^{p,q}(H^{i}\left( \frak{Ex}\left( \varphi\right) ,\mathbb{C}\right) ), &
\text{if \ \ }p+q<i \\
h^{p,q}(H^{i}(\mathfrak{Ex}\left( \varphi\right) ,\mathbb{C)}%
)-h^{d-p,d-q}(H^{2d-i}(\mathfrak{Ex}\left( \varphi\right) ,\mathbb{C)}), &
\text{if \ \ }p+q=i \\
0, & \text{if \ \ }p+q>i
\end{array}
\right.  \label{HPQ1}
\end{equation}
(The right-hand side of (\ref{HPQ1}) is therefore independent of the choice
of the resolution). Since $X$ is also a complete intersection, $L$ is $%
\left( d-2\right) $-connected (cf. \cite[Kor. 1.3]{HAMM}), and the local
Lefschetz theorem gives:
\begin{equation}
\begin{array}{l}
H^{i}\left( L,\mathbb{C}\right) \cong\mathbb{C}\text{,\ \ for }i\in\left\{
0,2d-1\right\} , \\
H^{i}\left( L,\mathbb{C}\right) =0\text{ ,\ \ for }i\notin\left\{
0,d-1,d,2d-1\right\} .
\end{array}
\label{HPQ2}
\end{equation}
Thus, for $i\in\left\{ 0,2d-1\right\} $, the only non-zero Hodge numbers are
\begin{equation}
h^{0,0}\left( H^{0}\left( L,\mathbb{C}\right) \right)
=h^{d,d}(H^{2d-1}\left( L,\mathbb{C}\right) )=1\ .  \label{HPQ3}
\end{equation}
By (\ref{HPQ1}), (\ref{HPQ2}) and (\ref{HPQ3}) we deduce {\small
\begin{equation*}
\begin{array}{l}
E\left( L;u,v\right) =\sum_{0\leq p,q\leq d}e^{p,q}(L)\ \,u^{p}\,v^{q}= \\
\, \\
=\sum\limits_{0\leq p,q\leq d}\ \left[ \left( h^{p,q}(H^{0}(L,\mathbb{C}%
))-h^{p,q}(H^{2d-1}(L,\mathbb{C}))\right) \right] \,\,u^{p}\,v^{q}+ \\
\, \\
+\sum\limits_{0\leq p,q\leq d}\ \left[ \left( -1\right) ^{d-1}\,\left(
h^{p,q}(H^{d-1}(L,\mathbb{C}))-h^{p,q}(H^{d}(L,\mathbb{C}))\right) \right]
\,\,u^{p}\,v^{q}= \\
\, \\
=\sum\limits_{0\leq p,q\leq d}\left[ (h^{p,q}(H^{0}(L,\mathbb{C}%
))-h^{p,q}(H^{2d-1}(L,\mathbb{C})))\right] \,u^{p}\,v^{q}+ \\
\, \\
+\sum\limits_{0\leq p,q\leq d}\left[ \left( -1\right)
^{d-1}(h^{p,q}(H^{d-1}(L,\mathbb{C}))-h^{d-p,d-q}(H^{d-1}(L,\mathbb{C})))%
\right] \,u^{p}\,v^{q}= \\
\, \\
=1-\left( uv\right) ^{d}+\left( -1\right) ^{d-1}\,[\sum\limits_{0\leq
p,q\leq d}\ h^{p,q}(H^{d-1}(L,\mathbb{C}))\,\,u^{p}\,v^{q}]+ \\
\, \\
+\left( -1\right) ^{d}\,[-\sum\limits_{0\leq p,q\leq d}\
h^{d-p,d-q}(H^{d-1}(L,\mathbb{C}))\,\,u^{p}\,v^{q}]= \\
\, \\
=1-\left( uv\right) ^{d}+\left( -1\right) ^{d-1}\left[ \sum \limits
_{\substack{ 0\leq p,q\leq d-1  \\ 0\leq p+q\leq d-1}}\ h^{p,q}(H^{d-1}(L,%
\mathbb{C}))\,\,u^{p}\,v^{q}\right] + \\
\, \\
+\left( -1\right) ^{d}\left[ \sum\limits_{\substack{ 1\leq p,q\leq d  \\ %
d+1\leq p+q\leq2d-1}}\ h^{d-p,d-q}(H^{d-1}(L,\mathbb{C}))\,\,u^{p}\,v^{q}%
\right]
\end{array}
\end{equation*}
} which proves the first assertion. Now setting
\begin{equation*}
\ell^{p,q}(L):=\text{ dim}_{\mathbb{C}}Gr_{\mathcal{F}^{\bullet}}^{p}\left(
H^{p+q}\left( L,\mathbb{C}\right) \right) ,
\end{equation*}
one has
\begin{equation*}
\ell^{p,q}(L)=\text{dim}_{\mathbb{C}}H^{q}\left( \mathfrak{Ex}\left(
\varphi\right) ,\Omega_{\widetilde{X}}^{p}\left( \text{log }\mathfrak {Ex}%
\left( \varphi\right) \right) \otimes\mathcal{O}_{\mathfrak{Ex}\left(
\varphi\right) }\right)
\end{equation*}
(cf. \cite[\S1]{STEENBRINK1} and \cite[\S3]{STEENBRINK3}). Obviously,
\begin{equation*}
\ell^{p,i-p}(L)=\sum_{q=0}^{d}h^{p,q}\left( H^{i}\left( L,\mathbb{C}\right)
\right)
\end{equation*}
for $i\geq p.$ If $\left( X,x\right) $ is, in addition, a \textit{rational}
singularity, then for all $i\geq1$ we have
\begin{equation}
\ell^{0,i}(L)=\text{dim}_{\mathbb{C}}H^{i}\left( \mathfrak{Ex}\left(
\varphi\right) ,\mathcal{O}_{\mathfrak{Ex}\left( \varphi\right) }\right)
=0=\ell^{i,0}(L)  \label{VANISH}
\end{equation}
because $\ell^{i,0}(L)\leq\ell^{0,i}(L),\,H^{i}(\widetilde{X},\mathcal{O}_{%
\widetilde{X}})=0$ and
\begin{equation*}
H^{i}(\widetilde{X},\mathcal{O}_{\widetilde{X}})\longrightarrow H^{i}\left( %
\mathfrak{Ex}\left( \varphi\right) ,\mathcal{O}_{\mathfrak{Ex}\left(
\varphi\right) }\right)
\end{equation*}
is surjective by \cite[Lemma 2.14]{STEENBRINK1}. Hence,
\begin{equation}
h^{j,0}\left( H^{i}\left( L,\mathbb{C}\right) \right) \ \overset{(\ref{PQP})%
}{=}\ h^{0,j}\left( H^{i}\left( L,\mathbb{C}\right) \right) \ \overset{(\ref
{VANISH})}{=}0,\ \ \ \text{for}\ 0\leq j\leq d\text{ and }i\geq1.
\label{VAN}
\end{equation}
This means that the $E$-polynomial of $L$ can be written as
\begin{equation}
\left\{
\begin{array}{l}
E\left( L;u,v\right) =1-\left( uv\right) ^{d}+ \\
\, \\
+\left( -1\right) ^{d-1}\left[ \sum\limits_{\substack{ 1\leq p,q\leq d-1  \\ %
2\leq p+q\leq d-1}}\ h^{p,q}(H^{d-1}(L,\mathbb{C}))\,\,u^{p}\,v^{q}\right] +
\\
\, \\
+\left( -1\right) ^{d}\,\left[ \sum\limits_{\substack{ 1\leq p,q\leq d-1  \\ %
d+1\leq p+q\leq2d-2}}\ h^{d-p,d-q}(H^{d-1}(L,\mathbb{C}))\,\,u^{p}\,v^{q}\,%
\right]
\end{array}
\right.  \label{3FORMEL}
\end{equation}
and formula (\ref{2FORMEL}) follows from (\ref{3FORMEL}) and (\ref{1FORMEL}%
).\hfill$\square$

\begin{remark}
(i) Let us now denote by $F_{f}$ the \textit{Milnor fiber} being associated
to the $\mathbf{A}$-$\mathbf{D}$-$\mathbf{E}$ singularity $(X_{f}^{\left(
d\right) },\mathbf{0).}$ As it is known (cf. \cite[Thm. 6.5]{MILNOR}), $%
F_{f} $ has the homotopy type of a bouquet of $d$-spheres, and its \textit{%
Milnor number}
\begin{equation*}
\mu \left( f\right) :=\mu \left( F_{f}\right) :=\text{ }\#\left\{ \text{of
these spheres}\right\} =\text{ dim}_{\mathbb{C}}(\mathcal{O}_{d+1}\,/\,(%
\tfrac{\partial f}{\partial x_{1}},\ldots ,\tfrac{\partial f}{\partial
x_{d+1}}))
\end{equation*}
is in each case equal to the subscript of the type under consideration.
According to the \textit{Sebastiani-Thom theorem} \cite{SEB-T} (see also
\cite[pp. 86-88]{DIMCA2}), the splitting $f=g+g^{\prime }$ (as in (\ref{HYP}%
)) gives rise to the construction of an homotopy equivalence between the
Milnor fiber $F_{f}$ and the join $F_{g}\ast F_{g^{\prime }}$ of the
corresponding Milnor fibers $F_{g}$ and $F_{g^{\prime }}.$ In particular,
this implies
\begin{equation}
\mu \left( f\right) =\mu \left( g\right) \cdot \mu \left( g^{\prime }\right)
=\mu \left( g\right)  \label{MILNN}
\end{equation}
(ii) For any isolated complete intersection singularity $(X,x)$ of pure
dimension $d,$ with link $L,$ Milnor fiber $F$ and Milnor number $\mu \left(
F\right) $, \textit{Steenbrink's invariant}
\begin{equation*}
s_{j}(X,x),\ \ 0\leq j\leq d,
\end{equation*}
is defined in \cite{STEENBRINK2} by regarding any $1$-parameter smoothing $%
\psi :\left( \frak{X},x\right) \rightarrow \left( \mathbb{C},0\right) $ of $%
(X,x)$ (with $\frak{X}_{0}=\psi ^{-1}\left( 0\right) \cong X$) and setting
\begin{equation*}
s_{j}(X,x):=\text{ dim}_{\mathbb{C}}Gr_{\mathcal{F}^{\bullet }}^{j}\,\mathbb{%
H}^{d}\left( \Phi _{\psi }^{\bullet }\left( \mathbb{C}\right) \right) ,
\end{equation*}
where $\mathcal{F}^{\bullet }$ denotes here the Hodge-filtration of the
highest hypercohomology group of the complex $\Phi _{\psi }^{\bullet }\left(
\mathbb{C}\right) $ of sheaves of vanishing cycles associated to $\psi .$
(For all $q,$ the direct image sheaves $\Phi _{\psi }^{q}\left( \mathbb{C}%
\right) =R^{q}\left( \vartheta _{t}\right) _{\ast }\underline{\mathbb{C}}_{%
\frak{X}}$ are defined on $\frak{X}_{0},$ with $\vartheta _{t}:\frak{X}%
_{t}\rightarrow \frak{X}_{0}$ denoting the restriction of the retraction $%
\vartheta :\frak{X}\rightarrow \frak{X}_{0}$ onto a fiber $\frak{X}_{t}.$ In
fact, the definition of $\Phi _{\psi }^{q}\left( \mathbb{C}\right) $ can be
made independent of the choice of the fiber $\frak{X}_{t}$ by passing to the
``canonical'' fiber $\frak{X}_{\infty }$ of $\psi .$ In this setting, the
fiber of the sheaf $\Phi _{\psi }^{q}\left( \mathbb{C}\right) $ over $x$ is
isomorphic to $\widetilde{H}^{q}\left( \frak{X}_{t,x},\mathbb{C}\right) ,$
where $\frak{X}_{t,x}$ is diffeomorphic to the Milnor fiber $F$). $%
s_{j}(X,x) $ is an upper semicontinuous invariant under deformations of $%
(X,x),$ \textit{does not} depend on the particular choice of $\psi $ (cf.
\cite[(1.8)-(1.10), and (2.6)]{STEENBRINK2})$,$ and
\begin{equation}
\mu \left( F\right) =s_{0}(X,x)+s_{1}(X,x)+\cdots +s_{d-1}(X,x)+s_{d}(X,x)
\label{MILN-S}
\end{equation}
On the other hand, taking into account the $\mathbb{Q}(-d)$-duality between $%
H^{d}(F,L,\mathbb{C})$ and $H^{d}(F,\mathbb{C}),$ and the exact MHS-sequence
\begin{equation*}
0\longrightarrow H^{d-1}(L,\mathbb{C)}\longrightarrow H^{d}(F,L,\mathbb{C}%
)\longrightarrow H^{d}(F,\mathbb{C})\longrightarrow H^{d}(L,\mathbb{C)}%
\longrightarrow 0,
\end{equation*}
one deduces the equalities
\begin{equation}
\begin{array}{l}
s_{j}(X,x)-s_{d-j}(X,x)= \\
\  \\
=\ell ^{j,d-j}\left( L\right) -\ell ^{j,d-j-1}(L)=\ell ^{d-j,j-1}\left(
L\right) -\ell ^{j,d-j-1}(L)
\end{array}
\label{SLE}
\end{equation}
\end{remark}

\begin{corollary}
\label{H11-COR}Let $X=X_{f}^{\left( 3\right) }$ be the underlying spaces of
the three-dimensional $\mathbf{A}$-$\mathbf{D}$-$\mathbf{E}$ singularities.
Then we have\emph{\ }
\begin{equation}
E\left( X\mathbb{r}\left\{ \mathbf{0}\right\} ;u,v\right) =\left(
uv-1\right) \,\left[ 1+\left( 1+h^{1,1}\left( H^{2}(L,\mathbb{C})\right)
\right) \,uv+\left( uv\right) ^{2}\right]  \label{H11-FOR}
\end{equation}
where
\begin{equation*}
\begin{tabular}{|c|c|c|c|c|c|}
\hline
$\mathbf{Types}$ & $\mathbf{A}_{n}$ & $\mathbf{D}_{n}$ & $\mathbf{E}_{6}$ & $%
\mathbf{E}_{7}$ & $\mathbf{E}_{8}$ \\ \hline\hline
$h^{1,1}\left( H^{2}(L,\mathbb{C})\right) $ & $
\begin{array}{c}
\, \\
\left\{
\begin{array}{l}
1,\ \text{\emph{for}\textit{\ }}n\text{\textit{\ }\emph{odd}} \\
0,\ \text{\emph{for}\textit{\ }}n\text{\textit{\ }\emph{even}}
\end{array}
\right. \\
\,
\end{array}
$ & $\left\{
\begin{array}{l}
1,\ \text{\emph{for}\textit{\ }}n\text{\textit{\ }\emph{odd}} \\
2,\ \text{\emph{for}\textit{\ }}n\text{\textit{\ }\emph{even}}
\end{array}
\right. $ & $0$ & $1$ & $0$ \\ \hline
\end{tabular}
\end{equation*}
\end{corollary}

\noindent\textit{Proof}. Formula (\ref{H11-FOR}) is nothing but (\ref
{2FORMEL}) for $d=3.$ So it remains to compute $h^{1,1}\left( H^{2}(L,%
\mathbb{C})\right) .$ Using the notation $\mu\left( f\right) :=\mu\left(
F_{f}\right) $ and$\ s_{j}(f):=s_{j}(X,\mathbf{0})$ for the singularity $(X,%
\mathbf{0)},$ the equalities (\ref{VANISH}), (\ref{VAN}) and (\ref{SLE})
give
\begin{equation}
\ell^{1,1}(L)=h^{1,1}\left( H^{2}(L,\mathbb{C})\right) =s_{2}(f\mathbf{)}%
-s_{1}(f\mathbf{)}  \label{H11-FORMULA}
\end{equation}
and $s_{0}(f)=s_{3}(f).$ Furthermore, by (\ref{MILN-S}),
\begin{equation*}
\mu\left( f\right)
=s_{0}(f)+s_{1}(f)+s_{2}(f)+s_{3}(f)=s_{1}(f)+s_{2}(f)+2s_{3}(f)\,.
\end{equation*}
In fact, since $(X,\mathbf{0})$ is a \textit{Du Bois singularity} (as it is
a rational isolated singularity)$,$ or equivalently, since $s_{3}(f)$ equals
the \textit{geometric genus} of $(X,\mathbf{0})$ (see \cite[\S4]{STEENBRINK3}%
, \cite[(2.17) and (3.7)]{STEENBRINK1})$,$ we have $s_{0}(f)=s_{3}(f)=0,$
i.e., $\mu\left( f\right) =s_{0}(f)+s_{1}(f).$ Now the splitting $%
f=g+g^{\prime}$ (as in (\ref{HYP})) leads to a ``Sebastiani-Thom formula''
for Steenbink's invariant; namely,
\begin{equation}
s_{j}(f)=s_{j-1}(g)  \label{SEB-T-ST}
\end{equation}
Applying Milnor's formula \cite[Thm. 10.5]{MILNOR} for the curve singularity
$(X_{g},\mathbf{0}),$ we obtain
\begin{equation}
\mu\left( g\right) =2\,\delta\left( g\right) -r\left( g\right) +1
\label{MILN-FOR}
\end{equation}
where
\begin{equation*}
r\left( g\right) :=\#\{\text{branches of the curve }X_{g}\text{ passing
through the origin}\}
\end{equation*}
and
\begin{equation*}
\delta\left( g\right) :=\#\{``\text{virtual'' double points w.r.t. }X_{g}\}=%
\text{ dim}_{\mathbb{C}}(\nu_{\ast}\mathcal{O}_{\widetilde{X_{g}}}\,/\,%
\mathcal{O}_{X_{g}})
\end{equation*}
with $\nu:\widetilde{X}_{g}\longrightarrow X_{g}$ the normalization of $%
X_{g} $. Note that this first number $r\left( g\right) $ is directly
computable because the only types for which $g\left( x_{1},x_{2}\right) $'s
are reducible, are $\mathbf{A}_{n}$'s$,$ for $n$ odd, with
\begin{equation*}
g\left( x_{1},x_{2}\right) =(x_{1}^{\frac{n+1}{2}}+\sqrt{-1}\
x_{2})\,(x_{1}^{\frac{n+1}{2}}-\sqrt{-1}\ x_{2}),
\end{equation*}
$\mathbf{D}_{n}$'s with
\begin{equation*}
g\left( x_{1},x_{2}\right) =\left\{
\begin{array}{ll}
x_{1}\,(x_{1}^{n-2}+x_{2}^{2}), & \text{if \ }n\ \text{ is odd} \\
\, & \, \\
x_{1}\,(x_{1}^{\frac{n}{2}-1}+\sqrt{-1}\ x_{2})\,(x_{1}^{\frac{n}{2}-1}-%
\sqrt{-1}\ x_{2}), & \text{if \ }n\ \text{ is even}
\end{array}
\right.
\end{equation*}
and $\mathbf{E}_{7}$ with
\begin{equation*}
g\left( x_{1},x_{2}\right) =x_{1}(x_{1}^{2}+x_{2}^{3})\ ,
\end{equation*}
while $\delta\left( g\right) $ can be read off from (\ref{MILN-FOR}) via the
Milnor number. Finally, since
\begin{equation}
s_{1}(f)\,\overset{(\ref{SEB-T-ST})}{=}\,s_{0}(g)=\delta\left( g\right)
-r\left( g\right) +1,\ \ \ s_{2}\left( f\right) \,\overset{(\ref{SEB-T-ST})}{%
=}\,s_{1}(g)=\delta\left( g\right) ,  \label{S1-2}
\end{equation}
(cf. \cite[(2.17), p. 526]{STEENBRINK1}), we may form the following
table:\bigskip\smallskip{\small
\begin{equation*}
\begin{tabular}{|c|c|c|c|c|}
\hline
\textbf{Types} & $\mu\left( f\right) =\mu\left( g\right) $ & $\ \ \ r\left(
g\right) $ & $s_{1}(f)=s_{0}(g)$ & $s_{2}\left( f\right)
=s_{1}(g)=\delta\left( g\right) $ \\ \hline\hline
$\mathbf{A}_{n},\ n\text{ odd}$ & $n$ & $
\begin{array}{c}
\text{\thinspace} \\
2 \\
\text{\thinspace}
\end{array}
$ & $\frac{n-1}{2}$ & $\frac{n+1}{2}$ \\ \hline
$\mathbf{A}_{n},\ n\text{ even}$ & $
\begin{array}{c}
\  \\
n \\
\
\end{array}
$ & $1$ & $\frac{n}{2}$ & $\frac{n}{2}$ \\ \hline
$\mathbf{D}_{n},\ n\text{ odd}$ & $n$ & $
\begin{array}{c}
\text{\thinspace} \\
2 \\
\text{\thinspace}
\end{array}
$ & $\frac{n-1}{2}$ & $\frac{n+1}{2}$ \\ \hline
$\mathbf{D}_{n},\ n\text{ even}$ & $
\begin{array}{c}
\  \\
n \\
\
\end{array}
$ & $3$ & $\frac{n-2}{2}$ & $\frac{n+2}{2}$ \\ \hline
$\mathbf{E}_{6}$ & $6$ & $1$ & $3$ & $3$ \\ \hline
$\mathbf{E}_{7}$ & $7$ & $2$ & $3$ & $4$ \\ \hline
$\mathbf{E}_{8}$ & $8$ & $1$ & $4$ & $4$ \\ \hline
\end{tabular}
\medskip
\end{equation*}
\smallskip} This table allows us to evaluate $h^{1,1}\left( H^{2}(L,\mathbb{C%
})\right) $ for all possible types via (\ref{S1-2}) and (\ref{H11-FORMULA}%
).\hfill$\square$\bigskip\medskip

\noindent \textbf{Proof of Theorem \ref{MAIN}: }It follows directly from the
explicit arithmetical data for each of the canonical resolutions given in
Lemma \ref{LEMMA1} and Proposition \ref{DISCR}, and from formulae (\ref
{ESTRING}), (\ref{eSTRING}), in combination with the formula (\ref{H11-FOR})
of Corollary \ref{H11-COR}.\hfill $\square $

\begin{final remarks and questions}
(i) Is the resolution algorithm (or a slight modification of it) extendible
to a wider class of three-dimensional Gorenstein terminal (or canonical)
singularities\thinspace ?\smallskip\ \newline
(ii) The $d$-dimensional generalization of Theorem \ref{MAIN} seems to be
feasible as the pattern of the local reduction of simple singularities
remains invariant (after all, adding quadratic terms does not cause very
crucial changes in the desingularization procedure), though the
investigation of the structure of the corresponding exceptional prime
divisors and of their intersections for the $\mathbf{D}$\textbf{-}$\mathbf{E}
$'s might be rather complicated\textbf{.\smallskip\ }\newline
(iii) Since the string-theoretic ``adjusting property'' of $E_{\text{str}}$%
-functions is of local nature and focuses solely on the \textit{singular loci%
} of the varieties being under consideration, it is clear how to treat of $%
E_{\text{str}}$ and $e_{\text{str}}$ in global geometric constructions with
prescribed $\mathbf{A}$-$\mathbf{D}$\textbf{-}$\mathbf{E}$-singularities. We
close the paper by giving some examples of this sort.
\end{final remarks and questions}

%
%
%

\section{Global geometric applications}

\noindent In view of Theorem \ref{MAIN}, the $E_{\text{str}}$-function of a
complex threefold $Y$ having only $\mathbf{A}$-$\mathbf{D}$\textbf{-}$%
\mathbf{E}$-singularities $q_{1},q_{2},..,q_{k}$ is computable provided that
one knows how to determine the Hodge numbers $h^{p,q}(H_{c}^{i}(Y,\mathbb{C}%
))$ of $Y,$ as we obtain:
\begin{align}
E_{\text{str}}\left( Y;u,v\right) & =E\left( Y\mathbb{r}\{q_{1},q_{2},%
\ldots,q_{k}\};u,v\right) \ +\ \sum_{i=1}^{k}\ E_{\text{str}}\left( \left(
Y,q_{i}\right) ;u,v\right) =  \notag \\
& =E\left( Y;u,v\right) \ +\ \sum_{i=1}^{k}\ \left( E_{\text{str}}\left(
\left( Y,q_{i}\right) ;u,v\right) -1\right)  \label{ESTR-GLOBAL}
\end{align}

\medskip

\noindent{}\textsf{(a)} \textbf{Complete intersections in projective spaces.}
A very simple closed formula for $e_{\text{str}}$ can be built whenever $\ Y$
is a (global) complete intersection in a projective space.

\begin{proposition}
Let $Y=Y_{\left( d_{1},d_{2},\ldots ,d_{r-3}\right) }$ $\ $be a
three-dimensional complete intersection of multidegree $\left(
d_{1},d_{2},\ldots ,d_{r-3}\right) $ in $\mathbb{P}_{\mathbb{C}}^{r}$ having
only $k$ isolated singularities $q_{1},q_{2},\ldots ,q_{k}$ of type $\mathbf{%
A}$-$\mathbf{D}$\textbf{-}$\mathbf{E}.$ Then its string-theoretic Euler
number equals
\begin{equation}
\begin{array}{l}
e_{\text{\emph{str}}}(Y)=\left[ \binom{r+1}{3}+\sum\limits_{\nu
=1}^{3}(-1)^{\nu }\binom{r+1}{3-\nu }\left( \sum\limits_{1\leq j_{1}\leq
..\leq j_{\nu }\leq r-3}d_{j_{1}}\,\cdots \,d_{j_{\nu }}\right) \right]
\left( \prod\limits_{j=1}^{r-3}\ d_{j}\right) + \\
\  \\
+\ \ \sum\limits_{i=1}^{k}\ \left[ e_{\text{\emph{str}}}\left(
Y,q_{i}\right) +\mu \left( Y,q_{i}\right) -1\right]
\end{array}
\label{ESTR-CI}
\end{equation}
where $\mu \left( Y,q_{i}\right) $ is the Milnor number of the singularity $%
\left( Y,q_{i}\right) $ and $e_{\text{\emph{str}}}\left( Y,q_{i}\right) $
can be read off from the Theorem\emph{\ \ref{MAIN}.}
\end{proposition}

\noindent\textit{Proof}. Considering a small deformation of $Y$ one can
always obtain a non-singular complete intersection $Y^{\prime}$ in $\mathbb{P%
}_{\mathbb{C}}^{r}$ having multidegree $\left( d_{1},d_{2},\ldots
,d_{r-3}\right) $. If we take a ball $B_{i}$ in $\mathbb{P}_{\mathbb{C}}^{r}$
centered at the point $q_{i},$ then, choosing $B_{i}$ small enough, $%
\overline{B_{i}}\cap Y$ is contractible and $\overline{B_{i}}\cap Y^{\prime}
$ can be identified with the (closed) Milnor fiber of the singularity $%
\left( Y,q_{i}\right) .$ $\widehat{Y}:=Y\mathbb{r(}\bigcup_{i=1}^{k}B_{i}%
\mathbb{)}$ and $\widehat{Y}^{\prime}:=Y^{\prime}\mathbb{r(}%
\bigcup_{i=1}^{k}B_{i}\mathbb{)}$ are homeomorphic. Therefore $e(\widehat{Y}%
)=e(\widehat {Y}^{\prime})$. Using Mayer-Vietoris sequence for the splitting
$Y=\widehat {Y}\cup\bigcup_{i=1}^{k}(\overline{B_{i}}\cap Y),$ on the one
hand, and for the splitting $Y^{\prime}=\widehat{Y}^{\prime}\cup%
\bigcup_{i=1}^{k}(\overline{B_{i}}\cap Y^{\prime}),$ on the other, we get $%
e\left( Y\right) =e(\widehat{Y})+k$ and
\begin{equation*}
e(Y^{\prime})\ =\ e(\widehat{Y}^{\prime})\ +\ k\ -\ \sum\limits_{i=1}^{k}\
\mu\left( Y,q_{i}\right) ,
\end{equation*}
respectively (see \cite[Ch. 5, Cor. 4.4 (ii), p. 162 ]{DIMCA2}). Hence,
\begin{equation*}
e\left( Y\right) =e\left( Y^{\prime}\right) \ +\ \sum\limits_{i=1}^{k}\
\mu\left( Y,q_{i}\right) \ .
\end{equation*}
The Euler number of $Y^{\prime}$ can be computed in terms of its multidegree
data either by determining the $\chi_{y}$-characteristic of $Y^{\prime}$ via
Riemann-Roch Theorem (see Hirzebruch \cite[\S2]{HIRZEBRUCH1}) or directly by
Gauss-Bonnet theorem, i.e., by evaluating the highest Chern class of $%
Y^{\prime}$ at its fundamental cycle (cf. \cite[p. 416]{GR-H} \& Chen-Ogiue
\cite[Thm. 2.1]{CH-O}), and is expressible by the closed formula:
\begin{equation*}
\begin{array}{c}
e\left( Y^{\prime}\right) =\left[ \binom{r+1}{3}+\sum\limits_{\nu=1}^{3}%
\,(-1)^{\nu}\,\binom{r+1}{3-\nu}\,\left( \sum\limits_{1\leq
j_{1}\leq\cdots\leq j_{\nu}\leq r-3}d_{j_{1}}\,\cdots\,d_{j_{\nu}}\right)
\right] \ \left( \prod\limits_{j=1}^{r-3}\ d_{j}\right) .
\end{array}
\end{equation*}
Now (\ref{ESTR-CI}) follows clearly from (\ref{ESTR-GLOBAL}).\hfill $%
\square\medskip$

\begin{examples}
(i) If $Y$ possesses only $\mathbf{A}_{1}$-singularities (i.e., ``ordinary
double points'' or ``nodes''), then the second summand in (\ref{ESTR-CI})
equals $2\,\#($nodes of $Y).$ Let us apply (\ref{ESTR-CI}) for some
well-known hypersurfaces $Y$ in $\mathbb{P}_{\mathbb{C}}^{4}$ with \textit{%
many} nodes. [$e_{\text{str}}(Y)$ is nothing but the Euler number of the
overlying spaces of the so-called (simultaneous) ``small resolutions'' of
the nodes of $Y$'s.].\medskip \newline
$\blacktriangleright $ \textbf{Schoen's quintic} \cite{SCHOEN1}\textbf{. }%
This is the quintic
\begin{equation*}
Y=\left\{ (z_{1}:\,\ldots \,:z_{5})\in \mathbb{P}_{\mathbb{C}}^{4}\ \left| \
\,\sum\limits_{i=1}^{5}\,z_{i}^{5}\,-\,5\,\prod_{i=1}^{5}\,z_{i}=0\right.
\right\}
\end{equation*}
having $125$ nodes, namely the members of the orbit of point $\left(
1:1:1:1:1\right) $ under the action of the group which is generated by the
coordinate transformations
\begin{equation*}
(z_{1}:\,\ldots \,:z_{5})\longmapsto (z_{1}:\zeta _{5}^{\alpha
_{1}}\,z_{2}:\,\ldots \,:\zeta _{5}^{\alpha _{4}}\,z_{5}),
\end{equation*}
where $\zeta _{5}=e^{\frac{2\pi \sqrt{-1}}{5}}$, $\sum\limits_{j=1}^{4}%
\alpha _{j}\equiv 0($mod $5).$ Hence, $e_{\text{str}}(Y)=-200+2\cdot
125=\allowbreak 50.\medskip $ \newline
$\blacktriangleright $ \textbf{Hirzebruch's quintic} \cite{HIRZEBRUCH2}%
\textbf{. }Let $\{\boldsymbol
{\Phi}(x,y)=\prod_{i=1}^{5}\Phi _{i}(x,y)=0\}$ be the equation of the curve
of degree $5$ in the real $(x,y)$-plane constructed by the five lines $\Phi
_{i}(x,y)=0,1\leq i\leq 5,$ of a regular pentagon:
\begin{figure}[h]
\begin{center}
\includegraphics[width=1.6509in,height=1.631in]{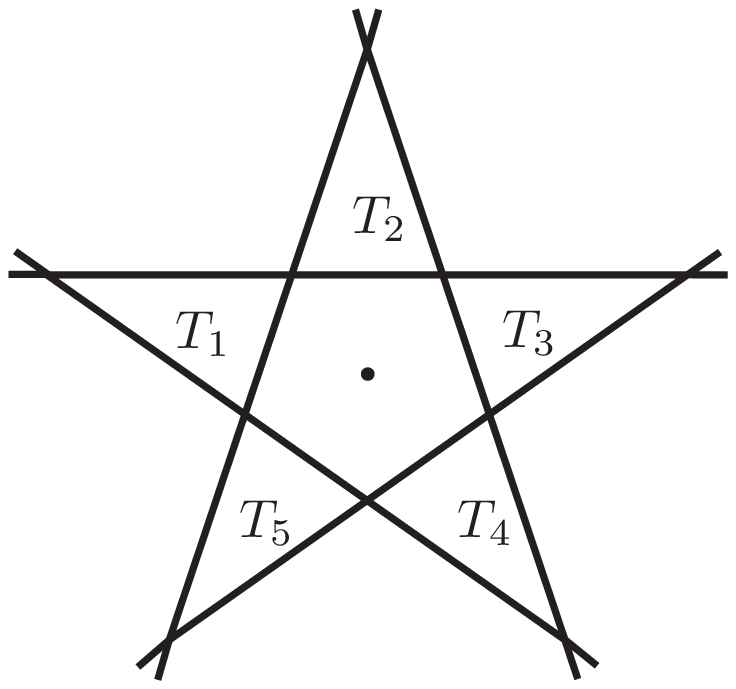}
\end{center}
\end{figure}
\newline
This real picture shows that both partial derivatives of $\boldsymbol{\Phi}$
vanish at the $10$ points of line intersections, as well as at one point $%
t_{i}$ at every triangle $T_{i}$ and at the center of the pentagon.
Moreover, by symmetry, one has $\boldsymbol{\Phi}(t_{i})=\boldsymbol{\Phi}%
(t_{j})$ for all $1\leq i\leq j\leq 5.$ The hypersurface $Y\subset $ $%
\mathbb{P}_{\mathbb{C}}^{4}$ obtained after homogenization of the
three-dimensional affine complex variety
\begin{equation*}
\left\{ \left( z_{1},z_{2},z_{3},z_{4}\right) \in \mathbb{C}^{4}\ \left| \
\right. \boldsymbol{\Phi}(z_{1},z_{2})-\boldsymbol{\Phi}(z_{3},z_{4})=0%
\right\}
\end{equation*}
has $10^{2}+5^{2}+1^{2}=126$ nodes. This means that $e_{\text{str}%
}(Y)=-200+2\cdot 126=\allowbreak 52.$\medskip \newline
$\blacktriangleright $ \textbf{Symmetric Hypersurfaces. }In $\mathbb{P}_{%
\mathbb{C}}^{5}$ with $(z_{1}:\,\ldots \,:z_{6})$ as homogeneous coordinates
we define the threefolds
\begin{equation*}
\left\{
\begin{array}{l}
Y_{1}:=\{(z_{1}:\,\ldots \,:z_{6})\in \mathbb{P}_{\mathbb{C}}^{5}\ \left| \
\sigma _{1}(z_{1},\ldots
,z_{6})=\,\sum\nolimits_{i=1}^{6}\,z_{i}^{3}=0\right. \}, \\
\, \\
Y_{2}:=\left\{ (z_{1}:\,\ldots \,:z_{6})\in \mathbb{P}_{\mathbb{C}}^{5}\
\left| \ \sigma _{1}(z_{1},\ldots ,z_{6})=\,\right. \sigma _{4}(z_{1},\ldots
,z_{6})=0\right\} , \\
\, \\
Y_{3}:=\left\{ (z_{1}:\,\ldots \,:z_{6})\in \mathbb{P}_{\mathbb{C}}^{5}\
\left|
\begin{array}{l}
\ \sigma _{1}(z_{1},\ldots ,z_{6})=\sigma _{5}(z_{1},\ldots ,z_{6})+ \\
+\sigma _{2}(z_{1},\ldots ,z_{6})\,\sigma _{3}(z_{1},\ldots ,z_{6})=0
\end{array}
\right. \right\} ,
\end{array}
\right.
\end{equation*}
where
\begin{equation*}
\sigma _{j}(z_{1},\ldots ,z_{6})=\sum_{1\leq \kappa _{1}<\kappa _{2}<\cdots
<\kappa _{j}\leq 6}\ z_{\kappa _{1}}\cdot z_{\kappa _{2}}\cdot \cdots \cdot
z_{\kappa _{j}},\ \ 1\leq j\leq 6,
\end{equation*}
denote the elementary symmetric polynomials with respect to the variables $%
z_{1},..,z_{6}.$ Obviously, $Y_{i}$'s are invariant under the symmetry group
$\mathfrak{S}_{6}$ acting on $\mathbb{P}_{\mathbb{C}}^{5}$ by permuting
coordinates. Moreover, since the first equation
\begin{equation*}
\sigma _{1}(z_{1},\ldots ,z_{6})=\,z_{1}+z_{2}+z_{3}+z_{4}+z_{5}+z_{6}=0
\end{equation*}
is linear, $Y_{i}$'s can be thought of as hypersurfaces in
\begin{equation*}
\mathbb{P}_{\mathbb{C}}^{4}=\left\{ (z_{1}:\,\ldots \,:z_{6})\in \mathbb{P}_{%
\mathbb{C}}^{5}\ \left| \ \sigma _{1}(z_{1},\ldots ,z_{6})=0\right. \right\}
.
\end{equation*}
The threefold $Y_{1}$ has $10$ nodes, namely the points of $\mathbb{P}_{%
\mathbb{C}}^{5}$ for which three of their coordinates are $1$ and the other
three are $-1$ (i.e., just the members of the $\mathfrak{S}_{6}$-orbit of $%
(1:1:1:-1:-1:-1)$). Correspondingly, $Y_{2}$ has $45$ nodes, and $Y_{3}$ has
$130$ nodes, $10$ constituting the $\mathfrak{S}_{6}$-orbit of $%
(1:1:1:-1:-1:-1),$ $90$ in the $\mathfrak{S}_{6}$-orbit of $(1:1:-1:-1:\sqrt{%
-3}:-\sqrt{-3})$ and $30$ more in the $\mathfrak{S}_{6}$-orbit of $(1:1:1:1:%
\sqrt{-3}-2:-\sqrt{-3}-2).$ The following table gives their special names,
their string-theoretic Euler numbers, as well as the main references for
further reading about their geometric properties. (Note that $Y_{1}$ and $%
Y_{2}$ attain exactly the upper bound for the cardinality of nodes for
cubics and quartics in $\mathbb{P}_{\mathbb{C}}^{4}$, respectively. $Y_{3}$
is, to the best of our knowledge, the quintic in $\mathbb{P}_{\mathbb{C}%
}^{4} $ with the largest known number of nodes).\medskip\
\begin{equation*}
\begin{tabular}{|c|c|c|c|}
\hline
\textbf{Threefolds} & \textbf{Name} & Ref. & $e_{\text{str}}$ \\ \hline\hline
$Y_{1}$ & \textbf{Segre's cubic} & \cite{SEGRE2} & $-6+2\cdot 10=14$ \\
\hline
$Y_{2}$ & \textbf{Burkhart's quartic} & \cite{BURKHARDT}, \cite{FINKELBERG}
& $-56+2\cdot 45=34$ \\ \hline
$Y_{3}$ & \textbf{van Straten's quintic} & \cite{VANSTRATEN} & $-200+2\cdot
130=60$ \\ \hline
\end{tabular}
\end{equation*}
\bigskip \newline
(ii) Let now $Y_{1},Y_{2}$ be the three-dimensional complete intersections
of two quadrics\smallskip\
\begin{equation*}
Y_{i}:=\left\{ \mathbf{z}=(z_{1}:z_{2}:..:z_{6})\in \mathbb{P}_{\mathbb{C}%
}^{5}\ \left| \ ^{t}\mathbf{z}\ M_{i}\ \mathbf{z}=\ ^{t}\mathbf{z}\
M_{i}^{\prime }\ \mathbf{z\ }=\ 0\right. \right\} ,\ i=1,2,
\end{equation*}
where
\begin{equation*}
\begin{array}{ccc}
M_{1}=\left(
\begin{array}{cccccc}
0 & 0 & 0 & 0 & 0 & 1 \\
0 & 0 & 0 & 0 & 1 & 0 \\
0 & 0 & 0 & 1 & 0 & 0 \\
0 & 0 & 1 & 0 & 0 & 0 \\
0 & 1 & 0 & 0 & 0 & 0 \\
1 & 0 & 0 & 0 & 0 & 0
\end{array}
\right) \ , & \  & M_{1}^{\prime }=\left(
\begin{array}{cccccc}
0 & 0 & 0 & 0 & 0 & 0 \\
0 & 0 & 0 & 0 & 0 & 1 \\
0 & 0 & 0 & 0 & 1 & 0 \\
0 & 0 & 0 & 1 & 0 & 0 \\
0 & 0 & 1 & 0 & 0 & 0 \\
0 & 1 & 0 & 0 & 0 & 0
\end{array}
\right) \ , \\
\  & \  & \  \\
M_{2}=\left(
\begin{array}{cccccc}
0 & 0 & 0 & 0 & 1 & 0 \\
0 & 0 & 0 & 1 & 0 & 0 \\
0 & 0 & 1 & 0 & 0 & 0 \\
0 & 1 & 0 & 0 & 0 & 0 \\
1 & 0 & 0 & 0 & 0 & 0 \\
0 & 0 & 0 & 0 & 0 & 1
\end{array}
\right) \ , & \  & M_{2}^{\prime }=\left(
\begin{array}{cccccc}
0 & 0 & 0 & 0 & 0 & 0 \\
0 & 0 & 0 & 0 & 1 & 0 \\
0 & 0 & 0 & 1 & 0 & 0 \\
0 & 0 & 1 & 0 & 0 & 0 \\
0 & 1 & 0 & 0 & 0 & 0 \\
0 & 0 & 0 & 0 & 0 & 0
\end{array}
\right) \ .
\end{array}
\end{equation*}
$Y_{1}$ and $Y_{2}$ have $q=(1:0:0:0:0:0)$ as single isolated point and
belong to a family of complete intersections which have been studied
extensively by Segre \cite{SEGRE1} and Kn\"{o}rrer \cite[pp. 38-51]{KNOERRER}%
. $\left( Y_{1},q\right) $ turns out to be an $\mathbf{A}_{5}$-singularity
and $\left( Y_{2},q\right) $ a $\mathbf{D}_{6}$-singularity. For both $Y_{1}$
and $Y_{2}$ the first summand in (\ref{ESTR-CI}) equals
\begin{equation*}
\left[ \binom{6}{3}-2\cdot 2\cdot \binom{6}{2}+3\cdot 2^{2}\cdot \binom{6}{1}%
-4\cdot 2^{3}\right] \ \left( 2^{2}\right) =0\ .
\end{equation*}
Hence, $e_{\text{str}}(Y_{1})=2+5-1=6\in \mathbb{Z},$ whereas
\begin{equation*}
\text{\ \ }e_{\text{str}}(Y_{2})=\tfrac{2633}{864}+6-1=8+\tfrac{41}{864}\in
\mathbb{QrZ\ }.
\end{equation*}
\bigskip
\end{examples}

\noindent {}\textsf{(b)} \textbf{Fiber products of elliptic surfaces over }$%
\mathbb{P}_{\mathbb{C}}^{1}$\textbf{. }Another kind of compact complex
threefolds having both $\mathbf{A}_{1}$ and $\mathbf{A}_{2}$-singularities
arises from a slight generalization of Schoen's construction \cite{SCHOEN2}.
Let $Z\rightarrow \mathbb{P}_{\mathbb{C}}^{1}$ and $Z^{\prime }\rightarrow
\mathbb{P}_{\mathbb{C}}^{1}$ denote two relatively minimal, rational
elliptic surfaces with global sections, and let $S$ (resp. $S^{\prime })$ be
the images of the exceptional fibers of $Y$ (resp. of $Y^{\prime }$) in $%
\mathbb{P}_{\mathbb{C}}^{1}.$ The fiber product
\begin{equation*}
Y:=Z\times _{\mathbb{P}_{\mathbb{C}}^{1}}Z^{\prime }\overset{\pi }{%
\longrightarrow }\mathbb{P}_{\mathbb{C}}^{1}
\end{equation*}
is a complex threefold with sigularities located only in the fibers
\begin{equation*}
Y_{s}=\pi ^{-1}\left( s\right) =Z_{s}\times Z_{s}^{\prime }
\end{equation*}
lying over points $s\in S^{\prime \prime }:=S\cap S^{\prime }.$ Since the
Euler number of any smooth fiber is zero, we have obviously
\begin{equation}
e\left( Y\right) =\sum_{s\in S^{\prime \prime }}\ e\left( Z_{s}\right) \
e\left( Z_{s}^{\prime }\right) .  \label{EULER-SCHOEN}
\end{equation}
We shall henceforth assume that $S^{\prime \prime }=\left\{
s_{1},s_{2},\ldots ,s_{\kappa }\right\} ,$ where for $1\leq i\leq \kappa ,\
Z_{s_{i}}$ is of Kodaira type $I_{b_{i}}$ (i.e., a rational curve with an
ordinary double point, if $b_{i}=1,$ and a cycle of $b_{i}$ smooth rational
curves, if $b_{i}\geq 2$), while $Z_{s_{j}}^{\prime }$ is of Kodaira type $%
I_{b_{j}^{\prime }}$, for all $j,$ $1\leq j\leq \nu ,$ $(\nu <\kappa <12),$
and of Kodaira type $II$ (i.e., a rational curve with one cusp), for all $j,$
$\nu +1\leq j\leq \kappa $. (See \cite[Thm. 6.2]{KODAIRA} for the
classification and Kodaira's notation of exceptional fibers). Under this
assumption, $Y$ is a $3$-dimensional Calabi-Yau variety with $%
b_{1}b_{1}^{\prime }+\cdots +b_{\nu }b_{\nu }^{\prime }$ $\mathbf{A}_{1}$%
-singularities (each of which contributing a $2$ as string-theoretic Euler
number) and $b_{\nu +1}+\cdots +b_{\kappa }$ $\mathbf{A}_{2}$-singularities
(each of which contributing a $\frac{7}{5}$ as string-theoretic Euler
number). Since $e\left( Z_{s_{i}}\right) =b_{i},$ for all $i,$ $1\leq i\leq
\kappa ,$ $e(Z_{s_{j}}^{\prime })=b_{j}^{\prime },$ for all $j,$ $1\leq
j\leq \nu ,$ and $e(Z_{s_{j}}^{\prime })=2,$ for all $j,$ $\nu +1\leq j\leq
\kappa $, the string-theoretic Euler number of $Y$ \ can be computed by (\ref
{ESTR-GLOBAL}) and (\ref{EULER-SCHOEN}), and can be written as follows:
\begin{equation}
e_{\text{str}}(Y)=2\left( \sum_{i=1}^{\nu }b_{i}b_{i}^{\prime }\right) +%
\frac{12}{5}\left( \sum_{i=\nu +1}^{\kappa }b_{i}\right)  \label{ESTR-SCHOEN}
\end{equation}

\begin{example}
Using Kodaira's homological and functional invariants (cf. \cite[\S8]
{KODAIRA}), as well as the normal forms of the corresponding Weierstrass
models (due to Kas \cite{KAS}), Herfurtner has shown in detail in \cite[cf.
Table 3, pp. 336-337 ]{HERFURTNER} the existence of relatively minimal,
rational elliptic surfaces $Z_{1}$ (resp. $Z_{2},Z_{3}$) with sections which
possess exactly four exceptional fibers having types $%
I_{1},I_{1},I_{5},I_{5} $ over the ordered $4$-tuple of points
\begin{equation*}
\left( \left( \tfrac{1+\sqrt{5}}{2}\right) ^{2},\left( \tfrac{1-\sqrt{5}}{2}%
\right) ^{2},0,\infty \right) \in (\mathbb{P}_{\mathbb{C}}^{1})^{4}
\end{equation*}
(resp. types $I_{1},I_{1},I_{2},I_{8}$ over $\left( -1,1,0,\infty \right)
\in (\mathbb{P}_{\mathbb{C}}^{1})^{4},$ resp. types $I_{1},I_{2},II,I_{7}$
over $(-\frac{9}{4},-\frac{8}{9},0,\infty )\in (\mathbb{P}_{\mathbb{C}%
}^{1})^{4}$). Hence,
\begin{equation*}
Y_{1}:=Z_{1}\times _{\mathbb{P}_{\mathbb{C}}^{1}}Z_{3},\text{ \ \ \ \ \ \ \
(resp. }Y_{2}:=Z_{2}\times _{\mathbb{P}_{\mathbb{C}}^{1}}Z_{3}\text{)},
\end{equation*}
has singularities only in the fibers over $0$ and $\infty ;$ more precisely,
it has five $\mathbf{A}_{2}$-singularities over $0$ and $35$ $\mathbf{A}_{1}$%
-singularities over $\infty $ (resp., two $\mathbf{A}_{2}$-singularities
over $0$ and $56$ $\mathbf{A}_{1}$-singularities over $\infty $).
Consequently, (\ref{ESTR-SCHOEN}) gives:
\begin{equation*}
e_{\text{str}}(Y_{1})=2\cdot 35+\frac{12}{5}\cdot 5=82\in \mathbb{Z}
\end{equation*}
whereas
\begin{equation*}
e_{\text{str}}(Y_{2})=2\cdot 56+\frac{12}{5}\cdot 2=116+\frac{4}{5}\in
\mathbb{QrZ}.
\end{equation*}
\end{example}

\bigskip\bigskip

\noindent \textbf{Acknowledgements. }We are indebted to E. Brieskorn for
pointing out the relation of \cite{ROCZEN1}, \cite{ROCZEN2} with the work of
P.J. Giblin, and to P.J. Giblin, who has sent to us a part of his thesis
\cite{GIBLIN} containing results on the canonical resolutions of
ADE-singularities.

We would like to express our warmest thanks to Nobuyuki Kakimi (University
of Tokyo) who informed us about some mistakes concerning our computations
for the discrepancy coefficients in a previous preprint-version of this
article.

We are also indebted to J. H. M. Steenbrink for drawing our attention to the
formula (\ref{SEB-T-ST}) which represents the analogue of Sebastiani-Thom
formula for the spectrum of the join of isolated hypersurface singularities.
Its use has simplified considerably our initial computations in \S 4.

\end{document}